\documentclass[11pt]{amsart}
\usepackage{}
\usepackage{amssymb}

\usepackage{amsmath}

\input xy
\xyoption{all}
\usepackage{soul}
\usepackage{amsfonts}
\usepackage{mathrsfs}

\usepackage{latexsym}
\usepackage{graphicx}
\usepackage{amscd,amssymb,amsmath,amsbsy,amsthm}
\usepackage[all]{xy}
\usepackage[colorlinks,plainpages,backref,urlcolor=blue]{hyperref}
\usepackage{verbatim}

\usepackage{tikz-cd}
\usepackage{tikz}
\usetikzlibrary{arrows,calc}
\usetikzlibrary{shapes,shadows,arrows}

\tikzset{
>=stealth',
help lines/.style={dashed, thick},
axis/.style={<->},
important line/.style={thick},
connection/.style={thick, dotted},
}

\usepackage{enumitem, hyperref}
\makeatletter
\def\namedlabel#1#2{\begingroup
    #2%
    \def\@currentlabel{#2}%
    \phantomsection\label{#1}\endgroup
}
\makeatother

\newcommand{\nc}{\newcommand}
\nc{\rnc}{\renewcommand}
\nc{\bb}[1]{{\mathbb #1}}
\nc{\bbA}{\bb{A}}\nc{\bbB}{\bb{B}}\nc{\bbC}{\bb{C}}\nc{\bbD}{\bb{D}}
\nc{\bbE}{\bb{E}}\nc{\bbF}{\bb{F}}\nc{\bbG}{\bb{G}}\nc{\bbH}{\bb{H}}
\nc{\bbI}{\bb{I}}\nc{\bbJ}{\bb{J}}\nc{\bbK}{\bb{K}}\nc{\bbL}{\bb{L}}
\nc{\bbM}{\bb{M}}\nc{\bbN}{\bb{N}}\nc{\bbO}{\bb{O}}\nc{\bbP}{\bb{P}}
\nc{\bbQ}{\bb{Q}}\nc{\bbR}{\bb{R}}\nc{\bbS}{\bb{S}}\nc{\bbT}{\bb{T}}
\nc{\bbU}{\bb{U}}\nc{\bbV}{\bb{V}}\nc{\bbW}{\bb{W}}\nc{\bbX}{\bb{X}}
\nc{\bbY}{\bb{Y}}\nc{\bbZ}{\bb{Z}}

\nc{\mbf}[1]{{\mathbf #1}}

\nc{\bfA}{\mbf{A}}\nc{\bfB}{\mbf{B}}\nc{\bfC}{\mbf{C}}\nc{\bfD}{\mbf{D}}
\nc{\bfE}{\mbf{E}}\nc{\bfF}{\mbf{F}}\nc{\bfG}{\mbf{G}}\nc{\bfH}{\mbf{H}}
\nc{\bfI}{\mbf{I}}\nc{\bfJ}{\mbf{J}}\nc{\bfK}{\mbf{K}}\nc{\bfL}{\mbf{L}}
\nc{\bfM}{\mbf{M}}\nc{\bfN}{\mbf{N}}\nc{\bfO}{\mbf{O}}\nc{\bfP}{\mbf{P}}
\nc{\bfQ}{\mbf{Q}}\nc{\bfR}{\mbf{R}}\nc{\bfS}{\mbf{S}}\nc{\bfT}{\mbf{T}}
\nc{\bfU}{\mbf{U}}\nc{\bfV}{\mbf{V}}\nc{\bfW}{\mbf{W}}\nc{\bfX}{\mbf{X}}
\nc{\bfY}{\mbf{Y}}\nc{\bfZ}{\mbf{Z}}
\nc{\bfa}{\mbf{a}}\nc{\bfb}{\mbf{b}}\nc{\bfc}{\mbf{c}}\nc{\bfd}{\mbf{d}}
\nc{\bfe}{\mbf{e}}\nc{\bff}{\mbf{f}}\nc{\bfg}{\mbf{g}}\nc{\bfh}{\mbf{h}}
\nc{\bfi}{\mbf{i}}\nc{\bfj}{\mbf{j}}\nc{\bfk}{\mbf{k}}\nc{\bfl}{\mbf{l}}
\nc{\bfm}{\mbf{m}}\nc{\bfn}{\mbf{n}}\nc{\bfo}{\mbf{o}}\nc{\bfp}{\mbf{p}}
\nc{\bfq}{\mbf{q}}\nc{\bfr}{\mbf{r}}\nc{\bfs}{\mbf{s}}\nc{\bft}{\mbf{t}}
\nc{\bfu}{\mbf{u}}\nc{\bfv}{\mbf{v}}\nc{\bfw}{\mbf{w}}\nc{\bfx}{\mbf{x}}
\nc{\bfy}{\mbf{y}}\nc{\bfz}{\mbf{z}}

\newcommand{\bS}{\mathbb{S}}
\newcommand{\G}{\mathbb{G}}

\nc{\mcal}[1]{{\mathcal #1}}
\nc{\calA}{\mcal{A}}\nc{\calB}{\mcal{B}}\nc{\calC}{\mcal{C}}\nc{\calD}{\mcal{D}}
\nc{\calE}{\mcal{E}} \nc{\calF}{\mcal{F}}\nc{\calG}{\mcal{G}}\nc{\calH}{\mcal{H}}
\nc{\calI}{\mcal{I}}\nc{\calJ}{\mcal{J}}\nc{\calK}{\mcal{K}}\nc{\calL}{\mcal{L}}
\nc{\calM}{\mcal{M}}\nc{\calN}{\mcal{N}}\nc{\calO}{\mcal{O}}\nc{\calP}{\mcal{P}}
\nc{\calQ}{\mcal{Q}}\nc{\calR}{\mcal{R}}\nc{\calS}{\mcal{S}}\nc{\calT}{\mcal{T}}
\nc{\calU}{\mcal{U}}\nc{\calV}{\mcal{V}}\nc{\calW}{\mcal{W}}\nc{\calX}{\mcal{X}}
\nc{\calY}{\mcal{Y}}\nc{\calZ}{\mcal{Z}}
\nc{\fA}{\frak{A}}\nc{\fB}{\frak{B}}\nc{\fC}{\frak{C}} \nc{\fD}{\frak{D}}
\nc{\fE}{\frak{E}}\nc{\fF}{\frak{F}}\nc{\fG}{\frak{G}}\nc{\fH}{\frak{H}}
\nc{\fI}{\frak{I}}\nc{\fJ}{\frak{J}}\nc{\fK}{\frak{K}}\nc{\fL}{\frak{L}}
\nc{\fM}{\frak{M}}\nc{\fN}{\frak{N}}\nc{\fO}{\frak{O}}\nc{\fP}{\frak{P}}
\nc{\fQ}{\frak{Q}}\nc{\fR}{\frak{R}}\nc{\fS}{\frak{S}}\nc{\fT}{\frak{T}}
\nc{\fU}{\frak{U}}\nc{\fV}{\frak{V}}\nc{\fW}{\frak{W}}\nc{\fX}{\frak{X}}
\nc{\fY}{\frak{Y}}\nc{\fZ}{\frak{Z}}
\nc{\fa}{\frak{a}}\nc{\fb}{\frak{b}}\nc{\fc}{\frak{c}} \nc{\fd}{\frak{d}}
\nc{\fe}{\frak{e}}\nc{\fFf}{\frak{f}}\nc{\fg}{\frak{g}}\nc{\fh}{\frak{h}}
\nc{\fri}{\frak{i}}\nc{\fj}{\frak{j}}\nc{\fk}{\frak{k}}\nc{\fl}{\frak{l}}
\nc{\fm}{\frak{m}}\nc{\fn}{\frak{n}}\nc{\fo}{\frak{o}}\nc{\fp}{\frak{p}}
\nc{\fq}{\frak{q}}\nc{\fr}{\frak{r}}\nc{\fs}{\frak{s}}\nc{\ft}{\frak{t}}
\nc{\fu}{\frak{u}}\nc{\fv}{\frak{v}}\nc{\fw}{\frak{w}}\nc{\fx}{\frak{x}}
\nc{\fy}{\frak{y}}\nc{\fz}{\frak{z}}

\theoremstyle{definition}

\newcommand{\tim}{\times}

\newcommand{\VV}{\mathcal{V}}
\newcommand{\A}{\mathbb{A}}
\newcommand{\Ao}{\mathbb{A}^n}
\newcommand{\B}{\mathbb{B}}

\newcommand{\N}{\mathbb{N}}
\newcommand{\Z}{\mathbb{Z}}

\newcommand{\R}{\mathbb{R}}
\newcommand{\Ss}{\mathbb{S}}
\newcommand{\C}{\mathbb{C}}

\nc{\hH}{{\mathbb H}}

\newcommand{\F}{\mathcal {F}}

\nc{\dddd}{d}
\let\d\dddd

\newcommand{\inj}{\hookrightarrow}

\newcommand{\g}{\mathfrak g}

\newcommand{\Ga}{\mathbb{G}_a}

\newcommand{\PP}{\mathbb{P}}

\DeclareMathOperator{\Hom}{Hom}
\DeclareMathOperator{\End}{End}

\DeclareMathOperator{\GL}{GL}

\DeclareMathOperator{\Gr}{Gr}

\DeclareMathOperator{\pt}{{pt}}

\DeclareMathOperator{\Spec}{Spec}

\DeclareMathOperator{\id}{id}

\DeclareMathOperator{\loc}{loc}

\DeclareMathOperator{\fac}{fac}
\DeclareMathOperator{\Rep}{Rep}

\newcommand{\Gm}{{\mathbb{G}_m}}

\newcommand {\Omit}[1]{}

\newcommand {\HH}{{\mathfrak H}}

\renewcommand {\Im}{\operatorname{Im}}

\newcommand{\rc}{\renewcommand}

\rc{\AA}{{\mcal A}}
\nc{\BB}{{\mcal B}} 
\nc{\CC}{{\mcal C}}
\nc{\DD}{{\mcal D}}
\nc{\EE}{{\mcal E}}
\nc{\FF}{{\mcal F}}
\nc{\GG}{{\mcal G}}
\nc{\II}{{\mcal I}}
\nc{\JJ}{{\mcal J}}
\nc{\KK}{{\mcal K}}
\nc{\LL}{{\mcal L}}
\nc{\MM}{{\mcal M}}
\nc{\NN}{{\mcal N}}
\nc{\OO}{{\mcal O}}
\nc{\QQ}{{\mcal Q}}
\nc{\RR}{{\mcal R}}
\rc{\SS}{{\mcal S}}
\nc{\TT}{{\mcal T}}
\nc{\UU}{{\mcal U}}
\nc{\WW}{{\mcal W}}
\nc{\ZZ}{{\mcal Z}}
\nc{\XX}{{\mcal X}}
\nc{\YY}{{\mcal Y}}

\nc{\sA}{{\mathsf A}}
\nc{\sB}{{\mathsf B}}
\nc{\sC}{{\mathsf C}}
\nc{\sD}{{\mathsf D}}
\nc{\sE}{{\mathsf E}}
\nc{\sF}{{\mathsf F}}
\nc{\sG}{{\mathsf G}}
\nc{\sH}{{\mathsf H}}
\nc{\sI}{{\mathsf I}}
\nc{\sJ}{{\mathsf J}}
\nc{\sK}{{\mathsf K}}
\nc{\sL}{{\mathsf L}}
\nc{\sM}{{\mathsf M}}
\nc{\sN}{{\mathsf N}}
\nc{\sP}{{\mathsf P}}
\nc{\sQ}{{\mathsf Q}}
\nc{\sR}{{\mathsf R}}
\nc{\sS}{{\mathsf S}}
\nc{\sT}{{\mathsf T}}
\nc{\sU}{{\mathsf U}}
\nc{\sV}{{\mathsf V}}
\nc{\sW}{{\mathsf W}}
\nc{\sX}{{\mathsf X}}
\nc{\sY}{{\mathsf Y}}
\nc{\sZ}{{\mathsf R}}
\nc{\sa}{{\mathsf a}}
\rc{\sb}{{\mathsf b}}
\rc{\sc}{{\mathsf c}}
\nc{\sd}{{\mathsf d}}
\nc{\ssf}{{\mathsf f}}
\nc{\sg}{{\mathsf g}}
\nc{\sh}{{\mathsf h}}
\nc{\sj}{{\mathsf j}}
\nc{\sk}{{\mathsf k}}
\nc{\sn}{{\mathsf n}}
\nc{\sr}{{\mathsf r}}
\nc{\su}{{\mathsf u}}
\nc{\sv}{{\mathsf v}}
\nc{\sw}{{\mathsf w}}
\nc{\sx}{{\mathsf x}}
\nc{\sy}{{\mathsf y}}
\nc{\sz}{{\mathsf z}}

\nc{\al}{{\alpha }}
\nc{\be}{{\beta }}
\nc{\ga}{{\gamma }}
\nc{\de}{{\delta }}
\nc{\ep}{{\varepsilon }}
\nc{\vap}{{\epsilon }}

\nc{\ze}{{\zeta }}
\nc{\et}{{\eta }}
\nc{\vth}{{\vartheta }}

\nc{\io}{{\iota }}
\nc{\ka}{{\kappa }}
\nc{\la}{{\lambda }}
\nc{\vpi}{{	\varpi		}}
\nc{\vrho}{{	\varrho		}}
\nc{\si}{{	\sigma 		}}
\nc{\ups}{{	\upsilon 	}}
\nc{\vphi}{{	\varphi 	}}
\nc{\om}{{	\omega 		}}

\nc{\Gam}{{\Gamma }}
\nc{\De}{{\Delta }}
\nc{\nab}{{\nabla}}
\nc{\na}{{\nabla}}
\nc{\Th}{{\Theta }}
\rc{\th}{{\theta }}
\nc{\La}{{\Lambda }}
\nc{\Ups}{{\Upsilon }}
\nc{\Om}{{\Omega }}

\rc{\Ao}{{	\A^1	}}
\nc{\Po}{{	\P^1	}}
\nc{\So}{{	S^1	}}

\nc{\All}{{	\forall		}}
\nc{\Exx}{{	\exists 	}}
\nc{\yy}{\infty}                       
\nc{\ys}{{  \frac{\infty}{2}  }}

\nc{\ii}{{i\in I}}
\nc{\nn}{{n\in \N}}
\nc{\ww}{{w\in W}}
\nc{\SES}[5]{{	0 @>>> {#1} @>{#2}>> {#3} @>{#4}>> {#5} @>>> 0	}}
\nc{\Ses}[3] {{	0 @>>> {#1} @>>>     {#2} @>>>     {#3} @>>> 0	}}
\nc{\pl}{{\oplus}}              		
\nc{\btim}{{\boxtimes}}
\nc{\ltim}{\ltimes}                  	%
\nc{\rtim}{\rtimes}			%
\nc{\ltr}{\triangleleft}        %
\nc{\rtr}{\triangleright}       %

\nc{\ten}{{	\otimes		}}            
\nc{\Lten}{{	\aa{L}\otimes	}}            
\nc{\Ltim}{{	\aa{L}\times	}}            
\nc{\Lcap}{{	\aa{L}\cap	}}            
\nc{\tenA}{	\bb{A}\ten	}
\nc{\tenB}{	\bb{B}\ten	}
\nc{\tenZ}{	\bb{\Z}\ten	}
\nc{\tenR}{	\bb{\R}\ten	}
\nc{\tenC}{	\bb{\C}\ten	}
\nc{\tenk}{	\bb{\k}\ten	}
\nc{\bten}{{\boxtimes}}         		

\nc{\con}{{ \CD@>>{\protect\cong}>\endCD }}  	
\nc{\conn}{{\CD    @<{\cong}<<  \endCD	}}  	
\nc{\Con}{{	\equiv		}}	
\nc{\appr}{{	\sim		}}	
\nc{\eqr}{{	\sim		}}	
\nc{\equi}{{	\sim		}}	

\nc{\fra}{ 	\frac	}     	
\nc{\ffr}[2]{{ 	\text{\footnotesize $\frac{#1}{#2}$	}	}}  

\nc{\ha}{{ \frac{1}{2} }}     		
	\nc{\half}{{ \frac{1}{2} }}    	

\nc{\ci}{{\circ}}               
\nc{\cd }{{\cdot}}            	
\nc{\cddd}{{\cdot\cdot\cdot}}	

\nc{\ox}{{	\OO_X		}}               
\nc{\omx}{{	\om_X		}}               
\nc{\Omx}{{	\Om_X^1		}}               
\nc{\Coh}{{	\CC oh		}}               %
\nc{\qcoh}{{	q\CC oh		}}               %

\nc{\xt}{{	X_*(T)		}}
\nc{\Xt}{{	X^*(T)		}}

\nc{\cfm}{{	co\fm		}}	

\nc{\tx}{	\text		}		
\nc{\df}{{ \protect\overset{ \text{def}}= 	}}		
\nc{\dff}{{ \ \df\				}}		
\nc{\inv}{{ {}^{-1}      }}			

\nc{\thh}{	^{\text{th}}	}                     	
\nc{\nd}{	^{\text{nd}}	}                     	
\nc{\rd}{	^{\text{rd}}	}                     	

\nc{\pmo}{{ 	\pm 1		}}
\nc{\mpo}{{ 	\mp 1		}}

\nc{\htt}{  \text{ht}}				

\nc{\emp}{{   \emptyset}}      			

\nc{\cowe}{{	\vee	}}			
\nc{\we}{{\wedge}}				
\nc{\wee}{{	\aa{\bullet}\wedge	}}		
\nc{\wetwo}{{     \pr\overset{2}\wedge       }}	

\nc{\limp}{{	\pr\underset {\leftarrow} \lim		}}	
\nc{\holimp}{{	\pr\underset {\leftarrow} {\tx{holim}}		}}	
\nc{\Limp}{{	\pr\underset {\leftarrow} {\bbb\lim}	}}	

\nc{\limi}{{	\pr\underset {\rightarrow}\lim		}}      
\nc{\holimi}{{
\pr\underset {\rightarrow}{\tx{holim}}	}}      
\nc{\Limi}{{	\pr\underset {\rightarrow}{\bbb\lim}	}}	

\nc{\llim}[1]{	 \bb{#1}\lim        	}   
\nc{\llimp}[1]{ \bb{#1}{ \pr\underset {\leftarrow} \lim       } }
\nc{\LLimp}[1]{ \bb{#1}{ \pr\underset {\leftarrow} {\bbb\lim} } } 	
\nc{\llimi}[1]{ \bb{#1}{ \pr\underset {\rightarrow}\lim       } }
\nc{\LLimi}[1]{ \bb{#1}{ \pr\underset {\rightarrow}{\bbb\lim} } }	
						
\nc{\ppp}{{ {\Bbb P}^1 }}            		
\nc{\ppn}{{ {\Bbb P}^n }}            		
				
\nc{\qlb}{{ \barr{{\Bbb Q}_l} }}      		
\nc{\ffq}{{  {\Bbb F}_q  }}           		
\nc{\ffp}{{  {\Bbb F}_p  }}           		
\nc{\tw}{   {}^{(1)}	}		

\nc{\AAb}{{ 	\AA b 		}}      		%
\nc{\Set}{{ 	\SS et 		}}      		%
\nc{\Top}{{ 	\TT op 		}}      		%
\nc{\cG}{{	\ch G		}}
\nc{\cB}{{	\ch B		}}
\nc{\cN}{{	\ch N		}}
\nc{\cT}{{	\ch T		}}
\nc{\cH}{{	\ch H		}}

\nc{\del}{{\partial }}
\nc{\delb}{{\barr\partial }}
\nc{\dd}[2]{	\fra{d{#1}}{d{#2}}		}
\nc{\ddel}[2]{	\fra{\del{#1}}{\del {#2}}	}
                                        
\nc{\hk}{{     \text{hyperk\"ahler} 	}}
\nc{\susy}{{\text{supersymmetry}}}

\nc{\ie}{{,\ \     \text{i.e.,}\ \ 	}}
\nc{\iif}{{\ \     \text{if}\ \ 	}}
\nc{\aand}{{\ \ \  \text{and}\ \ \ 	}}
\nc{\hence}{{\ \ \ \text{hence}\ \ \ 	}}
\nc{\while}{{\ \ \ \text{while}\ \ \ 	}}
\nc{\with}{{\ \ \  \text{with}\ \ \ 	}}
\nc{\oor}{{\ \     \text{or}\ \ 	}}
\nc{\foor}{{\ \     \text{for}\ \ 	}}
\nc{\suchthat}{{\ \     \text{such that}\ \ 	}}

\nc{\rk}{{\operatorname{rk}}}
\nc{\Ker}{{\operatorname{Ker}}}
\nc{\Coker}{{\operatorname{Coker}}}
\rc{\Im}{{ 	\text{Im} 	}}
\nc{\rank}{{	\ \text{rank} 	}}
\nc{\RHom}{{	\text{RHom}	}}
\nc{\HHom}{{	\text{$\HH$om}	}}
\nc{\RHHom}{{	\text{R$\HH$om} }}
\nc{\RGa}{{	\text{R$\Ga$}	}}
\nc{\EEnd}{{	\text{$\EE nd$}	}}
\nc{\AAut}{{	\text{$\AA ut$}	}}
\nc{\Ext}{{\operatorname{Ext}}}
\nc{\EExt}{{\operatorname{Ext}}}
\nc{\Tor}{{\operatorname{Tor}}}

\nc{\ord	}{{ \text{ord} }}			
\nc{\divv	}{{ \text{div} }}			

\nc{\timA} {{   \pr\underset{A}\tim             }}
\nc{\timB} {{   \pr\underset{B}\tim             }}
\nc{\timC} {{   \pr\underset{C}\tim             }}
\nc{\timG} {{   \pr\underset{G}\tim             }}
\nc{\timH} {{   \pr\underset{H}\tim             }}
\nc{\timN} {{   \pr\underset{N}\tim             }}
\nc{\timP}{{    \pr\underset{P}\tim             }}
\nc{\timQ}{{    \pr\underset{Q}\tim             }}
\nc{\timS} {{   \pr\underset{S}\tim             }}
\nc{\timT} {{   \pr\underset{T}\tim             }}
\nc{\timU} {{   \pr\underset{U}\tim             }}
\nc{\timV} {{   \pr\underset{V}\tim             }}
\nc{\timX} {{   \pr\underset{X}\tim             }}
\nc{\timY} {{   \pr\underset{Y}\tim             }}
\nc{\timZ} {{   \pr\underset{Z}\tim             }}

\nc{\ab}{{       ^{\text{ab}}   		}}
\nc{\af}{{       ^{\text{aff}}  		}}

\nc{\cod}{\text{codim}}	

\nc{\pr}{\protect}

\nc{\np}{{\newpage}}	

\nc{\lab}{	\label}
\nc{\npp}{{	\newpage\setcounter{page}{0}	}}
\nc{\setpa}{		\setcounter{part}		}
\nc{\setse}{		\setcounter{section}	}
\nc{\setsus}{		\setcounter{subsection}		}
\nc{\setsss}{		\setcounter{subsubsection}	}
\nc{\setpage}{		\setcounter{page}	}

\nc{\nfd}{ $$\text{ This version is preliminary and approximate, 
		             it is not for distribution. }$$	}

\nc{\noi}{{\noindent}}
\nc{\pf}{{	\noindent {\em Proof.}		}}
					
\nc{\epf}{ \fbox{\bf QED}	}

\nc{\heart}{{\tiny \cen{\tiny $\heartsuit $ }	}} 
\nc{\cont}{\tableofcontents}

\nc{\sbr}{{	\smallpagebreak	}}
\nc{\mbr}{{	\medpagebreak	}}
\nc{\bbr}{{	\bigpagebreak	}}
\nc{\bbb}{ 	\boldsymbol 	}

\nc{\bib}{		}
\nc{\bit}[1]{	\bibitem[#1]{#1} 		}

\rc{\b}{ 	\big         			}  
\nc{\lam}[1]{{ 	\text{\large $#1$	}	}}  
\nc{\smm}[1]{{ 	\text{\small $#1$	}	}}  
\nc{\fom}[1]{{ 	\text{\footnotesize $#1$	}	}}  
\nc{\tinm}[1]{{ \text{\tiny $#1$	}	}}  

\nc{\bu}{ \bullet         }  			
\nc{\bbu}{ \aa{\bbb \bullet}         }  	
\nc{\bus}{{	^\bullet	}}	 	
\nc{\bui}{{	_\bullet	}}	 	

\nc{\bem}{{	\begin{em}	}}
\nc{\eem}{{	\end{em} 	}}

\nc{\bbox}{{	\blackbox	}}	
\nc{\bx}{	\boxed	}		
\nc{\tbx}[1]{{\boxed{\tx{#1}}}}		

\nc{\mmbox}[1]{{	\mbox{$#1$}	}}	
\nc{\tbox}[1]{{		\mbox{\tx{#1}}	}}
 	
\nc{\ot}{		\leftarrow			}
\nc{\tto}{		\longrightarrow			}
\nc{\ott}{		\longleftarrow			}
\nc{\too}[1]{{		\aa{#1}\rightarrow			}}
\nc{\oot}[1]{{		\aa{#1}\leftarrow			}}
\nc{\ttoo}[1]{{		\aa{#1}\longrightarrow			}}
\nc{\oott}[1]{{		\aa{#1}\longleftarrow			}}

\nc{\Too}[2]{{		\aa{#1}{\bb{#2}\rightarrow}		}}
\nc{\TToo}[2]{{		\aa{#1}{\bb{#2}\longrightarrow}		}}
\nc{\ooTT}[2]{{		\aa{#1}{\bbb{#2}\longleftarrow}		}}
\nc{\toot}[2]{{		\aa{#1}{\bb{#2}\rightleftarrows}	}}
\nc{\ttoot}[2]{{	\aa{#1}{\bb{#2}\rightleftrightarrows}	}}

\nc{\ra}{{	\longrightarrow		}}

\nc{\raa}[1]{{	\aa{#1}\longrightarrow		}}
\nc{\laa}[1]{{	\aa{#1}\longleftarrow	}}	
\nc{\lra}{{\longrightarrow}}
\nc{\lla}{{\longleftarrow}}

\nc{\lr}{{\leftrightarrow}}     	
\nc{\lrs}{{\rightleftarrows}}     	

\nc{\imp}{{\Rightarrow}}        	
\nc{\impp}{{\Leftarrow}}        	
\nc{\eq}{{\Leftrightarrow}}        	
\nc{\impl}{{\Longrightarrow}}        	
\nc{\imppl}{{\Longleftarrow}}        	
\nc{\eql}{{\Longleftrightarrow}}        	
	\nc{\Ra}{{\Rightarrow}}         	
	\nc{\LRa}{{\Leftrightarrow}}        	

\nc{\injj}{{\pr	\hookleftarrow	}}    		

\nc{\sur}{{	\twoheadrightarrow	}}	
\nc{\surr}{{	\twoheadleftarrow	}}	
			
\nc{\mm}{{	\mapsto		}}     		
\nc{\mmm}{{	\leftarrow\shortmid }}		
\nc{\ainj}[1]{{\aa{#1}{\pr\hookrightarrow}	}}    	
\nc{\ainjj}[1]{{\aa{#1}{\pr\hookleftarrow}	}}    	

\nc{\asur}[1]{{	\aa{#1}\twoheadrightarrow	}}	
\nc{\asurr}[1]{{\aa{#1}\twoheadleftarrow	}}	
			
\nc{\amm}[1]{{	\aa{#1}\mapsto		}}     	
\nc{\ammm}[1]{{	\aa{#1}\leftarrow\shortmid }}	

\nc{\va}{{\uparrow}}              		

\nc{\syp}[1]{	^{ (#1) }		} 	
\nc{\up}[1]{	^{ (#1) }		} 	
\nc{\lp}[1]{	_{ (#1) }		}	
\nc{\hp}[1]{	^{ [#1] }		}	
\nc{\cle}{\preceq}		
\nc{\cl}{\prec}			
\nc{\cge}{\succeq}		
\nc{\cg}{\succ}			

\rc{\aa}{ 	\pr\overset 	}            

\nc{\indd}{{ ${} \ \ \ \ \  \ \        {} $	}}	
\nc{\inddd}{{ 	\indd\indd			}}	
\nc{\nnd}{{ 	\nn  \indd 			}}	
\nc{\nndb}{{ 	\nn  \indd $\bullet$		}}	
		
\nc{\bce}{	\begin{center}	}
\nc{\ece}{	\end  {center}	}
\nc{\cen}[1]{	\begin{center}	{\em  #1}	\end  {center}	}
\nc{\ce}[1]{	\begin{center}	{  #1}	\end  {center}		}

\nc{\bss}{{\backslash}}           		
\nc{\bs}{\bss}

\nc{\barr}{ 	\overline 	}      		
\nc{\ud}{	\underline	}		

\nc{\ti}{\tilde}              
\nc{\tii}{\widetilde}         

\nc{\hatt}{\widehat}				
\nc{\hata}{{	\bbb{ \hat{} }		}}	

\nc{\ch}{\check}              			
\nc{\cha}{{ 	\bbb{ \check{} }	}}      

\nc{\sub}{{	\subseteq	}}         
\nc{\subb}{{	\supseteq	}}         
\nc{\nsub}{{	\nsubseteq	}}         
\nc{\nsubb}{{	\nsupseteq	}}         %

\nc{\nin}{{	\notin	}}

\nc{\lb}{\langle}             				
\nc{\rb}{\rangle}
\rc{\l}{\langle}             				
\rc{\r}{\rangle}

\nc{\lB}{	\left(	}             			
\nc{\rB}{	\right)	}
\nc{\BBl}{{	\bbb{ \left( \right.}	}}             	
\nc{\BBr}{{	\bbb{ \left. \right)}	}}

\nc{\Pa}[2]{ {\lb} #1 {,} #2 {\rb} }				

\nc{\cD}[1]{ \tx{ $$\CD {#1} \endCD $$ }  }		

\nc{\sm} {		\left(		\smallmatrix	}	
\nc{\esm}{		\endsmallmatrix	\right)	}
\nc{\smat} {		\left(		\smallmatrix	}	
\nc{\esmat}{		\endsmallmatrix	\right)	}

\nc{\matr} {		\left[		\matrix	}	
\nc{\ematr}{		\endmatrix	\right]	}
\nc{\smr} {		\left[		\smallmatrix	}	
\nc{\esmr}{		\endsmallmatrix	\right]	}
\nc{\smatr} {		\left[		\smallmatrix	}	
\nc{\esmatr}{		\endsmallmatrix	\right]	}

\nc{\imat} {		\left.		\matrix	}	
\nc{\eimat}{		\endmatrix	\right.	}
\nc{\ism} {		\left.		\smallmatrix	}	
\nc{\eism}{		\endsmallmatrix	\right.	}

\nc{\ca}{		\left\{		\smallmatrix	}	
\nc{\eca}{		\endsmallmatrix	\right\}	}
\nc{\Ca}{		\left\{		\matrix		}	
\nc{\Eca}{		\endmatrix	\right.		}	
\nc{\eCa}{		\endmatrix	\right\}	}	

\nc{\com}{	\begin{diagram}	}
\nc{\ecom}{	  \end{diagram}	}

\nc{\tab}{	\begin{tabular}		}
\nc{\etab}{	\end{tabular}		}	

\nc{\Eq}{	\begin{equation}	}
\nc{\Eeq}{	\end{equation}	}
\nc{\aln}{	\begin{align}	}
\nc{\ealn}{	\end{align}	}

\nc{\pa}[1]{ 	\part{#1}		}
\nc{\pas}[1]{ 	\part*{#1}		}

\nc{\se}[1]{ 	\section{\bf#1}		}
\nc{\ses}[1]{ 	\section*{\bf#1}		}
\nc{\sus}{ 	\subsection		}
\nc{\sss}{ 	\subsubsection		}

\nc{\Lem}{ 	\subsection{Lemma}		}
\nc{\lem}{ 	\subsubsection{Lemma}		}
\nc{\slem}{ 	\subsubsection*{Lemma}		}
\nc{\sublem}{ 	\subsubsection{ Sublemma}	}
\nc{\ssublem}{ \subsubsection*{ Sublemma}	}

\nc{\Lemm}{ 	\subsection{Lemma}		}
\nc{\lemm}{ 	\subsubsection{Lemma}		}
\nc{\slemm}{ 	\subsubsection*{Lemma}		}
\nc{\sublemm}{ 	\subsubsection{ Sublemma}	}
\nc{\ssublemm}{ \subsubsection*{ Sublemma}	}

\nc{\Pro}{ 	\subsection{Proposition}	}
\nc{\pro}{ 	\subsubsection{Proposition}	}
\nc{\spro}{ 	\subsubsection*{Proposition}	}

\nc{\Cor}{ 	\subsection{Corollary}		}
\nc{\cor}{ 	\subsubsection{Corollary}	}
\nc{\scor}{ 	\subsubsection*{Corollary}	}

\nc{\Corr}{ 	\subsection{Corollary}		}
\nc{\corr}{ 	\subsubsection{Corollary}	}
\nc{\scorr}{ 	\subsubsection*{Corollary}	}

\nc{\Theo}{ 	\subsection{Theorem}		}		
\nc{\theo}{ 	\subsubsection{Theorem}		}
\nc{\stheo}{ 	\subsubsection*{Theorem}	}

\nc{\pretheo}{ 	\subsubsection{Pretheorem}	}

\nc{\rem}{ 	\subsubsection{Remark}		}
\nc{\srem}{ 	\subsubsection*{Remark}	}

\nc{\rems}{ 	\subsubsection{Remarks}		}

\nc{\srems}{ 	\subsubsection*{Remarks}	}

\nc{\Def}{ 	\subsection{Definition}		}
\nc{\ddef}{ 	\subsubsection{Definition}	}

\nc{\comm}{ 	\subsubsection{Comment}		}
\nc{\scomm}{ 	\subsubsection*{Comment}	}
\nc{\comms}{ 	\subsubsection{Comments}		}
\nc{\scomms}{ 	\subsubsection*{Comments}	}

\nc{\sclaim}{ 	\subsubsection*{Claim}	}

\nc{\nota}{ 	\subsubsection{Notation}	}

\nc{\sconj}{ 	\subsubsection*{Conjecture}	}

\nc{\ex}{ 	\subsubsection{Example}		}
\nc{\sex}{ 	\subsubsection*{Example}	}
\nc{\exs}{ 	\subsubsection{Examples}	}
\nc{\sexs}{ 	\subsubsection*{Examples}	}
\nc{\Ex}{ 	\subsection{Example}		}
\nc{\sEx}{ 	\subsection*{Example}	}
\nc{\Exs}{ 	\subsection{Examples}	}
\nc{\sExs}{ 	\subsection*{Examples}	}

\nc{\que}{ 	\subsubsection{Question}	}
\nc{\ques}{ 	\subsubsection{Questions}	}
\nc{\sque}{ 	\subsubsection*{Question}	}
\nc{\sques}{ 	\subsubsection*{Questions}	}

\nc{\bi}{	\begin{itemize}\item		}
\rc{\i}{	\item			}
\nc{\ei}{ \end{itemize}	} 
\nc{\ben}{	\begin{enumerate}\item		}
\nc{\een}{	\end{enumerate}			}
\nc{\ftt}[1]{{\footnote{#1}}}
\nc{\fttt}[1]{{$^($\footnote{#1}$^)$}}
\nc{\bftt}[1]{\footnote{#1}}

\nc{\Ff}[2]{ \fbox{#1$ $}\footnote{ \fbox{!}#2 }\fbox{$ $}		}
\nc{\f}[1]{ \fbox{$ $}\footnote{ \fbox{!}#1 }\fbox{$ $}		}
\nc{\ff}[1]{ \fbox{For}\footnote{ \fbox{Remove?}#1 }\fbox{get}		}

\nc{\Lb}{\LL_bullet}
\nc{\Lbb}{\LL_{\bu\bullet}}
\let\LL\calL

\nc{\tT}{{\tii T}}

\nc{\kk}{\Bbbk}
\rc{\k}{\kk}

\rc{\P}{\bb{P}}

\nc{\vb}{{v^\bu}}
\rc{\vb}{{\bar v}}
\nc{\vbb}{{v_\bu}}

\nc{\RRep}{{\RR ep}}

\nc{\bv}{{\bf v}}
\nc{\bbi}{{i_\bu}}
\nc{\bbv}{\barr{\bv}}

\nc{\bz}{{\bf 0}}
\nc{\bm}{{\bf m}}

\nc{\bV}{\barr V}
\nc{\bP}{{\barr P}}
\nc{\bp}{{\barr \fp}}
\nc{\Qs}{{Q^*}}
\nc{\bQ}{{\barr Q}}

\nc{\dop}{{ \aa{\cd}p }}
\nc{\doq}{{\aa{\cd}q}}

\nc{\sq}{{	\sqcup		}}

\nc{\txr}[1]{\textcolor{red}{#1}}

\nc{\tG}{{\tii G}}
\nc{\tP}{{\tii P}}
\nc{\tL}{{\tii L}}

\nc{\edd}{\end{document}}

\nc{\all}{{ ^{(\alpha)} }}
\nc{\bee}{{ ^{(\beta)} }}
\nc{\gaa}{{ ^{(\gamma)} }}
\nc{\nnnn}{ 		\syp{n}		}     
\nc{\nnn}{               \hp{n}	 	}     
\nc{\GK}{{  	G_\KK		}}
\nc{\GO}{{  	G_\OO		}}
\nc{\gK}{{  	\fg_\KK		}}
\nc{\gO}{{  	\fg_\OO		}}
\nc{\go}{{  	\fg_\OO		}}

\nc{\BK}{{  	B_\KK		}}
\nc{\BO}{{  	B_\OO		}}
\nc{\bK}{{  	\fg_\KK		}}
\nc{\bO}{{  	\fb_\OO		}}

\nc{\NK}{{  	N_\KK		}}
\nc{\NO}{{  	N_\OO		}}
\nc{\nK}{{  	\fn_\KK		}}
\nc{\nO}{{  	\fn_\OO		}}

\nc{\TK}{{  	T_\KK		}}
\nc{\TO}{{  	T_\OO		}}
\nc{\tk}{{  	\ft_\KK		}}
\nc{\tO}{{  	\ft_\OO		}}

\nc{\LK}{{  	L_\KK		}}
\nc{\LO}{{  	L_\OO		}}

\nc{\PK}{{  	P_\KK		}}
\nc{\PO}{{  	P_\OO		}}

\nc{\Kh}{\tx{Ka\"hler\ }}
\nc{\Khs}{\tx{Ka\"hler structure\ }}
\nc{\Khss}{\tx{Ka\"hler structures\ }}

\nc{\GKh}{\tx{Generalized Ka\"hler\ }}
\nc{\GKs}{\tx{Generalized Ka\"hler structure\ }}
\nc{\GKss}{\tx{Generalized Ka\"hler structures\ }}
\nc{\gKs}{\tx{Generalized Ka\"hler structure\ }}
\nc{\gKss}{\tx{Generalized Ka\"hler structures\ }}

\nc{\sYM}{\text{super Young-Mills\ }}
\nc{\tFT}{\text{topological Field Theory\ }}
\rc{\top}{\tx{topological\ }}
\rc{\Top}{\tx{Topological\ }}
\nc{\TFT}{\text{Topological Field Theory\ }}
\nc{\TQFT}{\text{Topological Quantum Field Theory\ }}
\nc{\TQFTs}{\text{Topological Quantum Field Theories\ }}
\nc{\QFT}{\text{Quantum Field Theory\ }}
\nc{\QFTs}{\text{Quantum Field Theories\ }}
\nc{\FT}{\text{Field Theory\ }}

\nc{\HM}{\text{Hitchin moduli\ }}
\nc{\Hf}{\text{Hitchin fibration\ }}
\nc{\Wi}{\text{Wilson\ }}
\nc{\Wo}{\text{Wilson operator\ }}
\nc{\Wos}{\text{Wilson operators\ }}
\nc{\tH}{\text{t'Hooft\ }}
\nc{\tHo}{\text{t'Hooft operator\ }}
\nc{\Ho}{\text{t'Hooft operator\ }}
\nc{\tHos}{\text{t'Hooft operators\ }}
\nc{\Hos}{\text{t'Hooft operators\ }}
\nc{\Sd}{\text{S-duality\ }}
\nc{\Ld}{\text{Langlands duality\ }}
\nc{\Be}{\text{Bogomolny equations\ }}
\nc{\He}{\text{Hecke\ }}
\nc{\Heo}{\text{Hecke operators\ }}
\nc{\Hem}{\text{Hecke modifications\ }}
\nc{\Bgs}{\tx{Bogomolny equations\ }}
\nc{\Bg}{\tx{Bogomolny equation\ }}
\nc{\Hk}{{\text{Hyperk$\ddot{a}$hler} }}

\nc{\eqq}{{ 	\ =\			}}
\nc{\Cy}{{ 	C_\yy	}}
\nc{\Ay}{{ 	A_\yy	}}
\nc{\Ly}{{ 	L_\yy	}}
\nc{\Cm}{{ 	C_m	}}
\nc{\Am}{{ 	A_m	}}
\nc{\Lm}{{ 	L_m	}}

\nc{\En}{{\tx{E$_n$}}}
\nc{\Eo}{{\tx{E$_1$}}}
\nc{\Et}{{\tx{E$_2$}}}
\nc{\Ey}{{ \tx{	E$_\yy$}	}}

\rc{\Cy}{{	Cat_\yy			}}
\nc{\Cey}{{	Cat^{ex}_\yy		}}
\nc{\Cyt}{{	\tii{Cat}_\yy		}}

\nc{\bGm}{{	\barr{G_m}			}}

\nc{\Uo}{{	U(1)	}}
\nc{\sqt}{{	\sqrt{2}	}}

\nc{\BGB}{{B\bss G/B}}

\nc{\Glb}{{ \barr{\GG_\la}	}}

\nc{\mut}{{	\mu_2	}}

\nc{\PPP}{\calP}
\nc{\HHH}{\calH}

\nc{\bii}{{ \tx{\bf i} }}

\nc{\sett}[1]{\setcounter{tocdepth}{#1}}

\nc{\timL}{\underset{L}\tim}
\nc{\GtimP}{G\underset{P}\tim}
\nc{\GtimL}{G\underset{L}\tim}

\nc{\RV}{\Rep(V)}
\nc{\RbV}{\Rep(\bV)}
\nc{\RF}{\Rep(F)}

\nc{\bchi}{\barr\chi}
\nc{\tchi}{\tii\chi}

\nc{\FS}{{\FF\SS}}

\nc{\GI}{{\G\tim \tx{I}}}
\nc{\VVI}{{\VV^\tx{I}}}

\nc{\HGI}{\HHH_{\bbG \times I}}

\nc{\GGG}{{	\GG(G)		}}
\nc{\GGGm}{{	\GG(\Gm)	}}
\nc{\tGGG}{{	\tii\GG(G)	}}

\nc{\GmO}{{		G_{m,\OO}		}}
\nc{\GmOn}{{		G_{m,\OO_n}		}}
\nc{\GmK}{{		\mathbb{G}_{m,\KK}		}}
\nc{\GmKr}{		G_{m,\KK}^{red}		}
\nc{\GmKu}{{		\protect\ud{G_{m,\KK}}		}}

\nc{\Gmd}{{	G_{m,\d}	}}
\nc{\Gmdu}{{	G_{m,\du}  }}
\nc{\Gmds}{{	G_{m,\ds}			}}
\nc{\Gmdsu}{{	G_{m,\dsu}  }}

\nc{\uQ}{\ud{Q}}

\nc{\tCC}{{\tii\CC}}
\nc{\CI}{{{C\tim I}}}
\nc{\CdI}{{ \CC_{\d\tim I} }}       \nc{\CCdI}{{ \CC_{\d\tim I} }}
\nc{\dI}{{{\d\tim I}}}
\nc{\HHCI}{{\HHH_{C\tim I}}}
\nc{\HHdI}{{\HHH_{\d\tim I}}}
\nc{\HCI}{{\HHH_{C\tim I}}}

\nc{\bc}[2]{{ \smat {#1}\\ {#2}\esmat }}

\nc{\sta}{\star}

\rc{\HH}{\calH}
\nc{\HCal}{{\HH_C^\al}}
\nc{\HCbe}{{\HH_C^\be}}

\nc{\Vn}{{V_\bu}}

\nc{\bRR}{{\barr\RR}}
\nc{\RRb}{{\barr\RR}}
\nc{\Vb}{{V_\bu}}
\nc{\FR}{{\FF^m\RR}}

\nc{\GXo}{{G_1\bs X_1}}
\nc{\GXt}{{G_2\bs X_2}}
\nc{\GXi}{{G_i\bs X_i}}

\nc{\RRs}{{\RR_\Qs}}
\rc{\RRb}{{\RR_\bQ}}

\nc{\bl}{{\bbb\fl}}

\nc{\BQ}{\bold Q}

\nc{\cdF}{{(d^*F)^{o}}}

\nc{\cp}{\aa{o}p}
\nc{\cq}{\aa{o}q}

\nc{\NI}{{\N[I]}}

\nc{\Gs}{{\G^2}}
\nc{\Gsc}{{G_{sc}}}
\nc{\Tsc}{{T_{sc}}}
\nc{\HGIc}{{\HH_\bbGI^c}}

\nc{\Vd}{{	\aa{\bbb .}V	}}

\nc{\glpgrass}{{\GG^P(I,Q)	}}

\nc{\uGG}{{ \ud{\GG} }}
\rc{\GGG}{{ \GG(G) }}
\nc{\uGGG}{{ \uGG(G) }}

\nc{\HHXI}{{	       \HH_{X\tim I}		}}
\nc{\HHCIal}{{	       \HH_{C\tim I}^\al		}}
\nc{\HHC}{{	       \HH_C	  		}}
\nc{\CCC}{{	       \CC_C	  		}}
\nc{\CCCI}{{	       \CC_{C\tim I}	  		}}
\nc{\HCIv}[1]{{	       \HH_{C\tim I}^{#1}		}}
\nc{\HC}{{	       \HH_C	  		}}
\nc{\HHCr}{{	       \HH_C^{reg}  		}}

\nc{\RanC}{{Ran(C)}}
\nc{\GGDe}{{\GG^{Gai}_{Ran(C)}}}
\nc{\HHle}{{ \HH_\CI^{\le}    }}

\nc{\GGRan}{\GG_{Ran(C)}}

\rc{\loc}{^{loc}}

\rc{\d}{{      d		      }}		
\nc{\ds}{{     d^*		      }}		
\nc{\dt}{{     d^2		      }}

\nc{\uX}{{ \protect	\ud{\fX}		}}
\nc{\fXu}{{	\ud\fX				}}
\nc{\ufX}{{	\ud\fX				}}
\rc{\sX}{{	\fX^*				}}
\nc{\suX}{{	\ud{\fX^*}			}}
\nc{\sXu}{{	\ud{\fX^*}			}}
\nc{\du}{{      \protect\ud{\d}			}}
\nc{\uds}{{	\protect\ud{\ds}			}}
\nc{\dsu}{{	\protect\ud{\ds}			}}

\nc{\dud}{{	\d_{\du}			}}

\nc{\Galp}{{		G_\al			}}

\nc{\uP}{{  	  	\ud{P}		}}
\nc{\buP}{{ 		\uP{}^\bu	}}
\nc{\bL}{{  		L{}^\bu		}}

\nc{\Aoz}{{	\A^1_z		}}
\nc{\Cx}{{	C_x		}}
\nc{\tTT}{{	\tii\TT		}}
\nc{\tGG}{{	\tii\GG		}}

\rc{\GGG}{{	\GG(G)		}}

\nc{\tGm}{{ \tii {\Gm}		}}
\nc{\tGmb}{{ \tii {\Gm}^b 	}}

\nc{\MMap}{{	\MM ap			}}
\nc{\RMap}{{	\tx{RMap}		}}
\nc{\RMMap}{{	\tx{R}\MM\ap		}}

\nc{\DDD}{{   \bbb\sD  }}

\nc{\Gk}{{	G_\ka	}}
\nc{\Aff}{{    \tx{Aff}		  }}
\nc{\tGL}{{	\tii{GL}	}}
\nc{\BGm}{{	     B(\Gm)  }}

\nc{\Kx}{{     \KK_x	      }}
\nc{\Kxn}{{     \KK_x^n	      }}

\nc{\GGt}{{  \GG^2	}}
\nc{\BGmg}{{ B(\Gm)^{gr}	}}

\nc{\detg}{{	 \det^{gr}	}}

\nc{\BunnxX}{{	Bun_{n,x} (X)		}}
\nc{\tGKx}{{	\tii\GK_x 		}}

\nc{\CtI}{{	       C\tim I			}}
\nc{\HHCtI}{{	       \HH_{C\tim I}		}}

\nc{\Xr}{{	\ZZ_{reg}  		}}
\nc{\Xm}{{	\ZZ\syp{m}	}}
\nc{\Xmr}{{	\ZZ\syp{m}_{reg}	}}
\nc{\XmD}{{	\ZZ\syp{m}_D	}}

\nc{\Ym}{{	\YY\syp{m}		}}
\nc{\Ymr}{{	\YY\syp{m}_{reg}	}}
\nc{\YmD}{{	\YY\syp{m}_D		}}

\nc{\PD}{{     \PP\sD	}}

\nc{\LS}{{ \tx{LSp}	}}
\nc{\cii}{{  \ch i }}

\nc{\fSch}{{	 \tx{fSch}		}}
\nc{\ffSch}{{	 \tx{f}^*\tx{Sch}	}}

\rc{\sh}{{	 \#		}}

\nc{\uDD}{{	\ud{\DD}	}}
\nc{\Cz}{{	C^o	}}

\nc{\bT}{{	\bbb{T}	}}
\nc{\bSz}{{	\barr{S_0}	}}

\nc{\bSSz}{{	\barr{\SS_0}	}}	
\nc{\bSal}{{	\barr{S_\al}	}}	
\nc{\bSla}{{	\barr{S_\la}	}}	
\nc{\bSd}{{	\barr{S_d}	}}	
\rc{\bS}[1]{{	\barr{S_{#1}}	}}	

\nc{\bTz}{{	\barr{T_0}	}}	
\rc{\bT}[1]{{	\barr{T_{#1}}	}}	
\nc{\bTal}{{	\barr{T_{\al}}	}}

\nc{\bSmz}{{	\barr{S^-_0}	}}

\nc{\bSp}{{	\barr{S_{p\cai} } 	}}

\rc{\com}{{	\ch\om		}}
\nc{\como}{{	\ch\om_1		}}
\nc{\cep}{{	\ch\ep		}}

\nc{\Czd}{{	\C[z]_{<d}	}}

\nc{\hpT}{{	^T			}}
\nc{\OGo}{{	\OO_\GG(1)		}}

\nc{\cQ}{\ch Q}

\nc{\GOm}{{G_{\OO_-} }}

\nc{\AJ}{\tx{AJ}}

\nc{\GNa}{{	(G/N)\protect\af	}}

\nc{\Fq}{{\F_q}}
\nc{\HHd}{\HH_{d} }

\rc{\Ga}{\Gam}

\nc{\matt}[4]{\left(\begin{\matrix}	{#1} & {#2} \\ {#3} & {#4} \end{matrix}\right)  }
\nc{\mat}{\begin{\bmatrix}}	
\nc{\emat}{\end{\bmatrix}}

\nc{\GmI}{\mathbb{G}_m^I}
\nc{\Gmt}{\mathbb{G}_m^2}
	
\nc{\CCI}{{\CC_I}}
	
\nc{\BD}{\tx{BD}}

\author[I.~Mirkovi\'c]{Ivan Mirkovi\'c}
\address{University of Massachusetts, Department of Mathematics and Statistics, 710 N. Pleasant St., Amherst MA 01003}
\email{mirkovic@math.umass.edu}

\author[Y.~Yang]{Yaping~Yang}
\address{School of Mathematics and Statistics, The University of Melbourne, 813 Swanston Street, Parkville VIC 3010, Australia}
\email{yaping.yang1@unimelb.edu.au}

\author[G.~Zhao]{Gufang~Zhao}

\address{Institute of Science and Technology Austria,
Am Campus, 1,
Klosterneuburg 3400,
\"Osterreich}\curraddr{School of Mathematics and Statistics, The University of Melbourne, 813 Swanston Street, Parkville VIC 3010, Australia}
\email{gufangz@unimelb.edu.au}

\nc{\aj}{{ \sA\sJ }}
\nc{\cth}{{\ch\theta}}

\rc{\bP}{{ \bold P }}

\nc{\bZ}{{\bold Z}}

\nc{\GGbd}{\GG^{BD}}
\nc{\GGBD}{\GG^{BD}}

\nc{\eCD}{\endCD}

\begin{document}

\title{
Loop Grassmannians of quivers and affine quantum groups
}

\date{\today}
\maketitle

\begin{flushright}
{\em To Alexander Beilinson and Victor Ginzburg}
\end{flushright}

\sett{2}

\begin{abstract}
We construct for each choice of a quiver $Q$, a cohomology theory $A$
and a poset $P$ a ``loop Grassmannian'' $\GG^P(Q,A)$.
This generalizes loop Grassmannians of semisimple groups
and the loop Grassmannians of based quadratic forms.
The addition of a 
```dilation'' torus
$\DD\sub \Gmt$ gives a
quantization
$\GG^P_\DD(Q,A)$.
This
construction is motivated by the program of introducing
an inner  cohomology theory in algebraic geometry adequate for the Geometric Langlands program
\cite{Mi} and  
on the construction of affine quantum groups from
generalized cohomology theories \cite{YZell}.
\end{abstract}

\tableofcontents

I.M. is very happy for opportunity to mention just a few transformative
effects of personalities of Sasha Beilinson and Vitya Ginzburg.
I.M.'s understanding of possibilities of being a mathematician
have been upturned through
Bernstein's talk at Park City
and
Beilinson's talks in Boston.
A part of the magic was that mathematics was alive,
high on ideas, low on ownership
and each talk would open 
for thinking some topic in mathematics, almost regardless of one's preparation. 
Before meeting Ginzburg, I.M. has  come to view him as a smarter twin brother
in mathematical tastes. Of  biggest  influence on I.M. was Ginzburg's paper
on loop Grassmannians that offered a new kind of mathematics,
orchestrated by an explosion of geometric ideas.

\setse{-1}
\se{
Introduction
}
For a semisimple algebraic group $G$ of ADE type,
the corresponding quiver
$Q$ is used to  study  representations
of $G$, its  loop group  $G((t))$
and their quantum versions.
Here we reconstruct 
from $Q$ the loop Grassmannian $\GG(G)$ of $G$. 
This construction produces a candidate for the loop Grassmannian
Kac-Moody extensions $\fg\af$ of  loop Lie algebras $\fg((t))$
and should provide a tool to study its representations.

\sss{
Loop Grassmannians associated to quivers
}
An advantage of the quiver approach is that it works
in large generality. It provides  a  ``loop Grassmannian''
$\GG_\DD^P(Q,A)$ associated to the data of an
arbitrary quiver $Q$, a cohomology theory $A$, a poset $P$ and a torus $\DD$ of dilations.
Intuitively, a quiver $Q$ should provide a ``grouplike'' object $G(Q,A)$
though at the moment we only see objects that should correspond to
(quantization of) its affinization.

A  cohomology theory $A$ gives a ``{\em cohomological schematization}'' functor $\fA(X)\dff
\Spec[A(X)]$ which assigns to a space $X$ the affine scheme $\fA(X)$
over the ring of constants of theory $A$.\fttt{
When $A$ is de Rham cohomology then $\fA_X$ can be viewed as an affinization
of the de Rham space $X_{dR}$ of $X$.
}
For us this simplifies
stacks (spaces with much  symmetry)
to classical geometry.
It takes the moduli of lines\ie the classifying space $\B\Gm$,
to a curve $\G=\fA(\B\Gm)$.   
Next, it turns
the moduli $\VV$ of finite dimensional  vector spaces
into the space of configurations
on the curve $\G$\ie the Hilbert scheme of points $\HH_\G=\sq_n\ \G\syp{n}$
of $\G$.

This configuration space  is then used as
the setting for the Beilinson-Drinfeld version
of the ``loop Grassmannian'' $\GG(Q,A)$
that intuitively corresponds to a (yet undefined) group
$G(Q,A)$.

Finally, we construct the space 
$\GG_\DD(Q,A)$ which should be the 
quantum loop Grassmannian of the  group 
$G(Q,A)$.
Here, one adds quantization  by letting a torus $\DD$
act on 
the extension correspondence for
representations of quivers
(and its cotangent stack). 
At this level there is a well defined related ``group theoretic'' object,
the
``affine quantum group'' that was constructed in \cite{YZ2} and  is denoted here by
$U_{\DD}(Q, A)$.

\sss{
Construction
}
 Since we avoid the group $G(Q,A)$ and its affinization, the construction
is  less standard.
We will argue that it is of  ``homological nature''.

The geometric ingredient is the technique of
{\em local projective spaces} 
from \cite{Mi}. This refers
to the notion of  {\em $I$-colored local vector bundles}
over a curve $C$\ie   
vector bundles over the $I$-colored  configuration space $\HHCI$
(the moduli of finite flat subschemes of  $C\tim I$),
that are in a certain sense ``local with respect to $C$''.

We actually start with a local line bundle $\sL$ over $\HHCI$ and we
induce it using a finite poset $P$ to a
local sheaf $Ind^P(\sL)$ over $\HHCI$.
The corresponding projective scheme $\bP(Ind^P\sL)$
is
 the projective spectrum of it symmetric algebra.
The (``projective'') {\em zastava space} $Z^P(\sL)$
is its ``local part'' $\bP\loc[Ind^P\sL]$.
Its fibers are obtained as collisions of products of fibers
at colored points $ai\in \HHCI$  (for
a point $a\in C$ and a color $i\in I$).
The collision happens inside
the projective scheme $\bP[Ind^P(\sL)]$ and the  ``rules of collision'' are
specified by the locality structure
on the line bundle $\sL$.

The local line bundles $\sL$ on $\HHCI$ are classified by 
symmetric bilinear forms $\QQ$ on $\Z[I]$, then the zastava space is denoted
$Z^P(I,\QQ)$.
However, we will concentrate on the case when $\QQ$ is associated to a quiver $Q$
(with vertices $I$).
In this case the topological ingredient of our construction
is an explicit reconstruction of
$\sL$ 
from the quiver $Q$ as
the Thom line bundle
of the moduli of representations
of a quiver $Q$.
This topological construction allows adding to the data a choice 
of a cohomology theory $A$,
to get the zastava space
$Z^P(Q,A)$.
Moreover we can upgrade to 
the {\em quantum version} $Z^P_\DD(Q,A)$ by
replacing the moduli of representations with the 
cotangent stack of the moduli of extensions of such representations, and
by switching
on an  action of a certain small torus $\DD$
from \cite{YZ2} on this cotangent stack.

Finally, we get the loop Grassmannian $\GG^P_\DD(Q,A)$
as a certain union of fibers of these zastava spaces.

\srem
The classical loop Grassmannians
of reductive groups are recovered when
the poset $P$ is a point (then $P$ is omitted  from notation).
In that case the fiber of
$Z^P(\sL)$
at any colored point is $\Po$.\fttt{
Here, we do not   pay attention to a choice of $P$
but when
$P=(1<\cddd<m)$ the fibers at
colored points are  $\P^m$ and 
$\GG^{1<\cddd<m}_{\DD}(Q, A)$
should ``correspond to level $m$'' in the sense that
the sections of the standard line bundle $\OO(1)$ on this object should be the same as the
sections of $\OO(m)$ in the case when $P=\pt$.
}

\sss{
``Quantum nature'' of $\GG^P_\DD(Q,A)$
}
The
loop Grassmannian $\GG(G)$  of a semisimple  group $G$ is a 
partial flag variety of $G\af$ so it has a known quantum version which is a 
non-commutative geometric object. For the $\GG_\DD(Q,A)$
construction this corresponds to 
the case when $A$ is the $K$-theory.
However, our incarnation $\GG_\DD(Q,A)$
is an object in standard geometry, and the hidden noncommutativity
manifests in its Beilinson-Drinfeld
form\ie when $\GG_\DD(Q,A)$ is extended to lie over a
configuration space. The configuration space is necessarily
ordered (``non-commutative'')\ie $\HH^n_{C\tim I}=(C\tim I)\syp{n}$
is replaced by
$(C\tim I)^n$. This has more connected components  
but this increase is ameliorated 
by a non-standard  
feature, a meromorphic braiding relating different connected components of the
configuration moduli.

We expect to have more explicit
descriptions of $\GG^P_\DD(Q,A)$ 
in terms of 
the graded algebra of sections of line bundles $\OO(m)$ or
in terms of the equation for the embedding into the projective space
corresponding to sections of $\OO(1)$.\fttt{
These embedding equations should be integrable hierarchies of differential equations
indexed by  $Q,P$ and $A$ since this is true in the classical case of $\GG(G)$.
}
In this paper we only do some preparatory steps
towards 
identifying the cases of $\GG^P_\DD(Q,A)$
with the 
classical loop Grassmannian of reductive groups.


This paper is   related to the work of Z. Dong
\cite{D18} that studies the relation between the Mirkovi\'c-Vilonen cycles in loop Grassmannians
and 
the quiver Grassmannian of representations of the preprojective algebra
(see \ref{Fixed points in MV cycles}). 


\sss{
Contents
}
In section \ref{Recollections on cohomology theories}
we recall the method of generalized cohomology theories.
Section
\ref{Loop Grassmannians and local spaces}
covers relevant aspects of classical loop Grassmannians
and how to  rebuild and generalize these in a ``homological'' way\ie
by turning the notion of  locality
into a construction.
In section \ref{Local line bundles from quivers}
we find a realization of these  ideas in the setting of quivers by 
constructing
local
line bundles 
on configuration spaces
from representations of quivers.
Finally, in section \ref{Loop Grassmannians GG P DD(Q,A) and quantum locality}
we get quantum generalization of the notion of local line bundles
and of the corresponding loop Grassmannians
using dilations on the cotangent bundle of moduli of extensions of
representations of a quiver.

Appendix \ref{Drinfeld's theory of classifying pairs}
completes the description 
of Cartan 
fixed points in intersections of closures of semi-infinite orbits in loop Grassmannians
(proposition
\ref{The T-fixed points}).
This observation was the starting point for our generalization
$\GG(Q,A)$ of loop Grassmannians $\GG(G)$.
Appendix
\ref{Comparison with the Thom line bundles in cite{YZell}}
compares computations of Thom line bundles of 
convolution diagrams in \ref{The A-Cohomology of the extension diagram}
and in \cite{YZell}. 
$$
$$

{\em Acknowledgments.}
I.M. thanks Zhijie Dong for  long term discussions on the material
that entered this work. We thank Misha Finkelberg for pointing out errors in  earlier versions.
His advice and  his insistence
have lead to a much better paper.
A part of the writing was done at the conference at IST (Vienna)
attended by all coauthors. We therefore thank the organizers of
the conference and the support of ERC Advanced Grant Arithmetic and Physics
of Higgs moduli spaces No. 320593. The work of I.M. was partially supported by NSF grants.

\se{
Recollections on cohomology theories
}
\lab{Recollections on cohomology theories}

\subsection{
Equivariant oriented cohomology theories}

An {\em oriented cohomology theory} is a contravariant functor $A$
that takes spaces $X$ to graded commutative rings $A(X)$
and  has certain properties such as the proper direct image.\fttt{
While the grading of a cohomology theory is fundamental
we will disregard  it in this paper.
}
For us, 
an oriented cohomological theory $A$ can be either
a topological cohomology theory or an algebraic cohomology theory.
In the first case the ``spaces'' are topological spaces,
and we will use the ones 
that are given by complex algebraic varieties.
In the second case the ``spaces'' mean  schemes  over
a given base ring $\k$.

Here we list some of the common properties of such theories $A$ that we will use.
First, $A$ extends canonically to pairs of spaces $A(X,Y)$
for $Y\sub X$. In particular we get cohomology $A_Y(X)\dff  
A(X,X-Y)$ of $X$ with supports in $Y$.
Such theory $A$ is functorial under flat
pullbacks and proper push forwards with usual properties
(homotopy invariance, projection formula, base change and the
projective bundle formula  \cite{LM,Le}).

Also, such  $A$ has an equivariant version
$A_\sG(X)$ defined as $\limp\ A(\XX_i)$
for  ind-systems of approximations $\XX_i$
of the stack $\sG\bs X$,
For this reason it is consistent to denote $A_\sG(X)$
symbolically as $A(\sG\bss X)$ even if 
we do not really
extend $A$ to category of stacks.

The basic invariants of  $A$ are the 
commutative ring of constants $R=A(\pt)$
and 
the 1-dimensional formal group 
$\bbG$ over $R$ with a choice  of a 
coordinate $\bl$ 
on $\bbG$ (called {\em orientation} of  theory $A$).

The geometric form of the theory $A$ is the functor 
$                                                                   
\fA$ from 
spaces to affine $R$-schemes defined by $\fA(X)=\Spec(A(X))$.
The $\sG$-equivariant version is again 
denoted  by the index $\sG$, it  yields ind-schemes
$\fA_\sG(X)=\Spec(A_\sG(X))$, also denoted $\fA(\sG\bs X)$, that lie above $\fA_G\dff \fA_\sG(\pt)$.
For instance the formal group
$\G$ associated to $A$ is $\fA_\Gm$ (approximations of $\B\Gm$ are given by $\P^\yy$,
the ind-system of finite projective spaces).

We denote $A_G=A_G(\pt)$ and $\fA_G=\fA_G(\pt)$.
For a torus $\sT$ let $X^*(\sT),X_*(\sT)$ be the dual lattices of characters and
cocharacters of $\sT$,
then $\fA_\sT=X_*(\sT)\ten_\Z\G$.
For a reductive group $\sG$ with a Cartan $\sT$ and Weyl group $W$,
$\fA_G$ is the categorical quotient $\fA_\sT//W$.
For instance for the Cartan $\sT=(\Gm)^n$ in $GL_n$, the  Weyl
group is the symmetric group $\fS_n$
and one has $\fA_\sT=\G^n$ while
$\fA_{GL_n}=\bbG^{(n)}$ is the symmetric power
$\G^n//\fS_n$ of $\G$.

\srems (0)
In the case when $\G$ is the germ of an algebraic group $\G_{alg}$
the equivariant $A$-cohomology has a refinement which gives indschemes over $\G_{alg}$.
All of our results extend to this setting and we will abuse the notation by allowing $\G$ to stand either for the formal group
or for this algebraic group.
For simplicity our formulations will assume that $\G_{alg}$ is affine --
the adjustment for the non-affine case
are clear from  the paper
 \cite{YZell} on
elliptic
curves (then $\G$ is an elliptic curve
and $\fA(X)$ is affine over $\G$ rather than affine). 
Either version  satisfies equivariant localization.


(1) For {\em algebraic} oriented cohomology theories
\lab{Algebraic oriented cohomology theories}
the basic reference is
 \cite[Chapter~2]{LM} (one can also use
 \cite[\S~2]{CZZ3} and  \cite[\S~5.1]{ZZ14}).\fttt{
 The terminology of ``algebraic cohomology''
 is also used by Panin-Smirnov for  a refinement of the formalism in which the
 theory is  bigraded (to  adequately encode the example of
 motivic cohomology). We will not be concerned with this version.
 }
Here, cohomology theory is  defined
on smooth schemes over
a given base ring $\kk$. 
However, such cohomology theory $A$ then extends
(with a shift in degrees) under the formalism  of
{\em oriented Borel-Moore homology}
to schemes over $\k$
that are of finite type and separable.\fttt{
What is called Borel-Moore homology here is not quite what this means in classical topology, however this is just a choice of terminology since the $A$-setting does  contain
the precise analogue of Borel-Moore homology.
For instance, for smooth $X$ the more appropriate version would be 
$BM_A(X)=\Th_A(TX)\inv$ in terms of the Thom bundle which is defined next.
}

\sus{
Thom line bundles
}

\lab{Thom line bundle}
When $V$ is a $\sG$-equivariant
vector bundle over $X$, the equivariant cohomology of $V$ supported in the zero section
$\Th_\sG(V)\dff A_\sG(V,V-X)$ is known to be
a line bundle over $\fA_\sG(X)$\ie
a rank one locally free module over $A_G(X)$, called the {\em Thom line bundle} of $V$.
Moreover, this is an ideal sheaf of an effective divisor in $\fA_\sG(X)$
called the {\em Thom divisor of $V$}
(see section 2.1 in
\cite{GKV95}).
As usual, one can think of this as the Thom 
line bundle $\Th(\sG\bs V)$ over
 $\fA(\sG\bs X)$ for the vector bundle $\sG\bs V$ over $\sG\bs X$.




\slem
(a)
Let  $\sV\to \mathsf{X}$ be a vector bundle
equivariant for a reductive
group $\sG$ with a Cartan $\sT$. Then $\Th_{\sT}(\sV)$ is 
the pull back 
of
$\Th_{\sG}(\sV)$ by a flat map $\fA_{\sT}\to \fA_{\sG}$
and $\Th_{\sT}(\sV)$ determines $\Th_{\sG}(\sV)$.

(b)
For a cohomology theory $A$ and a character $\eta$ of a torus $\sT$,
if $\eta$ is trivial, then $\Th_\sT(\eta)=\OO_{\fA_\sT}$, and otherwise
$\Th_\sT(\eta)$ is the ideal sheaf of the (Thom) divisor
$\Ker(\eta)\sub \sT$.

(c) For an extension of vector bundles $0\to V'\to V\to V''\to 0$
one has $\Th_\sG(V)\cong\ \Th_\sG(V')\ten\Th_\sG(V'')$.
So, $\Th_\sG$ is defined on the K-group of $\sG$-equivariant vector bundles over $X$.

\pf
These are proved in \cite[\S~2.1]{GKV95}. See also \cite[Proposition~3.13]{ZZ15}.
\qed

\sss{
$\Th_\sG(\sV)$ for a representation $\sV$ of $\sG$
}\label{sec:function l}
This is the case when $X$ is a point.
We can write $\Th_\sG(\sV)$ in terms of the character $ch(V)$.
First for a torus $\sT$ 
and $\fA_\sT\dff \fA_\sT(\pt)$ there is a unique homomorphism
$\bl:(\N[X^*(\sT)],+)\to (A_\sT,\cd)$
such that for any character $\chi$ of $\sT$, the function
$\bl_\chi$ is the composition
$
\fA_\sT\raa{\fA_\chi}
\fA_\Gm\raa{\bl}\A^1$.
Now, for  a reductive group $\sG$ with a Cartan $\sT$,
this restricts to a homomorphism
$\bl$ from $(\N[X^*(\sT)]^W,+)$ to $(A_\sG,\cd)$.  
Then the ideal $\Th_{\sG}(\sV)$ in functions on
$\fA_{\sG}=[X_*(\sT)\ten\G]//W$
is generated by
the function $\bl_{ch(\sV)}$ on $\fA_{\sG}$.
(By the preceding lemmas, it suffices to check this 
when $\sG$ is a torus, which in turn can be reduced to the case when $\sV$ is one dimensional
and $\sT=\Gm$.)


\sss{
Thom line bundles $\Th(f)$ of maps $f$
}
For a map of smooth spaces $f:\mathsf{X}\to \sY$
we have the tangent 
complex $T(f)=[T\mathsf{X}\to f^*T\sY]_{-1,0}$ on $\mathsf{X}$
and in degrees $-1,0$.
The line 
bundle $\Th(f)=\Th(T(f))$ on $\fA(\mathsf{X})$
is defined as the value of $\Th$  on the corresponding virtual vector bundle $f^*T\sY-T\mathsf{X}$.


\se{
Loop Grassmannians and local spaces
}
\lab{Loop Grassmannians and local spaces}
In 
\ref{Loop Grassmannians}
we recall loop Grassmannians 
and
in
\ref{T-fixed points}
we check the description of $T$-fixed points
in intersections $
(\barr{S_0} \cap \barr{S^-_{-\al}})^T
$
of closures  of semi-infinite orbits
in a loop Grassmannian.
This is used in  
the ``homological'' reconstruction and generalization of 
loop Grassmannians
in
\ref{A generalization GGP(I,QQ) of loop Grassmannians of reductive groups},
which is itself
 based on the formalism of {\em local spaces}
from 
\S~\ref{Local spaces over a curve}.

\sus{
Loop Grassmannians
}
\lab{Loop Grassmannians}
We start with the standard loop Grassmannians $\uGG(G)$.
Let
$\k$ be a commutative ring and let $\OO=\k[[z]]\ \sub\ \KK\ =\ \k((z))$
be the Taylor and Laurent series over $\k$,
these are functions on the indscheme $d$ (the formal disc)
and its punctured version $\ds=\d-0$.
For an algebraic group scheme $G$
we denote by $\GO\sub\GK$ its {\em disc} group scheme and 
{\em loop} group indscheme over $\k$, the   points over a $\k$-algebra $\k'$
are
$\GO(\k')=\ G(\k'[[z]])$
and $\GK(\k')=\ G\b(\k'((z))\b)$.
The standard loop Grassmannian  is  the ind-scheme given by the quotient
$
\uGG(G)=\GK/\GO
$.

We also notice that $\OO_-=\k[z\inv]$ defines group indscheme $\GOm\sub\GK$.
The congruence subgroups $K_\pm(G)$ are the kernels
of evaluations $\GO\to G$ and $\GOm\to G$ at $z=0$ and $z=\yy$.

\sss{
Global (Beilinson-Drinfeld) loop Grassmannians
$\GG_{\HHC}(G)$
}
\lab{Global loop Grassmannians}
Let $C$ be a smooth curve. 
For a finite subscheme $D\sub C$,
the first $G$-cohomology $H^1_D(C,G)$ of $C$ with the support at $D$
is the moduli indscheme
of pairs $(\TT,\tau)$ of a $G$-torsor $\TT$ over $C$ and its section
$\tau$ over $C-D$.
As $D$ varies in the Hilbert scheme of points $\HH_C$
one assembles these  into  
an ind-scheme  $
\GG(G)=\GG_{\HHC}(G)
$
over $\HH_C$,
with fibers $\GG_{\HHC}(G)_D=H^1_D(C,G)$.

A point of a smooth curve $c\in C$ defines a ``smooth formal curve''
$\hatt c$ and we equally get
$\GG_{\HH_{\hatt c}}(G)\to \HH_{\hatt c}$,\
which  is a restriction of $\GG_{\HHC}(G)\to \HHC$.
Moreover, a choice of a formal coordinate 
$z$ 
at $c$ identifies the fiber at $c$ with 
the standard loop Grassmannian
$$
\GG_{ \HH_{\hatt c} }(G)_c
\con\ \uGG(G)
.$$ 
One can also think of it
the compactly supported cohomology
$H^1_c(\hatt c,G)$.\fttt{
By {\em compactly supported cohomology} of  $X$ we mean  
the cohomology of a compactification $\barr X$ trivialized on the formal neighborhood
of the boundary of $X$ in $\barr X$.
}

\sss{
The Abel-Jacobi map
}
\lab{The loop Grassmannian of Gm}
Recall that on a smooth curve $C$ (hence also for  $C=\d$),
$\HH_C$ is a commutative monoid for addition of effective divisors.
Moreover, the {\em Abel-Jacobi map} is a map of monoids:
$$
\AJ_C:\HH_C\to \uGG(\Gm),\ \ \
\AJ_C(D)\dff  \OO_C(-D))
.$$
Precisely, $\AJ_C(D)$ consists of $D\in\HHC$,
the $\Gm$-torsor corresponding to the line bundle $\OO_C(-D)$
and  the canonical trivialization $1$ of $\OO_C(-D)$ off $D$
(which we often omit).

\srem
(0) There are two simple ways to realize all of $\uGG(\Gm)$.
One fixes the trivialization fixed at $1$ as above and the other
uses the trivial torsor $G\tim \d$. 
The transition involves a minus since for any $f\in\KK^*$ one has\
$f\inv:(f\OO_C,1)\con (\OO_C,f\inv)$.


\slem
\cite{CC} (see also \cite{Mi}).
The inclusion  $\d
\sub\HH_\d$
makes
$\HH_\d$  into the
commutative monoid affine indscheme
freely generated by the formal disc $\d$.
The Abel-Jacobi map 
$\d\to \uGG(\Gm)$ makes $\uGG(\Gm)$
into the commutative group indscheme
freely generated by  $d$.\qed

\srems
(1)
In \cite{Mi}, this is viewed as interpretation of
$\uGG(\Gm)$ as homology $\hH_*(d)$ of $d$ for a certain 
conjectural cohomology theory $\hH$.\fttt{
In general, the derived version of  homology $\hH_*(X)$
should be the free abelian commutative group object in derived stacks
freely generated by $X$.
}
The  above interpretations of 
$\uGG(\Gm)$ as both homology and the 
compactly supported cohomology (see \ref{Loop Grassmannians})
can then be viewed as  local Poincar\'e duality in
algebraic geometry.

(2)
A formal coordinate $z$ on $d$ gives a correspondence of
subschemes $D\in\HH_d$
and
monic polynomials $\chi_D$ in $\k[z]$ with nilpotent coefficients, such that
$\chi_D$ 
is an equation of $D$.
This gives a lift  of
the Abel-Jacobi map that  embeds $\HH_d$ into $\GmK$ by
sending $D\in \HH_d^n$ to  $\chi_D\inv$.
For instance for $n\in\N$ the divisor
$n[0]\dff\{z^n=0\}$ goes to $z^{-n}\in\GmK$.

(2) The group indscheme $\GmK$ has a factorization
$\Gm\tim z^\Z\tim K_+(\Gm)\tim K_-(\Gm)$ where the points of congruence subgroups
are $K_+(\Gm)(\k')=1+z\k'[[z]]$
and
$K_-(\Gm)(\k')$ is the invertible part of $1+z\inv\k'[z\inv]$\ie the part
where the coefficients are nilpotent
\cite{CC}.

\sus{
The $T$-fixed points in
semi-infinite varieties $\barr {S_\al}$
}
\lab{T-fixed points and semi-infinite orbits}
\lab{T-fixed points}

\sss{
Tori
}
\lab{Tori}
Let us restate the remarks in
\ref{The loop Grassmannian of Gm}
in the generality
of split tori $T\cong X_*(T)\ten_\Z\Gm$.
First, $\la\in X_*(T)$
gives $\KK^*\to T_\KK$,
and the image of a  coordinate $z$ on the disc is denoted $z^\la$.
This gives isomorphisms $X_*(T)\con \pi_0(\TK)$
and 
$X_*(T)\con \uGG(T)_{reduced}\con\pi_0(\uGG(T))$.
For $\la\in X_*(T)$ we denote
$
L_\la\dff\ z^{-\la}\TO
\in\uGG(T)
$
(independently of $z$), meaning a trivial $T$-torsor on $\d$
with a section $z^{
-\la
}$ on $\ds$.
The corresponding connected component $\uGG(T)_\la$ of $\uGG(T)$ is a 
$K_-(T)$-torsor.


\sss{
The ``semi-infinite'' orbits $S^\pm_\la$
}
\lab{The ``semi-infinite'' orbits S pm la}
Now let $G$ be reductive with a Cartan $T$.
Then
the $T$-fixed point subscheme $\uGG(G)^T$ is $\uGG(T)$.
A choice of  opposite Borel subgroups $B^\pm=TN^\pm$ yields orbits
$S^\pm_\la\dff N^\pm_\KK L_\la$ indexed by $\la\in X_*(T)$
(we often omit the super index $+$).
If $G$ is semisimple then $\uGG(G)$ is reduced and
these orbits provide two stratifications
of $\uGG(G)$.
The following is  well known:  

\slem
\cite{MV}
For $\la,\mu\in X_*(T)$ the following are equivalent:\
(0) $\barr{S_\la}\ni\mu$,\
(i) $\barr{S_\la}\subb S_\mu$,\
(ii) $S_\la$ meets $S_\mu^{-}$,\ and \
(iii) $\la\ge \mu$
(in the sense that $\la-\mu$ lies in the the cone $\cQ^+$ generated by the
coroots $\ch\al$  dual to roots $\al$ in $N$).
\qed

\sex
The loop Grassmannian of $G=SL_2$
is identified with the space $\LL$ of {\em lattices}
in $\KK^2=\KK e_1\pl \KK e_2$\ (the $\OO$-submodules
that
lie between two submodules of form $z^n\OO^2$)
and have volume zero.
Here.
$vol(L)\dff \ \dim(L/z^n\OO^2)-\dim(\OO^2/z^n\OO^2)$ for $n>>0$. 
For the standard Borel subgroup $B=TN$
we
have
$N=\left(\begin{smallmatrix} 1&*\\0&1\end{smallmatrix}\right)
$
and the coroot $\ch\al$ in $N$ is  $\ch\al(a)
=
\left(\begin{smallmatrix} a&0\\0&a\inv\end{smallmatrix}\right)
\in T$.
Then $X_*(T)=\Z\ch\al$
and $L_{n\ch\al}=z^{-n\ch\al}L_0$
is the lattice
generated by two vectors
$\l z^{-n}e_1,z^{n}e_2\rb$.
For a lattice $L\in\LL$ one has 
$L\in S_{n\ch\al}$ if 
$L\cap\KK e_1=z^{-n}\OO e_1$
and
$L\in \barr{S_{n\ch\al}}$ if  $L\ni z^{-n}e_1$.




\sss{
The $T$-fixed points
}
\lab{The T-fixed points}
Now let $G\subb B=TN$ be semisimple with a simply connected cover
$\Gsc\subb\Tsc$,
If $\al_i,\ \ii$, are simple roots  in $N$
then $\prod_\ii\ch\al_i:\GmI\con \Tsc$ defines the
{\em Abel-Jacobi embedding}
$\AJ^G_\d$ as the composition
$\HH_{\d\tim I}
\inj
\uGG(\Gm^I)
\con
\uGG(\Tsc)\sub
\uGG(\Gsc)
\sub \uGG(G)$,\
where\
$D=(D_i)_I\mm\
(\OO_C(-D_i))_I\mm\
(\ch\al_i[\OO_C(-D_i)])_I$.
In particular, for $\al=\sum_Ia_ii\in\NI$,
$$
\AJ^G(\al[0])=\
(\ch\al_i\b[\OO_C(-a_i[0]\b])_I
=\
(\ch\al_i[ z^{a_i}\OO_C)])_I
\aa{\tx{Remark \ref{The loop Grassmannian of Gm}}.0}=\
z^{-\al}
L_0=\
L_\al
.$$

The following  has been 
announced in \cite{Mi}.

\spro 
(a)
The image of the Abel-Jacobi embedding 
$\AJ^G_\d:\HH_{d\tim I}\inj\uGG(G)$
is the fixed point sub-indscheme
$(\barr{S_0})^T
$.

(b) For $\al\in\NI$, this identifies the connected component
$\HH^\al_{d\tim I}
\ =\ \prod_\ii\ \HH^{\al_i}_d
$ of
$\HH_{d\tim I}$
with
the intersection
of
$
\barr{S_0}$
with the connected
component
$\uGG(T)_{\al}
$ of $\uGG(T)$.

(c) Also, the moduli
$
\HH_{\al[0]}\sub\HHdI$
of all subschemes of the finite flat colored scheme $\al[0]$
is   identified with
$
(\barr{S_0} \cap \barr{S^-_{\al}})^T
$.\fttt{
So, the connected components of $\barr{S_0}^T$ are 
$
(\barr{S_0}\cap \barr{S^-_{-\al}})\cap \uGG(T)_{-\beta}$,
for $0\le\be\le\al$,\
identified with the moduli $\HH^\be_{\al[0]}$ of length $\be$
subschemes of $\al[0]$.
}

\nc{\cai}{{ \ch\al }}

\pf
We start with the proof for $G=SL_2$.
Parts (a-b) claim 
that
$\bSz$ meets the connected component $\uGG(T)_{p\cai}$ of
$\uGG(T)$ iff $p\ge 0$ and then the intersection is
$\AJ^G_\d(\HH^p_\d)$. 

The points of the negative congruence subgroup
$K_-(\Gm)\sub\GmK$ 
are the comonic polynomials $Q=
1+q_1z^{-1}+\cddd+q_sz^{-s}$ 
in $z\inv$ with nilpotent coefficients.
Now, the isomorphism
$K_-(\Gm)\con\uGG(T)_{p\cai}$
is given by $Q\mm\ \ch\al(Q)L_p
$
which is the lattice
$
\l Q\inv z^{-p}e_1,
Qz^{p}e_2
\r
$.
According to the example in
\ref{The ``semi-infinite'' orbits S pm la}
this is in $\bSz$ iff $(Qz^p)\inv\OO\ni 
z^0$. This means that $z^{p}Q\in \OO$\ie that  the powers of $z\inv$ in $Q$ are
$\le p$. Such $z^pQ$ form all monic polynomials in $z$ of degree $p$
with nilpotent coefficients.
So, all such $\ch\al(Qz^p)L_0$ form exactly $\AJ^G_\d(\HH^p_d)$.

For part (c)
the example in
\ref{The ``semi-infinite'' orbits S pm la}
says that 
$\barr{S^-_{m\ch\al}}$ consists of lattices $L$ that contain $z^{m}e_2$.
Now, for $D\in\HH^p_d$ with a monic equation
$P\in\k[z]$,
$\AJ^G_\d(D)=\ch\al(P)L_0=\l P\inv e_1,Pe_2\r$
lies in
$\barr{S^-_m}$ iff
$P\OO\ni z^m$\ie polynomial $P$ divides $z^m$.
This is equivalent to $D\sub m[0]$.

The proof in the general case is postponed to the
appendix
\ref{Proof of the proposition The T-fixed points}
\qed

\srem
The connected component 
of $\uGG(G)$
is $\uGG(\Gsc)$, so
the spaces 
$\bSz\subb\ 
\barr{S_0}\cap \barr{S^-_{-\al}}$, and their $T$-fixed points
do not depend on the center of $G$.

\sss{
The Kamnitzer-Knutson program of reconstructing MV-cycles
}
\lab{Fixed points in MV cycles}
Here we restate the proposition
and   recall one of the origins of this paper.
Consider a simply laced semisimple Lie algebra $\fg$
and its adjoint group $G$.
In \cite{BK} the irreducible components $\sC$ of the variety $\La$ of 
representations of  the {\em preprojective algebra} $\Pi$
of a Dynkin quiver $Q$ of $G$
are put into a canonical bijection with certain irreducible subschemes
$X_\sC$ of the corresponding loop Grassmannian $\uGG(G)$, called {\em MV-cycles}
\cite{MV}.

For any  representation $\Vd$ of the preprojective algebra $\Pi$ 
the moduli
$\Gr_{\Pi}(\Vd)$ 
of all $\Pi$-submodules of $\Vd$
is called the  {\em quiver Grassmannian}
of $\Vd$.

\sconj \cite{Mi}
For any irreducible component $\sC$ of $\La$,
and  a sufficiently generic representation $\Vd$ in $\sC$,
the cohomology of its quiver Grassmannian
$\Gr_{\Pi}(\Vd)$ 
is the ring of functions on
the subscheme  $X_\sC^T$ of points 
in the corresponding MV cycle  
$X_\sC$ in $\uGG(G)$
that are
fixed by a Cartan subgroup $T$ of $G$.

The grading on cohomology corresponds to the action of
loop rotations on $X_\sC^T$.

\srems
(0)
This is a version of 
a conjecture of Kamnitzer and Knutson
on equality of dimensions of
cohomology $H^*[Gr_{\Pi}(\Vd)]$ and of 
sections of  the line bundle $\OO(1)$
over the MV cycle $X_\sC$.

(1)
Zhijie Dong has constructed a  map in one direction in this conjecture
\cite{D18}.
 
(2)  The form of this conjecture is alike the
Hikita conjecture
in the 
symplectic duality
framework.

(3) The MV cycles are defined as 
irreducible components of  
intersections in $\uGG(G)$
of closures of semi-infinite orbits
$\bSz\cap \barr{S^-_{-\al}}$ for $\al\in\NI$.
The proposition
\ref{The T-fixed points}.c
will imply the following simplified 
version of the conjecture
that replaces on the loop Grassmannian side the individual MV cycles
with the intersections
$\bSz\cap \barr{S^-_{-\al}}$;
and  on the quiver side it 
degenerates the operators in the representation of $\Pi$ to zero:

\scor
Let $\al\in \NI$, then 
$(\bSz\cap \barr{S^-_{-\al}})^T$
is the spectrum of
cohomology of the quiver Grassmannian
$\Gr_{\Pi}(\Vd)$,\
where $\Vd$ is the zero representation of $\Pi$
of dimension $\al$.
Also, the grading on cohomology corresponds to the action of
loop rotations on the fixed point subscheme of the loop Grassmannian. 

\pf
For $\al=\sum\al_i i$ we have
$\Vd=\pl_\ii\ V_i$ with $\dim(V_i)=\al_i$.
The quiver Grassmannian $Gr_\Pi(\Vd)$ is then the product
$\prod_\ii\ Gr(V_i)$ of total Grassmannians
of components $V_i$.

Since $(\bSz\cap \barr{S^-_{-\al}})^T=\HH_{\al[0]}
=\prod_\ii\ \HH_{\al_i[0]}$ by the proposition
\ref{The T-fixed points}.c,\
it remains to notice that $H^*(Gr_p(n))$ can be calculated by
Carell's theorem as functions on the fixed point subscheme
$Gr_p(n)^e$ of a regular nilpotent $e$ on $\k^n$.
If we realize $\k^n$ and $e$ as $\OO(n[0])$ and the operator of multiplication by $z$,
we see that $Gr_p(n)^e$ is $\HH^p_{n[0]}$
(a subspace of $\OO(n[0])$ is $z$-invariant iff it is the ideal of a subscheme).  

Finally, the degree $2p$ cohomology corresponds to the $p$-power of $z$ which is the grading by
rotations of the disc $d$. 
\qed


\sus{
Local spaces over a curve
}
\lab{Local spaces over a curve}
The notion of local spaces has appeared
in \cite{Mirk} as a common framework for the factorization spaces of Beilinson-Drinfeld
and the factorizable sheaves of Finkelberg-Schechtman.

\sss{
Local spaces
}
\lab{Local spaces}
For a set $I$ and a smooth curve $C$
we consider the Hilbert scheme 
 $\HH_{C\tim I}\cong\ (\HH_C)^I$
 of {\em $I$-colored} points of $C$.\fttt{
 One could try replacing a curve  by a more general scheme
and $\HH$ by other notions of powers of a scheme like the Cartesian powers $\CC^n_C=C^n$.
}
Its connected components  
$\HH^\al_{C\tim I}\cong\ \prod_\ii\ \HH^{\al_i}_C$
are given by subschemes of length  $\al\in\N[I]$.
For  a space $Z$ over 
$\HH_{C\tim I}$ we denote the fiber at $D\in \HH_{C\tim I}$
by $Z_D$.

An {\em $I$-colored local space} $Z$ over $C$
is a space $Z$ over $\HH_{C\tim I}$,
together with a commutative and associative system of  isomorphisms
for disjoint $D',D''\in \HH_{C\tim I}$
$$
\io_{D',D''}:\
Z_{D'} 
\tim
Z_{D''} 
\con\
Z_{D'\sq D''} 
.$$ 
We have $Z_\emp=\pt$. When $\al=i\in I$  the connected component
$\HH^i_{C\tim I}$ is $C\tim i$.
We call the fiber
$Z_{ai}$ at $a\in C$ the ``$i$-particle at $a$''
and we think of
$Z$ as a {\em fusion diagram} for these particles.

\sex
A {\em factorization space} in the sense of Beilinson and Drinfeld
is a local
space $Z\to\HHCI$ such that the fibers $Z_D$ only depend on the formal neighborhood
$\hatt D$ of $D$ in $C$. These can be viewed as spaces
over the {\em Ran space} $\RR_C$, the moduli of finite subsets of $C$.

\srems
(1)
A {\em weakly local} structure is the case when the structure maps $\io$
are only embeddings.
Any weakly local space $Z$ has its
{\em local part} $Z\loc\sub Z$ which we define as the least
closed
local subspace
of $Z$ that contains all particles.
Explicitly, one first constructs $Z^{loc,reg}$ over
$\HH_\CI^{reg}$ so that 
at a discrete $D\in \HH_{C\tim I}$  the fiber is $ 
\prod_{ai\in D}\ Z_{ai}
$, then  $Z\loc$ is the closure in $Z$ of $Z^{loc,reg}$.

(2)
A {\em local structure} on a  (super) vector bundle $V$ over a local space $Z$ 
is an associative and commutative  system of isomorphisms  
$V|_{Z_{D'}} \bten V|_{Z_{D'}}\ \cong\  V|_{Z_{D'\sq D''} }
$. 
By the Segre embedding its projective bundle $\P(V)$ is a weakly local space. 
Its local part 
$\P(V)\loc$ is called the {\em local projective space}
$\P\loc(V)$ of a local vector bundle $V$.
%
%
%
%
%

(3) The notion of  ``locality structure'' is a version of the Beilinson-Drinfeld
``factorization structure'' where emphasis is changed slightly
to get a tool for producing
spaces such as $\P\loc(V)$.
However, the use of closure makes this construction existential rather than explicit.

(4) One would like to extend  this locality mechanism from a smooth curve to a formal disc $\d$
but this requires a supply of disjoint finite flat subschemes of $\d$.
This is expected (or known) to be doable
in terms of rigid geometry.

\sss{
Classification of local line bundles
on the colored Hilbert scheme
}
\lab{Classification of local line bundles}
\lab{Local line bundles and quadratic forms}

The following is a simplified version of the classification of factorizable line bundles
on the space of colored divisor in proposition 3.10.7 of \cite{BD}.

\slem
(a)
For a line bundle
$\sM=(\sM_i)_\ii$ over $\CI$ there is a unique  local line bundle
$\tii\sM$ over $\HHCI$
which agrees with $\sM$ on $\CI$
and whose locality structure maps extend to isomorphisms
across  the diagonals.

(b)
Isomorphism classes of local line bundles
on  $\HHCI$ are classified by pairs $(\sM,\QQ)$
of
a line bundle
$\sM=(\sM_i)_\ii$ over $\CI$ 
and a symmetric bilinear form  (``{\em sb-form}'')
$\QQ$  on $\Z[I]$
by $(\sM,\QQ)\mm \tii\sM(\QQ\De^\HH)$
where
$
\QQ\De^\HH$ is the divisor in $\HHCI$ given by
$
\sum_{i\le j} \QQ(i,j)\De_{ij}
$
for 
the  discriminant divisors
$\De^\HH_{ij}
\sub\HHCI$.

\srems
(0)
We also check the  same classification
for local super line bundles on the
``ordered configuration space''
$\CC_{C \tim I}\dff\ \sq_n\ (C\tim I)^n$.
A line bundle $\sM$ on $\CI$ and
an sb-form $\ka$ on $\Z[I]$ give
the corresponding local super line bundle 
by the ``same'' formula $
\tii\sM(\ka\De^\CC)
$ where the diagonals $\De^\CC_{ij}$ are now in $\CC_{C \tim I}$
and the parity of $\sM^i$ is  that of $\ka(i,i)$.
Local line bundles $\sL$ of the form $\tii \sM(-\ka\De^\CC))$
are characterized by
requiring that the locality maps
$\io^{ij}:
\sL^i\bten \sL^j\to \sL^{i,j}$
defined on $\CC^{i,j}_C=C^2$ and off $\De_C$,
have vanishing of order $\ka(i,j)$ along $\De_C$.


(1) A local super line bundle $\tii\sM(\ka\De^\CC)$ on $\CC_\CI$
descends to $\HHCI$ iff the quadratic form $\ka$ is even.
(Since the pull-back of the diagonal in
$\HH^{2i}_C=C\syp{2}$ to
$C^2$ is the double of the diagonal divisor,\
the pull back of $\OO_\HCI(\QQ\De^\HH)$ from $\HHCI$ to $\CC_\CI$
is $\OO_\CC(\ka\De^\CC)$ where $\ka$ is obtain from $\QQ$ by
doubling the numbers on the diagonal.)


\pf
(a)
In the setting of Cartesian powers,
the restriction $\tii\sM^{i_1,...,i_n}$ of $\tii\sM$
to the connected component $\CC^{i_1,...,i_n}_\CI=\prod_{k=1}^n (C\tim i_k)$
is simply $\bten_{k=1}^n\ M_{i_k}$.

In the setting of Hilbert powers
consider the tautological bundle  $\TT\raa{q}
\HHCI$, the fiber at $D\in\HHCI$ is the subscheme $D$ of $\CI$.
From
$\TT\sub\ 
\HHCI\tim(\CI)\raa{pr_2}\CI$
we have a line bundle $i^*pr_2^*\sM$ on $\TT$
and we define the line bundle $\tii\sM$ on $\HHCI$ as its Deligne direct image in line bundles,
along the map $q$. 
(For an abelian group 
$A$ one has the direct image of $A$-torsors along a finite flat map.)

(b) Let us write the proof in the more general case
of Cartesian powers.
For any $\sM,\QQ$, the line bundle 
$\tii M(\QQ\De^\CC)$
on $\CC_\CI$ is clearly local.
Conversely, let $L$ be any
local line bundle and denote $\sM=L|_{C\tim I}$.


For any $i,j\in I$ 
the locality isomorphism
$L^i\bten L^j\cong L^{i,j}$
is defined on $\CC^{ij}_\CI=(C\tim i)\tim (C\tim j)$ minus the diagonal. 
It extends to  an identification over all of $C^{ij}_\CI$:\ 
$L^i\bten L^j\cong L^{i,j}(-\ka(i,j))$
for a unique integer $\ka(i,j)$.

Now local line bundles $L$ and $\sM(\ka\De^\CC)$
on $\CC_\CI$ have been identified over $\CC^{\le 2}_\CI$.
However an isomorphism  over $\CC^{\le 2}_\CI$ extends uniquely to $\CC_\CI$ since the 
the remaining higher incidences  have codimension $\ge 2$.
\qed

\sus{
A generalization $\GG^P(I,\QQ)$
of loop Grassmannians of reductive groups
}
\lab{A generalization GGP(I,QQ) of loop Grassmannians of reductive groups}

This is a case
of the  local projective space construction
(remark \ref{Local spaces}(2)), in the 
setting of a local line bundle
$L$ on the configuration space $\HHCI$ of colored effective divisors on
a curve $C$. 
In \ref{Local line bundles from loop Grassmannians}
we notice
that for a semisimple group $G$ leads to such local line bundle $\sL$
as a restriction of the line bundle $\OO(1)$ on the loop Grassmannian $\uGG(G)$.
In \ref{Zastava spaces of local line bundles}-\ref{Grassmannians from based quadratic forms}
we will associate to  a based symmetric bilinear form
$(I,\QQ)$ and a poset $P$ the corresponding
zastava space $Z^P(I,\QQ)$
and the  loop Grassmannian $\GG^P(I,\QQ)$.
Finally, in \ref{Reconstruction of loop Grassmannians associated to groups}
we check that this is indeed a generalization of the corresponding spaces
for simply connected semisimple groups.

This reconstruction roughly says that
the loop Grassmannian $\GG(G)$
can be effectively reconstructed
from $\GG(T)$.
Previous results in this direction
include 
\cite{Zhu},
\cite{FK},
\cite{Se}.
The key observation is that the equations of 
the semi-infinite variety $\bSz$
in the projective space
$\P(\Ga[\bSz,\OO_{\uGG(G)}(1)]^*)=\P(\Ga[\bSz^T,\sL^*])$
are given by the locality structure
on $\OO_{\uGG(G)}(1)$.
So, the locality structure allows us to reconstruct
$\bSz$
and then also $\uGG(G)$ as a certain limit of
copies of $\bSz$.

\label{subsec:o1}

\sss{
Local line bundles from loop Grassmannians
}
\lab{Local line bundles from loop Grassmannians}
If $G$ is simple and simply connected then
$Pic[\uGG(G)]\cong \Z$ with the canonical generator $\OO(1)$
(given for instance by the divisor which is the complement of the
open Bruhat cell $\uGG^0$ in $\uGG(G)$).
To study $\OO(1)$ we will use   {\em factorizable line bundles} on various versions of
loop Grassmannians, these  are defined and compared in
\cite{TZ}.

We now choose two relevant versions of the Abel-Jacobi map.
For a smooth curve $C$  define the maps
$
\CCCI
\CD
@>x\mm\barr x>>
\HHCI
@>N>>
\HHC
\eCD
$.
For
 $x=(c_1,i_1,...,c_n,i_n)\in\CC^n_\CI=(\CI)^n$
(so, $c_p\in C,\ i_p\in I$),
let 
$\barr x=\sum_\ii ix_i\in\HHCI$ with
$x_i=\sum_{i_p=i}\ c_p\in\HH_C$.
Also,
for $D_i\in\HHC$, let $N(\sum_\ii iD_i)=\sum_I\ D_i$.
For a semisimple group $G$
consider the global  Abel-Jacobi map
$$
\aj^G_C\dff
\b(
\CC_\CI
\CD
@>x\mm\barr x>>
\HHCI
@>\AJ^G_C>>
\eCD
\GG_\HHC(G)\b)
,\ \ \
\aj^G_C(x)\dff\ \AJ^G_C(\barr x)=
\b(
N\barr x,\OO_C(-\barr x)
\b)
,$$
with
$
N\barr x\in\HH_C
$
and
$
\OO_C(-\barr x)
\in
\GG_\HHC(T)_{N\barr x}
$.
We also consider a version at one point $0\in C$ :
$$
\aj^G_\d\dff
\b(
\CC_\dI
\CD
@>x\mm\barr x>>
\HHdI
@>\AJ^G_\d>>
\eCD
\uGG(G)\b)
,\ \ \
\aj^G_\d(x)\dff\ \AJ^G_\d(\barr x)\dff
\OO_C(-\barr x)
.$$

\spro
Let
$G$ be simple and simply connected.
Then the pull back
$(\AJ^G_\d)^*\OO(1)$ of $\OO(1)$
via the above Abel-Jacobi map 
is the line bundle
$\OO_{\CC_\dI}(-\ka\De^\CC)$ on $\CC_\dI$
corresponding to the negative of 
the Cartan  form $\ka$ on $\Z[I]$.
It has a canonical extension
to
a local line bundle
$\OO_{\CC_\CI}(-\ka\De^\CC)$ 
on $\CCCI$.

\pf
We will use the curve 
  $C=\Ao$ which  contains the disc $\d=\hatt 0$.
Then   $\uGG(G)$ is the fiber
 $\GG_\CCC(G)_0$ of the Beilinson-Drinfeld Grassmannian $\GG_\CCC(G)$
 at the divisor $0\in C=\CC^1_C$.
We
use
the extension of
$\OO(1)$ to a factorizable line bundle $\OO^\CC(1)$ on 
the loop Grassmannian $\GG_\CCC(G)$. By
proposition 2.5 in  \cite{TZ}
this extension is unique once we
trivialize $\OO(1)$ at the origin of $\uGG(G)$
(this is essentially the section 3.4 in \cite{Zhu_PCMI}).

The restriction
$\LL$ of  $\OO^\CC(1)$ to 
$\GG_\CCC(T)\sub\GG_\CCC(G)$ inherits the structure of  a factorizable line bundle.
This will imply that
$(\aj^G_\d)^*\OO(1)$ is canonically identified with the restriction of 
$(\aj^G_C)^*\OO^\CC(1)$ to $\HHdI
\sub\HHCI$.
The point is that at $x\in\CCdI$, the fiber of
$(\aj^G_\d)^*\OO(1)$ is $\LL_{0,\OO_C(-\barr x)}$,
and
the fiber of
$(\aj^G_C)^*\OO^\CC(1)$ is $\LL_{N\barr x,\OO_C(-\barr x)}$.
The canonical identification
of these fibers of $\LL$ is a special case of the unique
descent
of any factorizable line bundle on $\GG_\CCCI(T)$ to the so called {\em rational
Grassmannian} $\GG_{rat}(T)$, this was proved in proposition 1.4 of \cite{TZ}.

We are now interested in the line bundle 
$(\aj^G_C)^*\OO^\CC(1)$ on $\CCCI$ with its  structure of
a local line bundle that it inherits
-- by definitions --
from
the 
factorization line bundle structure on $\OO^\CC(1)$.
To any factorization line bundle $L$ on $\GG_\CCC(G)$
one associates a  quadratic form $q_L\in QF(X_*(T))$
whose symmetric bilinear form $\ka_L(\la,\mu)=q_L(\la+\mu)-q_L(\la)-q_L(\mu)$
is given by the order of vanishing of  the locality structure
along the diagonals $\De_C\sub C^2$ corresponding to the pair $\la,\mu$.

The proof of the proposition 2.5 in \cite{TZ}
checks that if $G$ is simple and $\pi_1(G)=0$ the quadratic form
corresponding to $\OO^\CC(1)$
is the ``minimal'' integral invariant quadratic form $q_m$,
characterized by $q_m(\ch\al)=1$ for short coroots $\ch\al$,

Finally, in the simply connected case
this is the Cartan matrix $\ka$ in the sense that $\ka_{ij}\dff\l \al_i,\ch\al_j\r$ equals
$(\ch\al_i,\ch\al_j)$ for the sb-form $(-,-)$ given by the quadratic form $q_m$.
\qed

\sss{
Zastava spaces of local line bundles
}
\lab{Zastava spaces of local line bundles}
We can induce any local line bundle $\sL$ on $\HHCI$
along
a  poset $P$.
Consider the 
correspondence 
$$
(\HHCI)^P
\laa{\pi}
^P\HH\
\raa{\si}\
\HHCI
,$$
where the fiber $\si\inv\{D\}$
of $^P\HH$ at $D\in\HHCI$ 
is the set $\Hom(P,\HH_D)$ of maps $E$ of posets\ie
for $a\le b$ in $P$ one has $E^a\sub E^b\sub D$.
This gives a sheaf
$
Ind^P(\sL)\ \dff\ 
\si_*\pi^*(\sL^{\bten P})
$ on $\HHCI$.

\slem
$Ind^P(\sL)$  is  a local sheaf on $\HHCI$.
 If $P$ is an interval
$[m]=(1<
\cddd <m)$ then $Ind^P(\sL)$
is a local vector bundle.\fttt{
Flatness fails for the poset
$P=(0<a,b$).
}

\pf
If $D=D^1\sq D^2$ in $\HHCI$ then a representation of $P$ in $D$ is a pair
of representations in $D^i$'s.
Sheaf    $Ind^P(\sL)$ is a vector bundle if
$^P\HH$ is flat over \
$\HHCI$.
If $D$ is a length $n$ subscheme of a smooth curve then the length of
the Hilbert scheme $\HH_D$ is $2^n$.
This gives the case $m=1$. The general case follows by induction
in posets $i<\cddd<m$ as $i$ goes from $m$ to $1$.
\qed

Define  the  {\em $m\thh$ zastava space}
of $\sL$
as the local projective space
$$
Z^m(\sL)\dff\ \P\loc([Ind^{[m]}\sL]^*)
.$$

\srems (0)
One can extend the construction
to zastava spaces $Z^P(\sL)$ for finite posets $P$,
by replacing the projective space of $Ind^P\sL{}^*$
with 
 the projective spectrum of the symmetric algebra of $Ind^P\sL$.

(1)
Any symmetric bilinear form $\QQ$ on  $\Z[I]$
defines a local line bundle $\OO(\QQ)=\OO_\HHCI(\QQ\De^\HH)$
on  $\HHCI$ (\ref{Classification of local line bundles}), hence
also zastava spaces $Z^P(I,\QQ)\dff Z^P(\OO(\QQ))$.

\sex
When $P$ is a point we omit $P$ from notation.
Then $^P\HH\raa{\si}\HHCI$
is the relative Hilbert scheme $\HH_{\TT/\HHCI}$ for
the tautological bundle $\TT/\HHCI$\ie
the fiber $\TT_D$ at $D\in\HHCI$ 
is the Hilbert scheme $\HH_D$ of all subschemes of the finite scheme $D\ \sub \CI$.

If $D$ is a point $ai$ with  $a\in C$ and $\ii$
then $\HH_{ai}=\{\emp,ai\}$,
hence
$Ind(\sL)=\sL_\emp\pl\sL_{ai}=\k\pl\sL_{ai}$.
So, the particle at $ai$ is $\Po$ with two chosen points.
Then one is constructing
$Z(\sL)$ by colliding $\Po$'s according to the prescription
given by the locality structure on the line bundle $\sL$.

Similarly, when $P=[m]=\{1<\cddd<m\}$ then all
particles of $Z^m(\sL)= Z^{[m]}(\sL)$
are $\P^m$.
Actually, this $\P^m$ is naturally the $m\thh$ symmetric power of
the particle $\Po$, so $\Po$ embeds into $\P^m$  by Veronese embedding.

\sss{
Weak flatness of zastavas
}
\lab{Flatness of zastavas}
We start with a general statement about Kodaira embeddings.

\slem
Let $\LL$ be a relatively very ample line bundle on a projective
scheme $Y$ over a base variety $S$.
Let $X$ be  a projective subvariety over an open dense $U\sub S$.
Suppose that 
there exists a finite flat $S$-subscheme  $F
\sub Y
$ such that over $U$\
it lies inside $X$ 
\ and 
the restriction map
$\LL_{X/U}\to (\LL_{F/S})|_U$
is an isomorphism.\fttt{
We denote by $\LL_{F/S}$
the direct image of $\LL|_F$ to $S$ etc.
}
Then  the
closure $\barr X$ of $X$ in $Y$
is a closed subscheme of $\P(\LL_F{}^*)\cap Y$
and also  $\LL_{\barr X/S}=\LL_{F/S}$ is a vector bundle.

\pf
For a sheaf $\EE$ denote by $\bP(\EE)$ the projective spectrum of the symmetric algebra
of $\EE$.
The restriction map
$\LL_{Y/S}\to \LL_{F/S}$ is surjective since $\LL$ is very ample for
$Y/S$ and $F/S$ is finite flat.
Therefore, inside
$\bP( \LL_{Y/F})$ we have both   $Y$ and $\bP(\LL_{F/S})$.
Over $U$,
both contain $X$ since
$\bP(\LL_{F/S})|_U=\bP\LL_{X/U}$ and this contains $X$ since $\LL$ is very ample.
So, $Y\cap \bP(\LL_{F/S})$ also contains $\barr X$.

The restriction $\rho:\LL_{\barr X/S}\to \LL_{F/S}$ is surjective since
we have $F\sub \barr X\sub Y$
and  $\LL_{Y/S}\to \LL_{F/S}$ is surjective since
$\LL$ is very ample for $Y/S$ and $F/S$ is finite flat.
Moreover, if $S$ is a variety then $\Ker(\rho)$ is supported over the boundary $S-U$ of $U$.
If $X$ is a variety then  so is $\barr X$.
Moreover, $\Ker(\rho)=0$ since the line bundle $\LL$ has no sections on $\barr X$
that vanish on $X$.
So, $\LL_{\barr X/S}=\LL_{F/S}$ and this is a vector bundle since $F/S$ is finite flat.
\qed

\scor
Recall that the  zastava spaces $Z^m(\sL)$ 
lie inside the bundle of projective spaces 
$\P[Ind^{[m]}(\sL)^*]$ which carries the line bundle $\OO(1)$. 
Then the sheaf
$(Z(\sL)/\HHCI)_*\OO(1)$
is the vector bundle $Ind^{[m]}(\sL)$.

\pf
In the lemma we  choose
$S\subb U$ as $\HHCI\subb\HH^{reg}_\CI$.
Let $F/S$ be the map $\si$ in the correspondence
$\HHCI\laa{p}\ ^{[m]}\HH\raa{q}\HHCI$ from
\ref{Zastava spaces of local line bundles}

A local line bundle $\sL$ on $S$ gives a line bundle $\pi^*\sL^{\bten m}$ on $F$
and $(\pi^*\sL^{\bten m})_{F/S}$ is just $Ind^{[m]}(\sL)$.
Now, the projective bundle   $Y/S\dff\ \P(\pi^*\sL_{F/S})/S$ carries a very ample line bundle
$\LL=\OO(1)$ and contains $\barr X=\P\loc([\pi^*\sL]_{F/S})=Z^m(\sL)$
and $X=\P^{loc,reg}([\pi^*\sL]_{F/S})=Z^{m,loc}(\sL)$.
So, the claim 
$\LL_{\barr X/S}=\LL_{F/S}$ from the lemma
means that $(Z^m(L)/\HHCI)_*\OO(1)$ equals
$(^{[m]}\HH/\HHCI)_*\OO(1)$. However,
$\OO(1)|_{^{[m]}\HH}$ is just $\pi^*\sL^{\bten m}$
(a property of Kodaira embeddings), and 
$(^{[m]}\HH/\HHCI)_*\pi^*\sL^{\bten m}$
is $Ind^{[m]}(\sL)$.
\qed

\sss{
Grassmannians from based symmetric bilinear  forms
}
\lab{Grassmannians from based quadratic forms}
Let $\QQ$ be a quadratic form on  $\Z[I]$.
Its zastava space $Z^P(I,\QQ)\dff Z^P(\OO(\QQ))$
defines
the
{\em semi-infinite} space $S^P(I,\QQ)$ over $\HHCI$,
the fiber at 
$D\in\HHCI$ is the colimit (union)
of 
zastava fibers $Z^P(I,\QQ)_D$
at multiples of $D$
$$
S^P(I,\QQ)_D\ =\ 
\limi_n\
Z^P(I,\QQ)_{nD}
.$$

Finally, $\N[I]$ acts on $S^P(I,\QQ)$,
and the corresponding {\em loop Grassmannian} 
is defined as
$$
\GG^P(I,\QQ)\ \dff\ \ \Z[I]\tim_{\N[I]}\  
S^P(I,\QQ)
.$$

\sss{
Zastava spaces $Z^\al(G)$ for groups
}
\lab{Zastava spaces Zal(G) for groups}
We will recall these spaces from \cite{FM}
which compares several definitions.
The description  in terms of 
the Beilinson-Drinfeld Grassmannian
$\GG^{BD}(G)$ is from  section 6 of  \cite{FM}.
%

We will be interested in the connected component
$\GG^{BD}(\Gsc)$ of $\GG^{BD}(G)$. 
For a Borel $B$, the simple coroots identify $H=B/N$ with $\Gm^I$.
Therefore the singularities of a rational section $\tau$ of any $H$-torsor over $C$
define an $I$-colored divisor $D$.
We say that the  {\em singularity of $\tau$} is the degree 
$\deg(D)\in \Z[I]$.

For $\al\in \cQ^+$ we define 
$S^{BD}(\al)\sub\GG^{BD}(\Gsc) $
so that
the fiber $S_D(\al)$ 
at $D\in\HHC$ consists of all pairs
$(\TT,\tau)$ of a $G$-torsor $\TT$ over $C$ and its section
$\tau$ over $C-D$, such that the $B$-subtorsor $B\tau\sub \TT|_{C-D}$
extends to $C$ (necessarily as the closure $\barr{B\tau}$ in $\TT$),
and that the section $\tau$ of the $H$-torsor
$N\bs \barr{B\tau}$ has singularity $\al$.
If $\al\in\N[I]$ this means that the divisor $D$
is in $\HH^\al_\CI$.
At a single point $a\in C$ the fiber $S_a(\al)$ is the orbit $\NK L_\al$ from
\ref{The T-fixed points}.

Now a pair of opposite Borels $B^\pm=TN^\pm$ (with $B^+=B$) gives 
two semi-infinite stratifications
$S^{BD,\pm}(\al),\ \al\in \cQ^+$,\
of $\GG^{BD}(\Gsc)$.
For $\al\in\N[I]$ consider a version $\bZ^\al(G)$
of ``zastava'' for $G$ which is obtained by
pulling back the intersection
of closures  
$
\barr{
S^{BD,+}_0
}
\cap
\barr{
S^{BD,-}_{\al}
}$
to the connected component
$\HH^\al_\CI$ of the colored configuration space,
via the
sum map $\HHCI\to\HHC,\ (D_i)_\ii\mm\sum D_i$.
Our ``zastava'' space  $\bZ\to\HHCI$ is local
since the singularity of a section $\tau$ regular off
$D'\sq D''$ is the sum of contributions at $D'$ and $D''$.

When
$ai$ is a point in $C\tim I=\HH^i_\CI$
then the fiber
$(\bZ^i)_{ai}$ is 
$\Po$ with Cartan fixed points $L_0,L_{\ch\al_i}$
(this reduces to $G=SL_2$
and then it is as easily seen from the example in \ref{The ``semi-infinite'' orbits S pm la}).
Then 
by locality, the fiber at a regular divisor
$
D\in \HH^{\al,reg}_\CI
$
is isomorphic to 
$(\Po)^D$, and the  $T$-fixed points
in $\bZ^\al_D$ are all $L_{\be}$ with $0\le \be\le \al$ in $\cQ$.

This in particular shows that
the restriction
$Z^{\al,reg}\dff \bZ^\al\b|\HH^{\al,reg}_\CI$ is reduced.
The zastava space that we are interested in
is its closure  $Z^\al$ in
the scheme $\bZ^\al$.
By definition. it is a reduced local space. 

\nc{\ZZz}{\aa{\ci}\ZZ}
\srems
(0) [Comparison of terminology with \cite{FM}.]
The space $Z$ over $\HHCI$ is projective so we can call it {\em projective zastava}.
We also have open subspaces
$Z\subb\ZZ\subb \ZZz$.
What is called {\em open zastava} in  \cite{FM}
is
$\ZZz$ which is obtained by  replacing
$
\barr{
S^{BD,+}_0
}
\cap
\barr{
S^{BD,-}_{\al}
}$
with
$
S^{BD,+}_0
\cap
S^{BD,-}_{\al}
$. For $C=\Ao$ this is the space of {\em based
maps} into the flag variety
$Map[(\Po,\yy),(\BB,\fb_-)]$.
What is called {\em zastava} in \cite{FM} is the space $\ZZ$ for which one removes one closure
and uses
$
S^{BD,+}_0
\cap
\barr{
S^{BD,-}_{\al}
}$.
For $C=\Ao$ this is the space of {\em based quasimaps} $QMap[(\Po,\yy),(\BB,\fb_-)]$. We can call it
{\em affine zastava} since it is affine.

(1)
The diagonal  points of $\HH^\al_\CI$ are of the form $\al p=(a_ip)_\ii$ for
a point $p\in C$ and $\al=\sum_\ii a_ii\in\N[I]$.
The unions of central fibers of zastava spaces
$\cup_{\al\in\cQ_+}Z^\al_{\al p}\subb
\cup_{\al\in\cQ_+}\ZZ^\al_{\al p}$
are by definitions the semi-infinite varieties $\barr{S_0}\subb S_0$
in $\uGG(G)$.

Moreover,
one can view $\uGG(G)$ as the increasing union of
$\barr{S_{\be}}=z^{-\be} \barr{S_{0}}\cong \barr{S_0},\ \be\in\cQ^+$.
So, the whole loop Grassmannian is a certain direct limit
$\limi_{\cQ^+}\ \barr{S_0}$.
As $\al\in \N[I]$ acts on $\barr{S_0}$ by $z^\al$,
we can rewrite this limit as $\uGG(G)=\Z[I]\tim_{\N[I]}\barr {S_0}$.

\sss{
Reconstruction of loop Grassmannians $\GG(G)$
associated to groups
}
\lab{Reconstruction of loop Grassmannians associated to groups}
Here we check that
the Grassmannians $\GG^P(I,\QQ)$
from \ref{Grassmannians from based quadratic forms}
indeed generalize the loop Grassmannian $\uGG(G)$
of semisimple groups.

\stheo
Let $G$ be a simply connected
simple group which is simply laced.
Let $I$ index the simple roots and let $\QQ'$ be the modification of the Cartan matrix
by dividing by $2$ on the diagonal.\fttt{
The modification appears because we use the Hilbert scheme $\HHC$ rather than over
powers of curves $\CCC$ (see the remark \ref{Classification of local line bundles}(1)).
}
Then
$Z^\al(I,-\QQ'),\ S(I,-\QQ'),\ \GG(I,-\QQ')$
are naturally identified with the corresponding notions
$Z^\al(G),\ \barr{S_0}(G),\ \GG(G)$
for the group $G$.

\pf
We only need to prove equality $Z^\al(I,\QQ)=Z^\al(G)$
since
the other two spaces are obtained from zastavas in the same way
(compare  \ref{Grassmannians from based quadratic forms}
and the remark \ref{Zastava spaces Zal(G) for groups}.1).

Let
$\TT^\al $
be the tautological bundle over $\HH^\al_\CI$
so that the fiber at a colored divisor $D=\sum_\ii\ D_ii$
is the finite flat scheme $D=\sq_i\ D_ii\sub C\tim I$.
The interesting object is the relative Hilbert scheme $\HH(\TT^\al/\HH^\al_\CI)$, its fiber at $D$
is 
$\HH(D)\dff \prod_\ii\ \HH_{D_i}$,
where, $\HH(D_i)$ is the moduli of all subschemes of a finite flat scheme $D_i$.
The point is that  according to the proposition
\ref{The T-fixed points}.
the  fixed points  subscheme $(Z^\al(G))^T$ is exactly
$\HH(D)\dff \prod_\ii\ \HH_{D_i}$.\fttt{
While the quoted  lemma is stated at a single point of $C$,
we actually  need a version of the lemma
for the family
$Z^\al(G)\to \HH^\al_\CI$. This is easy using the moduli description of semi-infinite stratifications
from \ref{Zastava spaces Zal(G) for groups}.
}

Now we can apply 
the weak flatness lemma \ref{Flatness of zastavas}.
Choose $S\subb U$ to be $\HH^\al_\CI\subb \HH_\CI^{\al,reg}$.
For $Y/S$ we choose the pull back of the Beilinson-Drinfeld
Grassmannian
$\GG^{BD}(G)$ to $\HH_\CI$ (by the sum map).
It carries a very ample line bundle $\LL$ which is the pull back of
$\OO_\GG(1)$ from $\GG^{BD}(G)$.

For $F\sub Y$ take $(Z^\al(G))^T
=\HH(\TT^\al/\HH^\al_\CI)$. We have already argued that it is flat
in the proof of lemma 
\ref{Zastava spaces of local line bundles}.
Finally, let $X/U$ be
$Z^{\al,reg}(G)$,
the restriction of $Z^\al(G)$ to  $\HH_\CI^{\al,reg}$.

We still need to check that the restriction map
$\sL_{X/U}\to \sL_{F_U/U}$ is an isomorphism.
Due to locality of
$Z(G)$  
and of the line bundle $\LL$ on $Y$,
we only need the claim at colored divisors
$D$ which are points $ai$.

Here, $Z^{i}_{ai}$ is $\Po$ with $T$-fixed points
$\{L_0,L_{\ch\al_i}\}$. So, one just needs to check that the restriction of
$\sL$ to $Z^{i}_{ai}$ is $\OO_\Po(1)$.
Recall that $\sL=\OO_\GG(1)$ is $\OO_{\GG(G)}(\DD^G)$ for the divisor
$\DD^G$ in $\GG(G)$ defined as the boundary of the open orbit
$\GG^0(G)$ of the negative congruence subgroup.
So, it suffices to see that $Z^{i}_{ai}\cap \DD^G$ is a single point.
This claim reduces to the $SL_2$ subgroup $S_i$ of
$G$ corresponding to the simple root $\al_i$\ --\
since $D^G\cap\GG(S_i)=D^{S_i}$.
In the $SL_2$ case one easily checks explicitly in the lattice model
for the loop Grassmannian (see the example \ref{The T-fixed points}).

Now,  lemma \ref{Flatness of zastavas} guarantees that
for $Z^\al(G)=\barr{Z^{\al,reg}(G)}=\barr X$
and any
$D\in \HH^\al_\CI$ the restriction map
$\Ga[Z^\al(G)_D,\OO(1)]\to
\Ga[Z^\al(G)_D^T,\OO(1)]
$ is an isomorphism.
This implies that $Z^\al(G)$ is a local projective space for the local vector bundle
dual to $\LL_{Z^\al(G)^T/\HHCI}$ (we use  the notation from lemma \ref{Flatness of zastavas}).

Now $Z(G)^T$ has been identified with
the Hilbert scheme for $\TT/\HHCI$ and then the local vector bundles
$\LL_{Z^\al(G)^T/\HHCI}$ and $Ind(I,-\QQ')$
on the two spaces are identified by the proposition
\ref{Local line bundles from loop Grassmannians}.
Therefore, the corresponding local projective spaces
$Z(G)$ and $Z(I,-\QQ')$ are the same. 
\qed

\srem
The ``simply laced'' restriction
can be removed using the 
folding mechanism. Here this means
an action of a group $\ga$ on the data $(I,\QQ)$.
Furthermore,  if $\ga$ is allowed  to act on a curve above $C\tim I$
then one  would include loop Grassmannians for twisted affine Lie algebras
and the Galois action.

\sss{
The ``homological'' aspect of $\GG^P(I,\QQ)$
}
\lab{The ``homological'' aspect of GGP(I,QQ)}
The standard interpretation of loop Grassmannians $\uGG(G)$ is cohomological 
(\ref{Global loop Grassmannians}).
The construction \ref{A generalization GGP(I,QQ) of loop Grassmannians of reductive groups}
provides what one could view as a homological interpretation.
First, for a torus $T=\GmI$,
one  builds 
$\uGG(T)$ from the union of formal discs $\d\tim I$ in stages
$\d\tim I\mm\HH_{\d\tim I}\mm \uGG(\GmI)$.
Here, the configuration space $\HH_{\d\tim I}$
is the free commutative monoid generated by $\d\tim I$, and then 
one inverts the center of the disc
to get 
$\uGG(\GmI)$ as the free commutative group
on $\d\tim I$
(remark 0 in \ref{The loop Grassmannian of Gm}).
One repeats this procedure
for a reductive group $G$ with a Cartan $T=\GmI$
by adding a local line bundle
$\sL$ on $\HH_{\d\tim I}$,
to get $\uGG(G)$ (or more generally $\GG^P(I,\QQ)$).
First, the {\em positive part} of
the loop Grassmannian is the zastava space
$Z^P(I,\QQ)$
built using the monoid
$\HH_{\d\tim I}$
in \ref{Zastava spaces of local line bundles}.
Then
$\GG^P(I,\QQ)$ itself is obtained from
$Z^P(I,\QQ)$ in  
\ref{Grassmannians from based quadratic forms}
by inverting $\N[I]\sub\HH_{\d\tim I}$.
We mention that for a complete  curve $C$ a reconstruction of $Bun_G(C)$ from $C$
has been pursued long ago in
\cite{Feigin-Stoyanowski}.

\se{
Local line bundles from quivers
}
\lab{Local line bundles from quivers}


We know that local line bundles $\sL$ on configuration spaces
$\HHCI$ correspond to
quadratic forms $\QQ$ (\ref{Grassmannians from based quadratic forms})
and the forms $\QQ$ with
non-negative integer 
coefficients clearly
correspond to graphs.
In this  section we construct local line bundles directly from
graphs or quivers. The advantage is that such
construction
 extends to the quantum setting
(see \cite{YZ2} and \ref{DD-quantization of the monoid HGI} below).
In the quantum setting the ``commutative'' configuration space
$\HGI$ will be replaced with the ``non-commutative''\ie  ordered,
configuration space $\CC_{\G \tim I}\dff\ \sq_n\ (\G\tim I)^n$.
On the level of representations of quivers
the noncommutative configuration space
corresponds to passing to
complete
flags in representations.

We start with the  curve $\G$ which is  the 1-dimensional  group corresponding to
a cohomology theory $A$.\fttt{
Here, $\G$ is defined over the ring of constants $A(\pt)$ of the cohomology theory.
}
Then
the Hilbert scheme $\HH_{\G\tim I}$ of points in $\GI$
is obtained as the cohomological schematization
$\fA(\Rep_Q)$
of the moduli $\VV^I$ of $I$-graded finite dimensional vector spaces.


In \ref{Quivers} we recall various categories of representations of quivers
and their extension correspondences.
The cotangent complexes for these correspondences
are considered in
\ref{The cotangent versions of the extension diagram}.

The ``classical'' {\em local} and {\em biextension} line bundles
$\sL(Q,A)$ and $\LL(Q,A)$ on $\HGI$ and $(\HGI)^2$
are constructed as Thom line bundles of moduli of extensions of representations in
\ref{Thom line bundles Th(V) of vector bundles V}.
Here, $\sL(Q,A)$ can be defined directly 
from the incidence
quadratic form of the quiver $Q$.

In 
\ref{The A-Cohomology of the extension diagram} we calculate Thom line bundles associated to 
the cotangent correspondence and the effect of dilations.
Finally, in
\ref{DD-quantization of the monoid HGI}
we recall the construction of the quantum group $\UU_\DD(Q,A)$ from the cotangent correspondence and this 
leads us to select a choice of quantization of the above ``classical'' line bundles
from \ref{Thom line bundles Th(V) of vector bundles V}.

\srem
This section is largely a retelling of  the paper \cite{YZell}.
That paper is primarily concerned with the construction of quantum affine groups 
in the language of preprojective algebras
which is here viewed as the cotangent bundle
of the moduli $\Rep_Q$.
This ``symplectic'' setting allows to 
``quantize'' the notion of local
line bundles and the construction of  loop Grassmannian
from local line bundles.
The quantization comes from the action of the dilation torus $\DD$ on representations
(which is in turn defined by a choice of 
a {Nakajima  function} $\bm$ on the set of arrows of the double $\bQ$ of the quiver $Q$).\fttt{
While \cite{YZell}
deals with the case of elliptic cohomology,
some of its ideas appear in an earlier paper 
\cite{YZ1} which was only concerned with affine groups $\G$.
This allowed for a  trivialization of Thom line bundles
which accounts for a
 different presentation
of   functoriality of cohomology in that paper.
}

\subsection{
Quivers
}
\lab{Quivers}
Let $Q$ be a quiver with finite sets $I$ and $H$ of vertices
and arrows. 
For each arrow $h\in H$, we denote by $h'$ (resp. $h''$) 
the tail (resp. head) vertex of $h$.
The opposite quiver $Q^*=(I,H^*)$ has the same vertices
and the set of arrows $H^*$ is endowed with a bijection $*:H\to H^*$,  so that $h\mm h^*$
exchanges sources and targets. The double $\bQ$ of the quiver $Q$
has vertices $I$ and arrows $H\sqcup H^*$.

Let $\VVI$
be the moduli of finite dimensional 
$I$-graded vector spaces 
$V=\pl_\ii\ V^i$. 
Let $\Rep_Q$ be the moduli of representations 
of $Q$.
Its fiber at $V\in\VVI$ is the vector space
$\Rep_Q(V)$ of  representations
on $V$.
This is the sum over $h\in H$ of $Rep_Q(V)_h=\Hom(V^{h'},V^{h''})$.
We usually denote $v=\dim(V)\in\N[I]$ and let $G$ be $\GL(V)$
so that
the  connected component
$\Rep_Q(v)$ -- given by representations  of fixed dimension
$v\in\N[I]$ -- is $G\bs Rep_Q(V)$.

\sss{Dilation torus $\DD$}
\label{sss:dilation}
A choice of Nakajima's weight function $\bm:H\coprod H^*\to \bbZ$ gives an action of $\bbG_m^2$ on 
\[
T^* \Rep_Q(V)=\Rep_\bQ(V)=\Rep_Q(V)\pl\Rep_{Q^*}(V).
\]  
Elements $(t_1,t_2) $ act for each $h\in H$ on $\Rep_Q(V)_h$
by $t _1^{\bm_h}$ and on $\Rep_{Q^*}(V)_{h^*}$ by
$t _2^{\bm_{h^*}}$.
We also let $\Gmt$ act on the Lie algebra $\g$ of $\GL(V)$ by $t_1t_2$.

We choose a subtorus $\DD$ of $\Gmt$ and 
require that  the moment map for the $GL(V)$-action on
$T^*\Rep_Q(V)$ is $\DD$-equivariant. This means that on $\DD$ we have
$t_1^{\bm(h)} t_2^{\bm(h^*)}=t_1 t_2$\ for any $h\in H$.
In particular, the symplectic form on
$T^* \Rep_Q(V)$ has weight $t_1t_2$.

\sex
\leavevmode
\begin{enumerate}
\item
Nakajima's construction of quantum affine algebra associated to $Q$ uses $\DD=\bbG_m$, the diagonal torus in $\bbG_m^2$  (see \cite[(2.7.1), (2.7.2)]{Nak99}). 
Here the $\DD$-weight of the symplectic form on $T^* \Rep_Q(V)$ and on
$\fg$ is 2, and the condition on $\bm$ is $\bm(h)+\bm(h^*)=2$. If there are $a$ arrows
in $Q$ from vertex $i$ to $j$, we fix a numbering $h_1, \cdots,  h_{a}$ of these arrows, and let 
\[
\bm(h_p):=a + 2 - 2p, \,\ \bm(h_p^*):=-a  + 2p, \,\ \text{for}\,\  p=1, \cdots, a. 
\]
\item
In \cite{SV} the elliptic Hall algebra (the spherical double affine Hecke algebra of $\GL_{\infty}$) is obtained from the choice 
$\DD=\bbG_m^2$ and  $\bm=1$. \end{enumerate}

\sss{
The extension correspondence for quivers
}
\lab{The extension correspondence for quivers}
\label{subsec:QuiverCorresp}
The moduli 
$\RR=Rep_Q$ is given by pairs of $V\in\VVI$ and $a\in\Rep_Q(V)$.
We denote the elements of $\RR^m$ as
sequences $(V_\bu,a_\bu)$ of pairs  of $(V_i,a_i)\in\RR$.

Let $\FF^m$  be the moduli of $m$-step filtrations $F
=
(0=F^0\sub F^1\sub\cddd\sub F^m=V)$ on objects $V$ of $\VVI$.
Similarly, 
we consider the moduli of filtrations
$\FF^m\RR$ of representations, the objects
are triples of $V\in\VVI$,
representation $a$ of $Q$ on $V$ and
a compatible filtration $F\in\FF^m\RR(V)$.
We 
denote the fiber of $\FF^m$ at $V\in\VVI$  by
$\FF^m(V)$
and the fiber of $\FF^mRep_Q$
at
$F\in\FF^m(V)$
by $\FF^m\RR(F)=Rep_Q(F)$.

The fiber $\RR(V)$
of $\RR$ at $V\in\VVI$
is $\Rep_Q(V)$. Also, $\Rep_\Qs(V)=\RR(V)^*$
and for $\bRR=Rep_\bQ$ we have 
$\bRR(V)=T^*\RR(V)$.
By representations
on a sequence $V_\bu=(V_i)_{k=1}^m\in\VV^m$, we mean a sequence
of representations, say
$\RR(V_\bu)=
\pl_{k=1}^m\ \RR(V_k)$.

A decomposition $\bold{f}$ of $v\in\N[I]$ as $\bv_1+\cddd+\bv_m$
gives the connected components
$\FF^{\bold{f}}$ and $\FF^m\RR(\bold{f})$,  given by $\dim[Gr_F(V)]=\bold{f}$.  
The stabilizer $P$ of a chosen  $F\in\FF^\bold{f}(V)$ is a parabolic
in $G$ and then $\FF^{\bold{f}}\cong P\bs G$.

Now, the $m$-step {\em extension correspondence} for $\RR$ is
$$
\RR^m
\laa{p}
\FF^m\RR
\raa{q}
\RR
$$
where $p(V,a,F)=Gr_F(V, a)$ and
$q(V,a,F)=(V,a)$. 
The obvious splitting $\pl_1^m$ of $p$ is given by 
sending $(V_\bu,a_\bu)$ to $(\pl_1^m\ V_i,\pl_1^ma_i,F)$
for $F_p=\pl_{k=1}^p V_k$.

A filtration $F$ on  vector spaces $A,B$
defines
a filtration on $\Hom(A,B)$, 
where  operator $x$ is in $F_d$ if
$xF_pA\sub F_{p+d}B$ for all $p$.
In particular we get a filtration $F_d(A^*)=F_{-d-1}^\perp$
and the two filtrations on $\Hom(B,A)\cong Hom(A,B)^*$ coincide.

So, a filtration
$F\in\FF^m(V)$ induces a filtration on $\Rep_Q(V)\sub\End(\pl_{i\in I} V^i)$
with  $x\in F_d\Rep(V)$ if
$xF^pV^{h'}\sub F^{p+d}V^{h''}$.
Then $F_0\Rep_Q(V)$ is the  space $\Rep_Q(F)$
of
representations compatible with $F$  
and
$Gr_0^F\Rep_Q(V)=\Rep(Gr_FV)$.
Also, the two filtrations on $Rep_\bQ(V)=T^*\Rep_Q(V)$ coincide.

\sus{
Thom line bundles 
}
\lab{Thom line bundles Th(V) of vector bundles V}

\sss{
Classical Thom bundles for quivers
}
\lab{Classical Thom bundles for quivers}
As $\Rep_Q(V)$ is quadratic in $V$ we define its bilinear  version
$\Rep_Q(V_1,V_2)=\pl_{h\in H}\ \Hom(V_{1}^{h'},V_{2}^{h''})
$ for $V_i\in\VVI$.
Let $v_i=\dim(V_i)$ and denote $L=GL(V_1)\tim GL(V_2)$. 
Over
$\fA_L
=\G\syp{v_1}\tim \G\syp{v_2}$
we
define the line bundle
$$
\LL(Q,A)_{v_1,v_2}\dff\ \Th_{L}(\Rep_Q(V_1,V_2))
.$$

\slem
(a)
For a  quiver $Q=(I,H)$ the Thom line bundle
$\sL(Q,A)\dff \Th[\Rep_Q]$ is a local line bundle $\OO(-\QQ)$ on
$\fA(\Rep_Q)=\HH_{\G\tim I}$,
corresponding to the incidence quadratic form  $\QQ$ of the quiver.

(b)
The line bundle
$\LL(Q,A)_{V_1,V_2}$
on
$\G\syp{v_1}\tim \G\syp{v_2}$
is bilinear in $V_1,V_2$ in the sense that
for the addition map $\Ss:GL(U')\tim GL(U'')\inj GL(U'\pl U'')$ one has
$$
\LL_{U',V}\bten \LL_{U'',V}
\ \cong\
(\fA_{\Ss})^*(\LL_{U'\pl U'',V})
$$ 
and the same for $V$.


\pf
(a)
For $V\in \VVI$ and $G=GL(V)$,
as a $G$-module
$\RR(V)= \Rep_Q(V)$ is $\pl_H\ \RR_h(V)$
for
$\RR_h(V)=\Hom(V^{h'},V^{h''})$.
The corresponding  connected component of $\Rep_Q$ is
a vector bundle $G\bs \Rep_Q(V)$ over $\B(G)$.
Then the ideal  $\Th_G(\RR_h)$ in $\calO_{\mathfrak{A}_G}$
is generated by the function 
$\bl_{ch(\RR_h)}$ (defined in \S\ref{sec:function l}) corresponding to the character of $\RR(h)$.

A system of coordinates $x^i_s$ on each $V^i,\ \ii$,
gives a
Cartan $T$ in $G$ such that
a basis in $X^*(T)$ can be denoted by $x^i_s$.
If $i\raa{h}j$ then the character of $\RR_h$ is
$ch(\RR_h)=\ \sum_s\ \sum_t\ x^j_t(x^i_s)\inv$, hence
$\bl_{ch(\RR_h)}
=\ \prod_{s,t}\ \bl(x^j_t(x^i_s)\inv)
$
and
the same divisor is
given by
$\prod_{s,t}\
\bl(x^j_t)-\bl(x^i_s)
$ which
is the equation of the $(i,j)$-diagonal in $\fA_T\cong \prod_\ii\ \G^{\dim(V_i)}$.

(b)
$\Hom(U,V)$ is bilinear in $U$ and $V$. 
We use the obvious
observation that if 
$V_i$ is a module for $G_i$ for $1\le i\le n$,
then
$\Th_{\prod\ G_i}(\pl\ V_i)\ \cong
\bten\ \Th_{G_i}(V_i)$.
By multiplicativity of $\Th$ this reduces to
the claim that for
a representation $V$ of $G$,
$\Th_{G\tim G'}(V\bten\k)=\Th_G(V)\bten \OO_{\fA_{G'}}$.

For this 
we can assume that $G,G'$ are reductive
and then they can be replaced by their Cartans  $T,T'$.
Then we can also assume that $V_i$ are characters $\chi$ of $T$.
But then
$\Ker(\chi\bten\k)=\Ker(\chi)\tim T'$, and this implies
the claim. 
\qed

\sus{
Cotangent versions of the extension diagram 
}
\lab{The cotangent versions of the extension diagram}

\sss{
The (co)tangent functoriality 
}
\lab{The (co)tangent functoriality}

The {\em tangent complex} of a map of smooth spaces $f:\XX\to\YY$
is
$T(f)
=\
[T\XX
\raa{df}f^*T\YY]_{-1,0}
$
on $\XX$ and
the dual
{\em cotangent complex}
is
$T^*(f)=\ [f^*T^*\YY
\raa{d^*f}
T^*\XX
]_{0,1}
$.
When $f$ is an embedding these are the (co)normal bundles
$T(f)\cong N(f)$ and $T^*(f)=T^*_\XX\YY=N(f)^*$.
The {\em Thom line bundle of a map}
$f$ is
$\Th(f)\dff\Th[T(f)]=\
\Th(f^*T\YY)\Th(T\XX)\inv$.
For the map $f$ the direct image of $A$-cohomology 
takes the form of $f_*:\Th(f)\to A(\YY)$ \cite{GKV95}.

The {\em cotangent functoriality}
 associates to   
$f:\XX\to \YY$ the correspondence
\[
T^*\YY
\laa{\tii f}
f^*T^*\YY
\raa{d^* f}
T^*\XX.\]
Therefore, any correspondence 
$
A
\laa{p} 
C
\raa{q}
B$ of smooth spaces 
gives
two cotangent correspondences
$
T^*A
\laa{\tii p}
p^*T^*A
\raa{d^*p}
T^*C
\laa{d^*q}
q^*T^*B
\raa{\tii q}
T^*B$
that compose to the correspondence $p^*T^*A\tim_{T^*C}q^*T^*B$. \fttt{The fibered product has to be derived for the relevant base change to hold unless $d^*p, d^*q$ are transversal. } Say,  in the category of schemes this fibered product consists of all
$c\in C,\ \al\in T^*_{p(c)}A,\ \be \in T^*_{q(c)}B$
such that
$d^*p\al=d^*q\be$,
so
by passing to $(c,\al,-\be)$ we identify it with
$T^*_C(A\tim B)$.
Then  the cotangent version of the original correspondence is 
$T^*A
\laa{\dop}
T^*_C(A\tim B)
\raa{\doq}
T^*B$.

\sss{
Stacks
}
\lab{Stacks}
If $X$ is a smooth variety
with an action of a group $G$ then  $G\bs X$
is a smooth stack
whose tangent complex is $[\fg\to TX]_{-1,0}$
and the cotangent complex is $[T^*X\to\fg^*]_{0,1}$.

We will consider 
a map of smooth varieties $X_1\raa{f}X_2$ 
and $G_1\to G_2$ a  compatible  map of groups $G_i$ acting on $X_i$.
Then 
for $\XX_i=\GXi$
one gets $\ F:\XX_1\to \XX_2$.
We will calculate
its cotangent correspondence $T^*\XX_2\laa{\tii F}F^*T^*\XX_2\raa{d^*F}T^*\XX_1$.
First, the  Thom line bundles for stacky versions are the equivariant Thom
bundles for $f$ plus a change of equivariance
factor
$\Th_{G_1}(\fg_1/\fg_2)$ defined as
$\Th_{G_1}([\fg_1\to\fg_2]_{0,1})	$.

\slem
(a)
$\Th(F)=\Th_{G_1}(f)\ \ten\ \Th_{G_1}(\fg_1/\fg_2)	$.

(b)
$\Th(d^*F)\ =\
\Th_{G_1}(d^*f)
\ \ten\ \Th_{G_1}(\fg_1/\fg_2)
$.

(c)
$\Th(\tii F)\ =\ \Th_{G_1}(f)\ \ten\ \Th_{G_1}(\fg_1/\fg_2)$.

(d) 
The pull back map on cohomology
$A(T^*\XX_2)\raa{(\tii F)^*} A(F^*T^*\XX_2)$
is the
same as
$A(\XX_2)\raa{F^*} A(\XX_1)$.
Also, 
$(d^*F)^*$ is identity on
$A(\XX_1)$.

\pf
The (co)tangent complexes of spaces $\XX_i$ are calculated by
formulas $T(G\bs X)=G\bs (T_{G}X)$
and $T^*(G\bs X)=G\bs (T^*_{G}X)$,
where  $T_{G}X=
[\fg\tim X
\to TX]_{-1,0}$
and
$T^*_{G}X=
[
T^*X\to
\fg^*\tim X]_{0,1}$.

The map $F$ gives pull backs
$F^*(T\XX_2)
=\
\fra{X_1}{G_1}\tim_{\fra{X_2}{G_2}}\fra{T_{G_2}X_2}{G_2}
=\
\fra{f^*T_{G_2}X_2}{G_1}
$
and
$
F^*(T^*\XX_2)
=\
\fra{f^*T^*_{G_2}X_2}{G_1}
$.

(a)
The tangent complex $T(F)=\
[T\XX_1\to F^*T\XX_2]_{-1, 0}$
of the map $F$, 
comes from
(the $G_1$-quotient of)
the map of complexes
$
[\fg_1\to TX_1]
\to
[\fg_2\to TX_2]
$
given by $TX_1\raa{df}f^*TX_2$
and $\fg_1\to \fg_2$.
When we view this as
a bicomplex with horizontal and vertical degrees in $[-1,0]$,
then  $T(F)$ is its total complex
$[\fg_1
\to
TX_1
\pl \fg_2\to
f^*TX_2]_{-2,0}$, which   is an extension of
complexes
$[TX_1\raa{df} f^*TX_2]_{-1,0}=T(f)$
and
$[\fg_1\to\fg_2]_{-2,-1}$.
So
$\Th(F)$ is as stated.

Now, 
the cotangent correspondence can be written as 
$$
\CD
T^*\XX_2
@<\tii F<<
\ F^*(T^*\XX_2)
@>d^*F>>
T^*\XX_1
\\
@V=VV
@V=VV
@V=VV
\\
\fra{T^*_{G_2}X_2}
{G_2}
@<\
<<
\fra{f^*T^*_{G_2}X_2}
{G_1}
@>
>>
\fra{T^*_{G_1}X_1}
{G_1}
.\endCD
$$

(b)
Write $d^*F$ as $\cdF/G_1$
where
$\cdF:\
f^*[T^*X_2\to \fg^*_2]
\ \to\
[T^*X_1\to \fg_1^*]
$
is a map of complexes
viewed
as a bicomplex with all horizontal and vertical degrees in $[-1,0]$. 
So, $T(\cdF)$ is the total complex
$[f^*T^*X_2
\ \to\
T^*X_1\pl
\fg^*_2
\ \to\
\fg_1^*]_{-2,0}
$
which is an extension of complexes
$[f^*T^*X_2\
\raa{d^*f}
T^*X_1]_{-2,-1}$ and 
$[\fg_2^*\to \fg_1^*]_{-1,0}$.
So, 
$\Th(d^*F)=\Th_{G_1}(\cdF)=
\Th_{G_1}(d^*f)\ten
\Th_{G_1}([\fg_2^*\to \fg_1^*]_{-1,0})
=\Th_{G_1}(d^*f)\ten\Th_{G_1}([\fg_1\to \fg_2]_{0,1})
$
and then we use  invariance of the Thom line bundle under duality.

(c)
Denote
the complex
$T^*_{G_2}X_2=
[T^*X_2\to\fg_2^*]_{0,1}
$ by $\VV$ and let $\eta:f^*\VV\to\VV$, then 
the map
$\tii F$ 
is given by $\eta$ and the change of symmetry $G_1\to G_2$.
So, part (a) says that
$\Th(\tii F)=\Th_{G_1}(\eta)\ten \Th_{G_1}([\fg_1\to\fg_2]_{0,1})$.
Let us denote
$\pi:\VV
\to X_2$
and $\pi:f^*\VV
\to X_1$, then  
$T(\eta)=[T(f^*\VV)\raa{d\eta}\eta^*T(\VV)]_{-1,0}$
is 
$\pi^*T(f)$.
(One has $0\to\pi^*\VV\to T\VV\to \pi^*TX_1\to 0$
and $0\to\pi^*f^*\VV\to Tf^*\VV\to \pi^*TX_2\to 0$.  Now the map of complexes
is identity on subsheaves
$\pi^*f^*\VV\cong \eta^*\pi^*\VV$
and what remains is $\pi^*TX_2\to
\eta^*\pi^*TX_1=\pi^*f^*TX_1 $.)

So,
$\Th_{G_1}(\eta)=\Th_{G_1}(\pi^*T(f))$, 
and since $\pi$ is contractible this is $\Th_{G_1}(f)$.

(d)
After contracting complexes of vector bundles
the maps $\tii F$ and $d^*F$ become respectively
the map $X_2/G_2\laa{F}X_1/G_1$ and the identity on $X_1/G_1$.
\qed

\sus{
The $A$-Cohomology of the cotangent correspondence
for extensions
}
\lab{The A-Cohomology of the extension diagram}

We recall the  construction of \cite{YZ1} of a quantum group in
the above set up.
It originated from the study of affine quantum groups in \cite{Nak99} and \cite{SV}, and is closely related to \cite{KS}. 

\sss{
Connected components of the cotangent correspondence
}
Fixing $V\in\VVI$ and $F\in \FF^m(V)$, let $\Gr_F(V)=\oplus_{k=1}^m V_k$. 
We denote $G=\GL(V)$, $L=\prod_{k=1}^m\GL(V_k)$ the automorphism group of $Gr_F(V)$, and $P$ is a parabolic subgroup of $G$ with a Levi subgroup $L$.
Let $U$ be the unipotent radical of $P$. 
Denote the Lie algebras by $\fp,\fl,\fu$. 

These choices fix the connected component of
the correspondence
$\RR^m\laa{p}\FR\raa{q} \RR\ \to\ \VVI$ given by
\begin{equation}\label{eqn:QuiverV}
L\bs \RR(\Gr_F(V))
\laa{p}
P\bs\RR(F)
\raa{q}
G\bs \RR(V)
\ \to\
G\bs\pt
.\end{equation}

\sss{
Line bundles from the cotangent correspondence}
\lab{The multiplication induced by the cotangent version of the correspondence}

The extension correspondence  gives two cotangent correspondences
\begin{equation}\label{eqn:cotangent}
T^*\RR^m
\laa{\tii p}
p^*T^*\RR^m
\raa{d^*p}
T^*\FR
\laa{d^* q}
q^* T^*\RR
\raa{\tii q}
T^*\RR
.\end{equation}
These compose to a single correspondence
as in \ref{The (co)tangent functoriality}
which is the {\em cotangent correspondence}
of the extension correspondence.
We will not consider it since we are calculating here 
its effect on cohomology and
this is the composition of effects of
the above two simpler correspondences.

Let $V\in \VVI$ and $F\in \FF^m(V)$. 
We write the fiber of the correspondence \eqref{eqn:QuiverV}
over $F$ as 
$$\RR(\Gr_F(V))
\laa{\cp}
\RR(F)
\raa{\cq}
\RR(V).
$$ 
The connected component of the diagram \eqref{eqn:cotangent}  determined by $F$ takes the form 
\begin{equation}\label{eqn:cotangV}
T^*(L\bs\RR^m(Gr F))
\laa{\tii p}
p^*T^*(L\bs\RR^m(Gr F))
\raa{d^*p}
T^*(P\bs\FR(F))
\laa{d^* q}
q^* T^*(G\bs \RR(V))
\raa{\tii q}
T^*(G\bs \RR(V))
.\end{equation}

\slem
With notations as above (and the filtration on $\RR(F)$ as  in
\S~\ref{subsec:QuiverCorresp}) we have
$$
\Th(d^*p)\cong 
\Theta_{L}(\fg/\fp)\otimes\Theta_{L}(F_{-1}\RR(F))
\aand
\Th(\tii q)\cong \Th_{L}[\RR(V)/\RR(F)]\otimes\Th_{L}(\fg/\fp)\inv
.$$

\pf
According to the lemma \ref{Stacks}.b
$\Th(d^*p)
$ is $
\Th_{L}(d^*\cp)
\ten\Th_{L}([\fp\to\fl]_{0,1})
$.
The second factor is $\Th_{L}(\fu)$, since $\fu\cong (\fg/\fp)^*$
we can write it as
$\Th_{L}(\fg/\fp)$.
For the first factor, as $\cp:\RR(F)\to \RR(\Gr_F(V))$
we get $d^*\cp:
\RR(F)\tim \RR(\Gr_F(V))^*\to
\RR(F)\to\RR(F)^*
$, so up to a factor
$\RR(F)$ this is $\RR(\Gr_F(V))^*\inj\RR(F)^*$
with the quotient $[F_{-1}\RR(F)]^*$.
So, the first factor is
$\Th_{L}([F_{-1}\RR(F)]^*)=\
\Th_{L}(F_{-1}\RR(F))$.

\noi
Again,
by  lemma \ref{Stacks}.c\
$\Th(\tii q)$ is
\ $\Th_{P}(\cq)\ten\Th_{P}([\fp \to\fg ]_{0,1})
$\ for
the
embedding $[\RR(F)\raa{\cq} \RR(V)]_{-1,0}$. So, the 
first factor is $\Th_{L}[\RR(V)/\RR(F)]$
and the second is $\Th_{L}(\fg/\fp)\inv$.

\qed

\sss{
Dilations
}
\label{sss:dilation}
Recall the action of the dilation torus $\DD\subseteq \Gmt$ from \S~\ref{sss:dilation}. 
The  weight
of the first $\Gm$-factor
on  $\RR$
is prescribed by $\bm$ while the second factor acts trivially.
Then the  $\DD$-action on $T^*\RR$ is uniquely determined
by asking that the natural symplectic form on $T^*\RR$ has weight $t_1t_2$.
We  denote the $\DD$-character of weight $t_1t_2$ by $\omega$,
so that the  
 $\DD$-action on $T^*\RR$ is twisted by $\omega$. This gives rise to the following twisted version of \eqref{eqn:cotangV}, 
\begin{align}\label{eq:twist}
T^*(L\bs\RR^m(Gr F))\otimes \omega &
\laa{\tii p}
p^*T^*(L\bs\RR^m(Gr F))\otimes \omega
\raa{d^*p}
T^*(P\bs\FR(F))\otimes \omega\\&
\laa{d^* q}
q^* T^*(G\bs \RR(V))\otimes \omega
\raa{\tii q}
T^*(G\bs \RR(V))\otimes \omega \notag
.\end{align}
The maps in the above diagram are equivariant with respect to $\DD$. 

Now we analyze  $\DD$-action on the relative tangent complexes of $F$, $d^*F$, and $\widetilde F$.
Lemma~\ref{Stacks} applies to induced actions on cotangent bundles. 
When working $\DD$-equivariantly we need to add an $\omega$-twist.
This applies to the  Lie algebra factors in Lemma~\ref{Stacks}.
that come from the cotangent complexes. 
On the other hand, the Lie algebra factors
that come from the change of symmetry are not affected
as they only carry the adjoint action. 
To simplify notations, for any group $H$, we denote $H\times\DD$ by $\widetilde H$.
Then 
\[
\Theta_{\DD}(d^*F)\cong \
\Th_{\widetilde{G}_1}(d^*f)
\ \ten\ \Th_{\widetilde{G}_1}((\fg_1/\fg_2)\otimes\omega)
\aand
\Th_{\DD}(F)=\Th_{\widetilde{G}_1}(f)\ \ten\ \Th_{\widetilde{G}_1}(\fg_1/\fg_2)
.
\]
Therefore, 
$\Th_{\DD}(\tii F)\ =\ \Th_{\widetilde{G}_1}(f)\ \ten\ \Th_{\widetilde{G}_1}(\fg_1/\fg_2)$.

\slem
With notations above:\
$\Th(d^*p)\cong 
\Theta_{\widetilde{L}}(\fg/\fp\otimes\omega)\otimes\Theta_{\widetilde{L}}(F_{-1}\RR(F))$
\ and \
$\Th(\tii q)\cong \Th_{\widetilde{L}}[\RR(V)/\RR(F)]\otimes\Th_{\widetilde{L}}(\fg/\fp)\inv$. 

\pf
According to the lemma \ref{Stacks}.b
$\Th(d^*p)
$ is $
\Th_{\widetilde{L}}(d^*\cp)
\ten\Th_{\widetilde{L}}([\fp\to\fl]_{0,1}\otimes\omega)
$.
The second factor is $\Th_{\widetilde{L}}(\fu\otimes\omega)$, since $\fu\cong (\fg/\fp)^*$
we can write it as
$\Th_{\widetilde{L}}(\fg/\fp\otimes\omega)$.

Again, by  lemma \ref{Stacks}.c, 
$\Th(\tii q)$ is
$\Th_{\widetilde{P}}(\cq)\ten\Th_{\widetilde{P}}([\fp \to\fg ]_{0,1})$ for
the linear embedding $[\RR(F)\raa{\cq} \RR(V)]_{-1,0}$. So, the 
first factor is $\Th_{\widetilde{L}}[\RR(V)/\RR(F)]$
and the second is $\Th_{\widetilde{L}}(\fg/\fp)\inv$, as it comes from the change of symmetries. 

\qed

\sus{
$\DD$-quantization of the monoid $(\HGI,+)$
}
\lab{DD-quantization of the monoid HGI}
Here we recall  the construction from \cite{YZell} of a deformation $(\Coh(\HGI), \star)$ of the convolution
on the monoid $(\HGI,+)$.
The quantum group $U_{\DD}(Q,A)$ and its positive part $U^+_{\DD}(Q,A)$ were constructed in \cite{YZ1, YZ2}, 
as algebra objects in $(\Coh(\HGI),\star)$, and hence in particular as $R$-algebras. 

\sss{
Local and biextension line bundles $\sL_\DD(Q,A)$ and $\LL_\DD(Q,A)$
}
\label{sec:local and biextension}

These will be upgrades of  $\sL(Q,A)$
and $\calL(Q,A)$ 
from
\S~\ref{Classical Thom bundles for quivers}.
They will be constructed as special cases of
line bundles associated to cotangent correspondences of extension moduli \eqref{eqn:cotangV}. 
\begin{description}
\item[Case 1] 
The {\em biextension} line bundle $\LL=\LL_\DD(Q,A)$
comes from $m=2$\ie the
$2$-step filtrations $\mathcal{F}^2(V)$ of $V$.
For  $\Gr_F(V)=V_1\oplus V_2$
\[
\calL_{V_1, V_2}:=\Th(d^*p)\otimes \Th(\tii q)
\]
is a line bundle on $
\fA_{\widetilde{L}}\cong \fA_{G(V_1)}\times\fA_{G(V_2)} \times \fA_{\DD}$.

\item[Case 2] Our  quantum version $\sL=\sL_\DD(Q,A)$
of the local line bundle
$\sL(Q,A)$ depends on a choice of a {\em type}
of a complete flag $F\in\FF^m(V)$ which is
$\bold{f}=\dim(\Gr_{F}(V))\in (\bbN^I)^m$.
Then  
\[
[\sL_\DD(Q,A)]_{V, \bf{f}}\ \dff\
\Th(d^*p)\otimes \Th(\tii q)
\]
is a line bundle on $\fA_{\widetilde{L}}=\G^{|V|} \times \fA_{\DD}$, where $|V|=\sum_{i\in I} \dim(V^i)$ (here  the Levi subgroup $L$ is a Cartan in
$GL(V)$). It is called  {\em the local line bundle}. 
\end{description}

One easily sees that the restrictions of  "quantum objects"
$\sL_\DD(Q,A)$ and $\LL_\DD(Q,A)$
to $0\in\fA_\DD$
are
the classical Thom line bundles
$\sL(Q,A)$ and $\calL(Q,A)$
from \S\ref{Classical Thom bundles for quivers}.

\sss{
Convolutions and biextensions
}
We recall the monoidal structure $\star$ on 
coherent sheaves on $\HHH_{\G\tim I}\tim\fA_{\DD}$
(over the base scheme $\fA_{\DD}$) from \cite{YZell}.

For a smooth curve $C$, $\HHCI$
is a commutative monoid
freely generated by $C$. The
operation $\Ss:\HH_\CI\tim\HH_\CI\to \HH_\CI$  is the  addition of divisors
 (``symmetrization'').
Since it is a finite map it defines a convolution operation
on the abelian category $\Coh(\HHCI)$
of coherent sheaves
by $\FF\ast\GG=\Ss_*(\FF\bten\GG)$.

The following definition can be found in \cite[P126]{P}. 
A line bundle $\calL$ over $(\HHCI)^2$ is
a {\em biextension} (or {\em Poincare} line bundle)
if we have the following isomorphism of line bundles
\begin{align*}
& a_{x_1, x_2; y}: \calL_{x_1+x_2, y} \to \calL_{x_1, y}\otimes  \calL_{x_2, y}, \\
& a_{x; y_1, y_2}: \calL_{x, y_1+ y_2} \to \calL_{x, y_1}\otimes  \calL_{x, y_2}, 
\end{align*}
which satisfy the following cocycle conditions
\begin{align*}
& \text{(i)} \,\  a_{x_1+x_2, x_3; y} \circ ( a_{x_1, x_2; y} \otimes \id)=a_{x_1, x_2+x_3; y}\circ (\id \otimes a_{x_2, x_3; y}), \\
& \text{(ii)} \,\  a_{x; y_1+y_2, y_3} \circ ( a_{x; y_1, y_2} \otimes \id)=a_{x; y_1, y_2+y_3}\circ (\id \otimes a_{x; y_2, y_3}),\\
& \text{(iii)} \,\ (a_{x_1, x_2; y_1}\otimes a_{x_1, x_2; y_2}) a_{x_1+x_2; y_1, y_2}=
(a_{x_1; y_1, y_2}\otimes a_{x_2; y_1, y_2}) a_{x_1, x_2; y_1+y_2}. 
\end{align*}
This is equivalent to a central extension of the monoid
$(\HHCI,+)$ (or its  group completion) by $\Gm$.

Now $\LL$  twists the convolution on
$\Coh(\HHCI)$
to another monoidal structure
$\FF\sta\GG\dff\ \Ss_*[(\FF\bten\GG)\ten\LL]$.

From now on the curve $C$ will be  $\G=\fA_{\Gm}$.

\slem
The line bundle $\LL=\LL_\DD(Q,A)$ on $(\HHCI)^2\times \fA_{\DD}$ defined
in \S \ref{sec:local and biextension}.1
is an
$\fA_{\DD}$-family of biextension line bundles. This gives a
``$\DD$-twisted'' convolution
on $\Coh(\HGI\tim\fA_{\DD})$ 
by
$$
\FF\sta\GG\dff\ \Ss_*[(\FF\bten_{\fA_{\DD}}\GG)\ten\LL]. 
$$

\pf
We need to check that the quantum version of   $\calL$ is still
a biextension.
Notice that
the quantum version, has an extra factor
$\Theta_{\DD}(\fg/\fp)$.
However since 
for $m=2$, the space $\fg/\fp$ is of the form
$ \Hom(V_1, V_2)$,
the argument in the proof of  Lemma \ref{Classical Thom bundles for quivers}.b 
applies again.
\qed

\spro 
\cite[Theorem A, Theorem 3.1]{YZell}

(a) $(\Coh(\HGI\tim\fA_{\DD}), \star)$ is a monoidal category with a meromorphic braiding which is symmetric. 
The unit is the structure sheaf on
$\HHH^0_{\bbG\times I}\tim\fA_{\DD}$.

(b)
The structure sheaf on
$\HGI\tim\fA_{\DD}$ is an algebra object in this category.

For any $\tau\in \fA_\DD$, we denote by $\calL_\tau$ the restriction of $\LL$ to $\tau\in \fA_{\DD}$, and 
$\FF\sta_\tau\GG:=\ \Ss_*[(\FF\bten_{\fA_{\DD}} \GG)\ten \LL_\tau]$.

\srem
One way to motivate the $\LL$-twisted convolution of
coherent sheaves on $(\HGI)^2\tim\fA_{\DD}$ is to  notice that
when the cohomology theory $A$ extends to constructible sheaves,
then for a constructible
$\FF$ on a space $X$,  the cohomology $A(\FF)$ is a coherent sheaf on $\fA(X)$.
In this case the $A$-cohomology functor intertwines the
convolution of constructible sheaves on $\Rep_Q$ and the 
$\LL$-twisted convolution of
coherent sheaves on $\fA(\Rep_Q)=(\HGI)^2\tim\fA_{\DD}$. (This
follows as in the proof of lemma
\ref{The multiplication induced by the cotangent version of the correspondence}.)

\sss{
Quantum groups $U^+_{\DD}(Q,A)\subset U_{\DD}(Q,A)$
}\label{sus:quanLocLine}

Now we consider the set up of \S~\ref{The multiplication induced by the cotangent version of the correspondence} 
with $V\in \VVI$ and $F\in\FF^m(V)$. 
Let $\bold{f}=\dim(\Gr_{F}(V))\in (\bbN^I)^m$ be the type of the filtration $F$. 
Applying the cohomology theory $A$ to the diagram \eqref{eqn:cotangV},
 we have the following multiplication map associated to $\bf{f}$: 
\begin{equation}\label{eqn:multip}
m_{\bf{f}}:=(\tilde{q}_*)\circ(d^*q^*)\circ(d^*p_*)\circ(\tilde{p}^*):
\mathbb{S}_*( \Th(\tilde{q})\otimes \Th(d^*p))
\to 
A_\DD(T^*(G\bs \RR(V)))\cong\calO_{\fA_G\tim\fA_{\DD}}, 
\end{equation}
where $\mathbb{S}: \fA_{L}\to  \fA_{G}$ is the symmetrization map. 

Let $Sph(V)$ be the set of types $v$ of filtrations  in  $\FF^m(V)$ consisting of complete flags
(so $m=|v|\dff \sum_{i\in I}v^i$). 
We define $U^+_{\DD}(Q,A)$ so that on the connected component $\fA_G\tim\fA_{\DD}$, 
$$
(U^+_{\DD}(Q,A))_V\dff\ \sum_{\bold{f}\in Sph(V)} \hbox{Image} (m_{\bf{f}})
\ \
\ \sub\
\
A_{\tG}=\calO_{\fA_G\tim\fA_{\DD}}
.$$

\slem
\begin{enumerate}
\item 
The coherent sheaf
$U^+_{\DD}(Q,A)$ is an ideal sheaf on $\fA_G\tim\fA_{\DD}$. 
\item $U^+_{\DD}(Q,A)$ is an $\fA_{\DD}$-family of  algebras in
the monoidal categories $(\Coh(\HGI),\star_\tau)$. 
\end{enumerate}

\begin{proof}
As each $m_{\bold{f}}$, for $\bold{f}\in Sph(V)$,   is a morphism of coherent sheaves, the image $\hbox{Image} (m_\bold{f})$ is a coherent subsheaf  in 
$\calO_{\fA_G\tim\fA_{\DD}}$. 
Since  $Sph(V)$ is a finite set, $(U^+_{\DD}(Q,A))_V$ is a sum of finitely many coherent subsheaves,  
so it is itself a coherent subsheaf of $\calO_{\fA_G\tim\fA_{\DD}}$. A coherent subsheaf  of the structure sheaf is a sheaf of ideals, hence so is $(U^+_{\DD}(Q,A))_V$.

For (2), the algebra structure on $(U^+_{\DD}(Q,A))_V$ is defined using $m_{\bf{f}}$, where $F$ is the 2-step filtrations in \S\ref{sec:local and biextension} \textbf{Case 1}. 
\end{proof}

The sheaf $U^+_{\DD}(Q,A)$ on $\fA_{G}\times \fA_{\DD}$ is denoted by $\calP^{sph}$ in \cite{YZ1}, since it is the spherical subalgebra of the cohomological Hall algebra of preprojective algebra.

The affine quantum group $U_{\DD}(Q,A)$ associated to the quiver $Q$ and the cohomology theory $A$ is defined in \cite{YZ2} as the Drinfeld double of $U_{\DD}^+(Q,A)$. 
The quantization parameters of $U_{\DD}(Q,A)$ are given by $\fA_{\DD}$. 
This Drinfeld double was constructed in \cite{YZ2} using a comultiplication and 
a bialgebra pairing on an extended version of $U_{\DD}^+(Q,A)$. $U_{\DD}^+(Q,A)$ itself also has a coproduct but 
in the meromorphic braided tensor category $(\Coh(\HGI), \star)$ \cite{YZell}. The affine quantum group $U_{\DD}(Q,A)$ acts on the corresponding $A$-homology of the Nakajima quiver varieties (see \cite{YZ1, YZ2}),  generalizing a construction of 
Nakajima \cite{Nak99}.


\sss{
Local line bundle from zastava
and $U^+_{\DD}(Q,A)$
}
\label{sss:modify}
In this section we
associate to a quiver $Q$ 
a different integral form     $\underline{U}^+_{\DD}(Q,A)$
of the  quantum group $U^+_{\DD}(Q,A)$. 
There is an algebra homomorphism $U^+_{\DD}(Q,A)\to \underline{U}^+_{\DD}(Q,A)$, which becomes an isomorphism after a certain localization to be explained below. 

Recall from \S\ref{subsec:o1} that for a semisimple simply laced group  of simply connected type,
the restriction from the loop Grassmannian
gives
a local line bundle
on $\HHCI$
related to the
Cartan quadratic form. 
On the other hand, a quiver $Q$
produces the local line bundle
$\sL(Q,A)$ 
corresponds to the incidence quadratic form
of $Q$ (\S~\ref{Classical Thom bundles for quivers}).
So, when $Q$ is the corresponding Dynkin graph 
the quadratic forms differ on the diagonal.
Then the  new quantum group 
$\underline{U}^+_{\DD}(Q,A)$
will be related to the Cartan quadratic form.



\medskip

To define
$\underline{U}^+_{\DD}(Q,A)$ as an algebra object in a monoidal category,  we
modify the local and biextension line bundles from \S~\ref{sec:local and biextension}.
We follow the notations from \S~\ref{sss:dilation}. 
Let  $V\in \VVI$ and $F\in\FF^m(V)$. 
Let $\bold{f}=\dim(\Gr_{F}(V))\in (\bbN^I)^m$ be the type of the filtration $F$. Let $G$ be the automorphism group of $V$ and $P$ be the parabolic subgroup preserving $F$, and $L$ the automorphism group of $\Gr_{F}(V)$. 
Consider the $\tilde{L}$-representation $(\fu\oplus\fg/\fp)\otimes\omega$. 

Now assume that $\bf{f} \in Sph(V)$\ie it is a type of a complete flag. 
By the invariance under the symmetric group
$\Theta_{\tilde{L}}((\fu\oplus\fg/\fp)\otimes\omega)$ (as a line bundle on $\fA_L\times\fA_{\DD}$),
is obtained from pullback of a line bundle on $\fA_G\times\fA_{\DD}$, which by an abuse of notation is denoted by $\Theta_{\tilde{G}}((\fu\oplus\fg/\fp)\otimes\omega)$.
Twisting by this line bundle, we get from \eqref{eqn:multip}  the following map 
\begin{equation} 
m_{\bf{f}}:=(\tilde{q}_*)\circ(d^*q^*)\circ(d^*p_*)\circ(\tilde{p}^*):
\mathbb{S}_*( \Th(\tilde{q})\otimes \Th(d^*p)\otimes \Theta_{\tilde{L}}((\fu\oplus\fg/\fp)\otimes\omega))^{-1}
\to 
\Theta_{\tilde{G}}((\fu\oplus\fg/\fp)\otimes\omega)^{-1}.
\end{equation}
We define $(\underline{\sL}_{\DD}(Q,A))_{V,\bf{f}}$ to be  $\Th(\tilde{q})\otimes \Th(d^*p)\otimes \Theta_{\tilde{L}}((\fu\oplus\fg/\fp)\otimes\omega)^{-1}$.
We define $\underline{U}^+_{\DD}(Q,A)$ so that on the connected component
of
$\fA_G\tim\fA_{\DD}$ containing $V$, 
$$
(\underline{U}^+_{\DD}(Q,A))_V\dff\ \sum_{\bold{f}\in Sph(V)} \hbox{Image} (m_{\bf{f}})
\ \
\ \sub\
\
\Theta_{\tilde{G}}((\fu\oplus\fg/\fp)\otimes\omega)^{-1}
.$$

Similarly, we have the modified biextension line bundle $\underline{\LL}_\DD(Q,A)$
coming from  the
$2$-step filtrations $\mathcal{F}^2(V)$ of $V$.
For  $\Gr_F(V)=V_1\oplus V_2$
\[
\underline{\calL_{V_1, V_2}}:=\Theta_{\widetilde{L}}(\fg/\fp\otimes\omega)^{-1}\otimes\Theta_{\widetilde{L}}(F_{-1}\RR(F))\otimes \Th(\tii q)
\]
is a line bundle on $
\fA_{\widetilde{L}}\cong \fA_{G(V_1)}\times\fA_{G(V_2)} \times \fA_{\DD}$.

Restrictions of these  "quantum objects"
$\underline{\sL}_\DD(Q,A)$ and $\underline{\LL}_\DD(Q,A)$
to $0\in\fA_\DD$
gives line bundles on $\HHCI$ and $\HHCI^2$, denoted by 
$\underline{\sL}(Q,A)$ and $\underline{\calL}(Q,A)$ respectively.

Similar to \S~\ref{sus:quanLocLine},  $\underline{\LL}_\DD(Q,A)$ defines a family of monoidal structure on $\HGI$, denoted by  $(\Coh(\HGI),\underline{\star}_\tau)$.

\spro
Let $Q$ be a simply laced Dynkin quiver, with  Cartan quadratic form
$\underline{\QQ}$ and the   simply connected  group $G$.
\begin{enumerate}
\item  The local line bundle $\underline{\sL}(Q,A)$ on $\HHCI$ is isomorphic to
$AJ^*\OO_{\uGG(G)}(1)$. 
\item  $\underline{U}^+_{\DD}(Q,A)$ is an $\fA_{\DD}$-family of  algebras in
the monoidal categories $(\Coh(\HGI),\underline{\star}_\tau)$. 
\end{enumerate}
\pf
The proof of (1) is similar to that of Lemma~\ref{Classical Thom bundles for quivers}. 
By Lemma~\ref{Classical Thom bundles for quivers}, the local line bundle $\sL(Q,A)$ is associated to the incidence quadratic form $\QQ$ of $Q$. On the component associated to $\bf{f} \in Sph(V)$,  tensoring
$\sL(Q,A)$  with $\Theta_{G}(\fu\oplus\fg/\fp)$ we get $\underline{\sL}(Q,A)$, which is associated to  the Cartan quadratic form  $\underline{\QQ}$. On the other hand, when $G$ is the
 simply connected group whose Dynkin diagram is
$Q$ then $AJ^*\OO_{\uGG(G)}(1)$ is associated to $\underline{\QQ}$
by Proposition~\ref{subsec:o1}.

Proof of (2) is similar to  \S~\ref{sus:quanLocLine}. The algebra structure on $(\underline{U}^+_{\DD}(Q,A))_V$ is defined using $m_{\bf{f}}$, where taking $F$ to be a  2-step filtration gives $\underline{\LL}_\DD(Q,A)$. 
\qed

\srems (1) 
There is an algebra homomorphism $U^+_{\DD}(Q,A)\to \underline{U}^+_{\DD}(Q,A)$.
Topologically, for $V\in\VVI$ the map $U^+_{\DD}(Q,A))_V\to (\underline{U}^+_{\DD}(Q,A))_V$ is induced by \[z^*z_*:\Theta_{\tilde L}((\fu\oplus\fg/\fp)\otimes\omega)\to A_\DD(T^*(L\ \RR^m(\Gr F)))\] where $z$ is the zero section of the vector bundle on $T^*(L\backslash \RR^m(\Gr F))$ which is the pullback of $(\fu\oplus\fg/\fp)\otimes\omega$ form $\tilde{L}\backslash \pt$.
In particular, the map becomes an isomorphism after inverting the Euler class of $(\fu\oplus\fg/\fp)\otimes\omega$. 

(2)
The shuffle formula for $\underline{U}^+_{\DD}(Q,A)$ is similar to that of $U^+_{\DD}(Q,A)$ given in \cite[\S~1.2]{YZ2}. \footnote{Placing the
numerator of the factor $\fac_1$ in Equation \cite[(2)]{YZ1}  on the denominator to get the corresponding factor in the shuffle formula for  $\underline{U}^+_{\DD}(Q,A)$.
The homomorphism from (1) is on the level of shuffle algebras the
multiplication by the Euler class of $(\fu\oplus\fg/\fp)\otimes\omega$.}

(3) The integral form 
$\underline{U}^+_{\DD}(Q,A)$
has originally  appeared in \cite[Appendix~B]{BFN}. The isomorphism
can be seen directly by comparing the shuffle formula for $\underline{U}^+_{\DD}(Q,A)$  with the formula \cite[(3.6.3)]{FT}. 

\se{
Loop Grassmannians $\GG^P_\DD(Q,A)$ and quantum locality
}
\lab{Loop Grassmannians GG P DD(Q,A) and quantum locality}
In the preceding section 
\ref{Local line bundles from quivers}
we have attached to a quiver $Q=(I,H)$ and a cohomology theory $A$,
a local line bundle $\sL(Q,A)$
on the colored configuration space $\HGI$ of the curve $\G$
given by $A$
(lemma \ref{Classical Thom bundles for quivers}.a).
Modified as in \S~\ref{sss:modify},  in \ref{A generalization GGP(I,QQ) of loop Grassmannians of reductive groups}
the local line bundle $\underline{\sL}(Q,A)$ can be used to 
produce a ``loop Grassmannian'' $\GG(Q,A)$ over $\HGI$.

The local line bundle $\sL(Q,A)$ is closely related to  the biextension line bundle
$\LL(Q,A)$ from lemma \ref{Classical Thom bundles for quivers}.b.
In section \ref{Local line bundles from quivers}
we have also recalled the construction
of the affine quantum groups
$\underline{U}^+_\DD(Q,A)$
and used this to select the ``correct'' quantizations
$\underline{\sL}_\DD(Q,A)$ and 
$\underline{\LL}_\DD(Q,A)$ of the above line bundles,
on the basis of relation to 
this quantum group
(\ref{sec:local and biextension}).

While pieces $\underline{\sL}(Q,A)_\al$ of the classical local line bundle depend on $\al\in\NI$
parameterising  connected components of $\HGI$, 
the pieces
$\underline{\sL}_\DD(Q,A)_\bii$
of the quantum version  depends on
a choice of $\bii\in I^\N$ (\ref{sec:local and biextension}).
This really means that we are dealing with the
non-commutative 
(ordered) configuration spaces
$\CC=\CC_{\G\tim I}=
\sq\ (\G\tim I)^n$, so that
each $\al\in \NI$ is refined to all $\bii=(i_1,...,i_n)\in I^n$ with $\sum\ i_p=\al$.
The connected components given by all
refinements $\bii$ of the same $\al$ are  related by
the meromorphic braiding from
\cite{YZell}. So, the information carried by all refinements 
$\bii$ of $\al$ is (only) generically equivalent.

All together, $\GG_\DD(Q,A)$ can still be constructed by the same prescription
as in the case of $\GG(Q,A)$.
However, the local line bundle $\underline{\sL}_\DD(Q,A)$
now
lives on the larger (``non-commutative'')
configuration space $\CC_{\G\tim I}$.
The zastava space $Z_\DD(Q,A)$
over $\CC=\CC_{\G\tim I}$ is first defined  generically in $\CC$ 
where fibers are products of projective lines.
Then the singularities of the locality structure prescribes how fibers
degenerate. Finally, passing from the zastava space to loop Grassmannian
is given by the procedure of extending the
free monoid on $I$ to the free group on $I$.

All together, the key difference  in the quantum case is seen in the configuration space.
It has  more  connected  components 
(but they are related by braiding),
and the singularities of
locality structure (hence also the notion of locality)
are now the diagonals shifted by the quantum parameter.

\sss{
The ``classical'' loop Grassmannians
$\GG^P(Q,A)$
}
The choice of $A$ influences the space $\GG^P(Q,A)$ only 
through the curve
$\G$.
Whenever $\G$ is a formal group, then the orientation
$\bl$ of $A$ identifies $\G$  with the coordinatized formal disc
$d$.

However, since the loop Grassmannian $\uGG(\Gm)$
is the free commutative group indscheme generated by $d$
the group law on
$d$ given by  $A$ induces a commutative ring structure
on the loop Grassmannian $\uGG(\Gm)$.
This is the group algebra of
the group $\G$ taken in algebraic geometry.

\srem
The universal Witt ring
has the same nature,  it is 
the homology $\hH_*(\Ao,0)$  of the multiplicative monoid
$(\Ao,0,;\cd)$ in pointed spaces.
Observations of this nature have already been made in
\cite{BZ, Str, JN}.

\sus{
Quantization
shifts diagonals
}
\label{Quantization of convolution as shift of diagonals}

Any Thom line bundle is the ideal sheaf  of the corresponding Thom divisor.
While the Thom divisor corresponding to $\sL(Q,A)$ is a combination of
diagonals of $\HH=\HHdI$,
the quantization  shifts these  diagonals in the configuration
space $\CC=\CC_{\d\tim I}$.

We first examine how an added action of a torus $\DD$ affects
the Thom divisor in general
(\ref{Deformation of a Thom line bundle from an additional torus DD}),
and then we specialize this to
the local line bundle $\underline{\sL}_\DD(Q,A)$ in
\ref{m-shifted locality}.

\sss{
Deformation of a Thom divisor
from an additional   torus $\DD$
}\lab{Deformation of a Thom line bundle from an additional   torus DD}

For a representation
$E$ of a product $\tG=G\tim\DD$
we can view the line bundle
$\Th_\tG(E)$ on $\fA_G\tim\fA_\DD$
as a family of line bundles
$\Th_\tG(E)_\tau$ (for $\tau\in\fA_\DD$),
on $\fA_G\aa{\tau}\inj \fA_G\tim\fA_\DD$.
If  $E$
contains no trivial characters of a Cartan $T$, we will see that this deformation
lifts to divisors.

First, consider the case when $G$ is  a torus $T$ and
$E=\chi\bten\ze\inv$ for characters $\chi,\ze$ of $T,\DD$
(so $\chi\ne 0$). Then  for any $\tau\in\fA_\DD$,
the restriction $\Th_\tT(E)_\tau$ to $\fA_T$
is the ideal sheaf of the divisor
\[\Ker(\fA_{\chi\bten\ze\inv})
\cap [\fA_T\tim \tau]
\ =\
\fA_\chi\inv(\fA_\ze(\tau)) \subset \fA_{T}.\] 
Here $\chi: T\to \Gm$ induces the homomorphism $\fA_{\chi}: \fA_{T}\to \fA_{\Gm}=\G$ as in \S \ref{sec:function l}, and 
$\fA_\chi\inv(\fA_\ze(\tau))$ is a divisor in $\fA_{T}$. 
For $\tau=0$ this is the divisor 
$\Ker(\fA_\chi)$ whose ideal sheaf is  $\Th_T(\chi)$
and in general
$\fA_\chi\inv(\fA_\ze(\tau))$
is its torsor 
which we think of as a shift of $\Ker(\fA_\chi)=\fA_\chi\inv(0)$ by $\fA_\ze(\tau)\in\G$.



Now for any reductive group $G$ with a Cartan $T$ and Weyl group $W$,  
we decompose  $E$ according to $\DD$-action as
$E=\oplus_{\ze\in X^*(\DD)}\ (E_\ze\bten\ze\inv)$, for some $G$-modules $E_{\ze}$.  
Then $\Th_G(E_\ze)$ is the ideal sheaf of some divisor, denoted by $D(E_\ze)$,
in $\fA_G=\fA_T//W$. As $T$-representations, we have the decomposition $E_{\ze}|_{T}= \oplus_{\chi\in X^*(T)} [E_{\ze}: \chi] \chi$. 
Therefore, the divisor $D(E_\ze)$ is a sum over
$\chi\in X^*(T)$ of divisors $ [E_\ze:\chi]\cd\Ker(\fA_\chi)$.
Now,  for any $\tau\in\fA_\DD$,\  $\Th_\tG(E)_\tau$ is the ideal sheaf of the shifted divisors $D(E_\ze)+\fA_\ze (\tau)$  of $D(E_\ze)$.

\sss{
Quantum 
diagonals
}
\lab{m-shifted locality}
In our quiver setting, each $h\in H\sq H^*$
defines (via the Nakajima function $\bm$)
a character  $\mu_h\in X^*(\DD)$,
by which $\DD$ acts on the component $Rep_\bQ(V)_h$ of $Rep_\bQ(V)$.

For $_1V,\ _2V$ in $\VVI$
for each $\ii$ choose coordinates $_sx^i_p$ on $_sV^i$ hence a decomposition of $_sV^i$
into lines
$_sV^i_p$.
This gives Cartans $T_s\sub G_s=GL(_sV)$
with  a basis $_sx^i_p$ of $X^*(T_s)$.
Then on the line $\Hom(_1V^i_p,_2V^j_q)$,
the torus $\tT\dff T_1\tim T_2\tim\DD$
acts by $_2x^j_q\cd{}(_1x^i_p)^{-1}\cd\mu_h$, so its Thom divisor
is given by vanishing of 
$\fA_{_2x^j_q}+\fA_{\mu_h}-\fA_{_1x^i_p}$
in $\fA_{T_1 \tim T_2\tim\DD}$.
Therefore,
the Thom divisor of the $T_1\tim T_2\tim\DD$-module
$\Rep_\bQ(_1V,_2V)_h$
is
the shifted
diagonal
$$
\De^{v_1,v_2}_{h}(\tau)
\dff \De^{v_1,v_2}_{h',h''}+\ (0, \tau_h) \subset \fA_{T_1} \tim \fA_{T_2}
.$$
Here $\tau_h=\fA_{\mu_h}(t)$ depends on $h$, and $\De^{v_1,v_2}_{h',h''}\subset  \fA_{T_1} \tim \fA_{T_2}$ is the diagonal divisor defined by vanishing of $\prod_{p, q}(\fA_{_2x^j_q}-\fA_{_1x^i_p})$, 
and the shift $\De^{v_1,v_2}_{h',h''}+\ (0, \tau_h)$ means that for $\ _2V^j$
we use the
embedding of $\G=\fA_\Gm$ into $\fA_{GL(\ _2V^j)}$
via $\Gm=Z(GL(\ _2V^j))$ and the corresponding
addition action of $\G$ on $\fA_{GL(\ _2V^j)}$.

Consider the diagonal $\De^{v_1, v_2}_{i}$ of $\G^{v_1^i} \tim \G^{v_2^i}$ given by the vanishing of $\prod_{p, q}(\fA_{_1 x^i_p}-\fA_{_2 x^i_q})$. 
Let $\De^{v_1, v_2}_{i}(\tau)\dff   \De^{v_1, v_2}_{i} + (\tau, 0)$, where $\tau=\fA_{\omega}(t)$.  The character  $\omega\in X^*(\DD)$ is as before.
The shift $\De^{v_1,v_2}_{h',h''}+\ (\tau, 0)$ means that for $\ _1V^j$
we use the
embedding of $\G=\fA_\Gm$ into $\fA_{GL(\ _1V^i)}$
via $\Gm=Z(GL(\ _1V^i))$ and the corresponding
addition action of $\G$ on $\fA_{GL(\ _1V^i)}$, given by the vanishing of $\prod_{p, q}(\fA_{_1 x^i_p}-\fA_{_2 x^i_q}+\fA_{\omega})$.

We will say that for $\tau\in \fA_\DD$,
and $D_s=(D_s^i)_\ii\in \G^{|v_s|}$ for $s=1,2$;
the pair $(D_1,D_2)$ is {\em $(\bm,\tau)$-disjoint} if
$(D_1,D_2, \tau)$  and $(D_2,D_1,\tau)$ do not lie in any of the
shifted diagonals
$\De^{v_1,v_2}_{h}(\tau), \De^{v_1, v_2}_{i}, \De^{v_1, v_2}_{i}(\tau)$. Equivalently, for any $i\in I$,
the divisors $D_1^i\pm \tau$ and  $D_2^i$ are disjoint, $D_1^i$ and  $D_2^i$ are disjoint;
for each $h: h'\to h''$ in $\barr H$,
$D_2^{h''}\pm \tau_h$ and $D_1^{h'}$ are disjoint.

\subsection{
Quantum locality 
}
\lab{Quantum locality}
We will now consider locality in the setting of the (non-commutative) monoid
$\CC_{\G\tim I}$ freely generated by $\G\times I$.

\rc{\bbS}{\varpi}

Let $\CCI$ be the free monoid on $I$,
so elements are ordered sequences
$\gamma=i_1i_2\cdots i_N$ of elements in $I$.
The product of
$\gamma=i_1i_2\cdots i_N, \gamma'=j_1 j_2 \cdots j_{N'}$,
is the concatenation  
$\gamma+\gamma'=i_1i_2\cdots i_Nj_1 j_2 \cdots j_{N'}$. 

Let $\calC_{\G\times I}=\sq\ (\G\tim I)^n=
\sqcup_{\gamma\in \CCI} \G^{\gamma}$ be the
ind-scheme monoid freely generated by $\G\times I$,
with connected components labeled by $\CCI$. 
The natural projection, from the free monoid to the free commutative monoid
is denoted 
$\bbS:  \calC_{\G\times I} \to  \calH_{\G\times I}$.

We will use the notation $\LL=\LL_\DD(Q,A)$
both for the biextension line bundle defined on $\HGI^2\tim\fA_\DD$ in
\S\ref{sec:local and biextension} and also for its pull back to $\CC_{\G\tim I}^2\tim\fA_\DD$.
For
$\gamma', \gamma''\in \CCI$.
we denote by $\calL_{\gamma', \gamma''}$ its restriction to the 
 component $\G^{\gamma'}\times \G^{\gamma''} \times \fA_{\DD}$, 
 
\sss{$\fm$-locality
}
An {\em $\bm$-locality structure} on
a vector
bundle $K$ on $\calC_{\G\times I}\times \fA_{\DD}$ is a consistent system of isomorphisms
$$
(K_{\gamma_1,\tau}
\ \boxtimes \
K_{\gamma_2,\tau}) \otimes \calL_{\gamma_1, \gamma_2} \cong 
K_{\gamma_1+\gamma_2, \tau}\foor
\tau\in\fA_\DD
.$$

Any $\bm$-locality structure on $K$ implies an algebra structure on $K$ in 
the monoidal category $\big(\Coh(\calC_{\G\times I}\times \fA_{\DD}), \star \big)$
(by the biextension property of $\calL$).
In this way
an $\bm$-locality structure on $K$ is the same as a structure of a $\star$-algebra, whose multiplications are isomorphisms.\fttt{
Notice that this is stronger than the standard definition of
locality which only requires such isomorphism over the regular part of the configuration space
where
$\LL$ happens to trivialize by \ref{Quantization of convolution as shift of diagonals}.
}

\sex
The line bundle
$\sL=\sL_\DD(Q,A)$ constructed component-wise in \S\ref{sec:local and biextension} 
\textbf{Case 2}, is a line bundle over $\calC_{\G\times I} \tim\fA_{\DD}$ and it
has a natural $\bm$-locality structure.  So is the modification $\underline{\sL}_\DD(Q,A)$ as in \S~\ref{sss:modify}.
We will write the proof  only generically:

\slem
The line bundle $\sL$ on $\calC_{\G\times I} \times \fA_{\DD}$ defined in \S\ref{sec:local and biextension} 
\textbf{Case 2}
has the property that $(D_1,D_2)\in (\calC_{\G\times I})^2$ is {\em $(\bm,\tau)$-disjoint} for $\tau\in\fA_\DD$, then
there is a canonical identification of fibers
$$
\sL_{D_1+ D_2,\tau}
\ \cong\
\sL_{D_1,\tau}
\ \ten\
\sL_{D_2,\tau}
.$$

\pf
Let $V=V_1\pl V_2$ in $\VVI$
and 
$G_i=GL(V_i)$ and $G=GL(V)$.
Choose a  Cartan $T_i$ in $G_i$.
Then
$
\Rep_\bQ(V)\cong \Rep_\bQ(V_1)\oplus \Rep_\bQ(V_2)\oplus \Rep_\bQ(V_1, V_2)\oplus \Rep_\bQ(V_2, V_1)
$
gives
$$
\Th_\tG[\Rep_\bQ(V)]\ten \Th_\tG[\Rep_\bQ(V_1)]^{-1}\ten
\Th_\tG[\Rep_\bQ(V_2)]\inv
\ \cong\
\Th_\tG[\Rep_\bQ(V_1,V_2)]
\ten
\Th_\tG[\Rep_\bQ(V_2, V_1)]
.$$
Now the disjointness condition implies that the last two factors have canonical trivializations at $(D_1,D_2,\tau)$.
A similar statement holds for $\Theta_{\widetilde{L}}(\fg/\fp\otimes\omega)$, and $\Th_{\widetilde{L}}(\fg/\fp)\inv$, where 
$\fg/\fp=\oplus_{i=1}^2\fg_i/\fp_i$. 

Therefore, $\calL_{D_1, D_2}=\Th(d^*p)\otimes \Th(\tii q)$ has a canonical trivialization when $(D_1,D_2)$ is $(\bm,\tau)$-disjoint. 
The claim now follows from the identification 
$(\sL_{D_1,\tau}
\ \boxtimes \
\sL_{D_2,\tau}) \otimes \calL_{D_1, D_2} \cong 
\sL_{D_1+D_2, \tau}$. 
\qed

\srem The quantum local line bundle $\sL$ is in a sense a localization of
the quantum group $U^+_\DD(Q,A)$ to
the noncommutative configuration space $\CC_{\G\tim I}$.
By its definition the $\al$-weight space $U^+_\DD(Q,A)(\al)
$ is a sum of contributions from
all refinements $\ga\in\CCI$ of a given $\al\in\NI$.\fttt{
One formal way to say it is that 
$U^+_\DD(Q,A)_\al$  is the smallest subsheaf on $\G\syp{\al}$ such that it pull back to
all refinements $\CC_\ga$ contains $\sL_\ga$.
}

In the classical case $\DD=1$, for all $\ga$ above $\al$
$\sL_{\ga}$
are the same, so the sum
$U^+_\DD(Q,A)$ is the line bundle $\sL$. However, upon quantization
there is a genuine dependence on $\ga$ and one has to take the sum of all contributions in order to construct a  subalgebra.


\sex
In the case when $I$ is a point (the ``$sl_2$-case'') then $\CCI=\N[I]=\N$ hence 
 $\calC_{\G\tim I}$ is  the system $\sqcup_{n\in\N} \G^n$ of Cartesian powers of $\G$. 
Then $\underline{\sL}_{n}=\bbS_{n}^*(\underline{U}_{\DD}^+(Q, A)_{n})$.

\sss{
Some expectations
}
The above
construction of loop Grassmannians is of ``existential'' nature,
with hidden difficulties of explicit computations.
We hope to ameliorate this difficulty
by some equivalent descriptions.
Our construction is based on  ``abelianization''
(as we construct sections of  $\OO(1)$ on the loop Grassmannian
from the same objects for a Cartan subgroup)
and on locality (as we interpret
equations of the projective embedding of the Grassmannian
as locality conditions).

We would like to describe these equations
in more standard terms by constructing a central extension of
the quantum group $U_\DD(Q,A)$ and its action on sections of $\OO(1)$.
Here, the central extension should appear as one extends
the ``quantum local''  line bundle $\underline{\sL}_\DD(Q,A)$ 
from the analogue $\CC_{\G\tim I}$ of $\HGI$
to an analogue of $\GG(T)$.

One could also try to
construct the graded algebra of section of line bundles $\OO(m)$
by choosing the poset $P$ in $\GG^P_\DD(Q,A)$ to be $1<\cddd<m$.

\appendix

\se{
Loop Grassmannians with a condition
}
\lab{Drinfeld's theory of classifying pairs}

We recall a general technique
providing modular description of some parts of loop Grassmannian.
This allows us 
to finish
(in \ref{Proof of the proposition The T-fixed points})
the proof
of
proposition
\ref{The T-fixed points}
 on $T$-fixed points in closures of semi-infinite orbits.

\sss{
Moduli of finitely supported maps
}
\lab{Moduli of generically trivialized maps}

Here we recall some elements of Drinfeld's notion of loop Grassmannians with
a geometric  (``asymptotic'') condition, This material will be covered in more details elsewhere.
We will fix a smooth curve $C$.

We are interested  in various moduli  of $G$-torsors
over a curve $C$ that are local spaces over $C$.
As observed  by Beilinson and Drinfeld, 
the relevant spaces
$Y$ are
usually of the form 
$\MM_\YY(C)$, the moduli of
finitely supported\ie {\em generically trivialized} maps 
into  
some pointed stack $(\YY,\pt)$ built from $G$.
(We usually omit the point $\pt$ from notation.)

\sss{
The subfunctor $\uGG(G,Y)\sub\uGG(G)$ given by ``condition $Y$''
}
\lab{Subfunctor GG(G,Y) given by ``condition Y''}
Let $C$ be  a smooth connected curve with the generic point $\eta_C$. 
Let $G$ be an algebraic group and $(Y,y)$ a pointed scheme with a $G$-action on $Y$.
This gives a pointed stack $(\YY,*)$ with $\YY=G\bs Y$. 
Consider the moduli of maps of pairs $Map[(C,\eta_C),(\YY,*)]$.
Denote by $\uGG(G,Y)$ the space over $\HHC$
with the fiber at $D\in\HHC$ given by the maps $f\in 
Map[(C,\eta_C),(\YY,*)]
$ that are defined off
$D$. This is a {\em factorization space}
(\ref{Local spaces over a curve}).


If 
the orbit $Gy$ is open in $Y$ and its boundary $\del(Gy)$ is  a union of divisors $Y_i,\ \ii$,
to any $f\in Map[(C,\eta_C),(\YY,*)]$ one can associate
an $I$-colored finite subscheme $f\inv(\del*)\dff (f\inv G\bs Y_i)_\ii$.
Then we define 
$\uGG(G,Y;I)$ to be 
$Map[(C,\eta_C),(\YY,*)]$
considered as a space over $\HHCI$.
This is an $I$-colored local space (\ref{Local spaces over a curve}).

\sexs
(a) When $Y$ is a point, 
$\uGG(G,\pt)$ is the loop Grassmannian $\uGG(G)$.

(b) When $G=\Gm$ and $(Y,y)=(\Ao,1)$, $I$ is a point and
$\uGG(\Gm,\Ao,I)=Map[(C,\eta_C),(\Gm\bs\Ao,*)]$ is the 
space of effective divisors on $C$\ie the Hilbert scheme $\HHC$.

\slem
Let the scheme $Y$ be separated. 

(a) 
$\uGG(G,Y)$ is a subfunctor of $\uGG(G)$.
If $Y$ is also affine, then
$\uGG(G,Y)$ is closed in $\uGG(G)$.

(b) 
For a subgroup $K\sub G$ 
the intersection with 
$\uGG(K)\sub\uGG(G)$ reduces the condition $Y$ to the condition
$\barr{Ky}\sub Y$ :
$$\uGG(G,Y)\cap\uGG(K)
\ =\
\uGG(K,\barr{Ky})
.$$


\sss{
The closure of $S_0$ 
}
\lab{The closure of semi-infinite orbits}
\lab{The closure of S0}
It is well known that $G/N$ is quasi-affine\ie it is an open part of
its affinization
$\GNa$. We will consider it with the base point $y=eN$.

\spro Let $G$ be of   simply connected  type.

(a)
The scheme of $T$-fixed points $\uGG(G,\GNa)^T$
is $\HH_{d\tim I}$.


(b) 
The closure 
$\barr{S_0}$
is the reduced part $\uGG[G,\GNa]_{red}$
of the loop Grassmannian with the condition
$\GNa$.

\pf
(a)
The fixed points $\uGG(G)^T$ are known to be $\uGG(T)$ so
$
\uGG(G,\GNa)^T=
\uGG(G,\GNa)\cap\uGG(T)$.
This has been identified in the lemma 
\ref{Subfunctor GG(G,Y) given by ``condition Y''}.b
with
$\uGG(T,\barr{Ty})$ where 
$y$ is the base point $eN$ of $\GNa$.
So,  $Ty=B/N$.  When $G$ is simply connected
$\prod_\ii\ \ch\al_i:\GmI\con T\cong B/N$. This extends to an identification of  
the closure of $B/N$ in $\GNa$ (a $T$-variety) with $(\Ao)^I$
(a $\mathbb{G}_m^I$-variety).
Now, 
$\uGG(T,\barr{Ty})\cong\
\uGG(\Gm,\Ao)^I$ is identified with $\HH_{d\tim I}$ in the example 2
in \ref{Subfunctor GG(G,Y) given by ``condition Y''}.

(b)
Since $\GNa$ is affine, 
$\uGG(G,\GNa)$ is closed in $\uGG(G)$ 
(lemma 
\ref{Subfunctor GG(G,Y) given by ``condition Y''}.a).
Since $\uGG(G,\GNa)$  contains $\uGG(G,G/N)=\uGG(N)=S_0$,
its reduced part contains $\bSz$.
Since the stabilizer
of the base point of $\GNa)$ is $N$,
$\uGG(G,\GNa)\sub\uGG(G)$ is $\NK$-invariant.
Then the reduced part
$\uGG(G,\GNa)_{red}$ has a stratification by $\NK$-orbits $S_\la$ 
for $\la\in X_*(T)$ such that
$L_\la$ lies in $\uGG(G,\GNa)_{red}$. 

We have, according to 
Part (a), that  $\uGG(G,\GNa)^T$
 is $\HH_{d\tim I}$.  So,  $\uGG(G,\GNa)_{red}=\bSz$. 
\qed

\sss{
Proof of the proposition \ref{The T-fixed points}
}
\lab{Proof of the proposition The T-fixed points}
(a) According to proposition
\ref{The closure of S0}.a we have
$\bSz=\uGG(G,\GNa)_{red}$.
Therefore, $\bSz^T\sub \uGG(G,\GNa)^T$ which
is $\HHdI$ by proposition
\ref{The closure of S0}.b.

To see that 
$\HH_{d\tim I}\sub \uGG(G)$
lies in $\bSz$ we 
denote by $G_i\sub G$ the connected 3-dimensional subgroup
corresponding to $\ii$.
Then $\bSz$ contains the corresponding object $\bSz(G_i)$ for $G_i$,
and since we have already
checked the proposition
for $SL_2$ this is  $\AJ^{G_i}(\HH_d)$\ie $\A^GJ(\HH_{d\tim i})$.

It remains to prove that $\bSz^T\sub \uGG(T)$ is closed under the product
in $\uGG(T)$
(because $\AJ(\HHdI)$ is the product of all $\AJ(\HH_{d\tim i})$).
However,
the product in $\uGG(T)$ can be realized using fusion in $\GG(T)$.\fttt{
For $C=\A^1$ we have a canonical trivialization
of $\GG(G)\to
\HHC$ over $C=\HH^1_C$, as $\uGG(G)$.
Now, consider the pull-back $\GG_{C^2}(G)\dff C^2\tim_{C\hp{2}}\GG_\HCI(G)$
of 
the restriction of
$\GG_\HCI(G)$ to  $\HH^2=C\sp{2}$.
The locality identifies it over $C^2-\De_C$ with 
the constant bundle $\uGG(G)^2$.
By {\em fusion of $u,v\in\uGG(G)$} we mean the limit (when it exists) over the diagonal
of the constant section
$(u,v)$ which is
defined off the diagonal.
}
So, 
it suffices to notice that
$\bSz$ is the fiber at a point $a=0$ in a curve $C=\Ao$ 
of a factorization space 
$\barr{\GG(G,G/N)}$ which is defined as the 
closure of the factorization subspace
$\GG(N)\cong\GG(G,G/N)\sub \GG(G)$.

(b-c) The part (a) 
of the proposition
\ref{The T-fixed points} gives a factorization 
of $\bSz^T$ as a product $\HHdI\cong\ \prod_\ii\ (\HH_d)^I$
over contributions from all $\ii$. 
One therefore also has such factorization for
$\barr{S^-_{\al}}^T$ and obviously for  the connected components $\GG(T)_\be$.
This reduces
parts (b) and (c) of the proposition
to the $SL_2$ case. This case has already been checked
by explicit calculation
following
proposition \ref{The T-fixed points}.
\qed

\se{
Calculation of  Thom line bundles from \cite{YZell}
}\lab{Comparison with the Thom line bundles in cite{YZell}}
In \cite{YZell}, one uses a different convolution diagram. The only essential difference is the map $\iota$ described below. We check that it gives the same Thom line bundle 
as the calculation in \S\ref{The multiplication induced by the cotangent version of the correspondence} which used the dg cotangent correspondence. We will recall  without character formulas how computations of Thom line bundles were made in \cite{YZell}. 
For calculational reasons one uses an extra variety $X=G\tim_PY $ for $Y=\Rep_Q(\Gr_F(V))$
and then a nonlinear map $\io$ accounts for the difference between ambiental embeddings
$T^*_GX\sub T^*X$ and $T^*_LY\sub T^*Y$.

The notations are as in  \S\ref{The multiplication induced by the cotangent version of the correspondence}. 
Denote the elements of
$Y=\Rep(V_\bu)=\pl_{k=1}^m \Rep_Q V_k$ and $Y^*= \Rep_{Q^*} (V_\bu)$
by $y$ and $y^*$.
The moment map $\mu:T^*Y\to\bl^*\cong\fl$ is given by the projection
of the commutator to $\fl$
$$
\mu(y,y^*)=\ [y,y^*]_\fl
\dff\
\Big(\sum_{\{h\in H,\ h'=i\}}
y_hy^*_h-\sum_{\{h\in H: h''=i\}}
y_h^* y_h
\Big)_{i\in I}. 
$$

The story in \cite{YZell} is told in terms of
singular subvarieties $\mu^{-1}(0)\subseteq T^*Y$ 
(for a group $L$ acting on a smooth variety $Y$)
and the functoriality of cohomology is constructed in terms of
ambiental smooth varieties $T^*Y$.
The difference here is that
we derive  the cotangent correspondence mechanically from the original
correspondence. For instance this makes the associativity
of multiplication follow manifestly from associativity of the extension correspondence.


Let $W=G\times_P\RR(F)$ with projection to $X'=\RR(V)$. Let $Z:=T^*_W(X\times X')$. 
We have the following correspondence in \cite[Section 5.2]{YZ1}
\begin{equation}\label{equ:Lagr corresp} 
\xymatrix@R=1.5em{
G\times_{P} T^*Y
\ar@{^{(}->}[r]^(0.75){\iota}&T^*X& Z \ar[l]_{\phi}\ar[r]^{\psi} &T^*X'
}
\end{equation}
the maps are the natural ones, which we further describe below.

Let $U$ be the unipotent radical of $P$. 
Denote the Lie algebras by $\fp,\fl,\fu$.
Denote the natural projections 
by 
$\pi: P\to L$,
$\pi: \mathfrak{p}\to \mathfrak{l}$ and
$\pi':\mathfrak{p}\to \fu$.

For any associated  $G$-bundle $\EE=G\tim_PE$ we denote the fiber at the origin by
$\EE_\bz=E$.
Then $T^*X\cong G\tim_P(T^*X)_\bz$ and
the $L$-variety
$(T^*X)_\bz$
is (by \cite[Lemma 5.1 (a)]{YZ1})
\begin{align}
(T^*X)_\bz
\dff
\{ (c, y, y^*) \mid c\in \mathfrak{p}, (y, y^*)\in T^*Y,\  \text{such that}\
\mu(y, y^*)=\pi(c)\}. \label{eq:T*X}
\end{align}

\slem
\label{Map iota}
(a) We have an isomorphism of $L$-varieties
$\fu\tim T^*Y\cong (T^*X)_\bz$  over $G/P$
by 
$(u, y, y^*)\mapsto (u+\mu(y, y^*), y, y^*)$.

(b) This makes
$T^*X$
into a $G$-equivariant vector bundle over $G/P$, the sum of
$T^*(G/P)$ and $G\tim_PT^*Y$.

\begin{proof}
In (a) the inverse
map is $(c, y, y^*)\mapsto (\pi'(c), y, y^*)$.
In (b) we use
$T^*(G/P)\cong G\tim_P\fu$.
\end{proof}

The map $G\times T^*Y\to T^*X$ defined as $(g, y, y^*) \mapsto (g, \mu(y, y^*), y, y^*)$ induces  a well-defined map $\iota: G\times_P T^*Y \to T^*X$.
By  \cite[Lemma 5.1]{YZ1}, we have the isomorphism
\[
Z:=T^*_W(X\times X')\ \cong\
G\times_P\ \Rep_\bQ(F)
\]
with 
$
\psi(g, x, x^*) \mapsto\ ^g(x, x^*) 
$ for $g\in G$ and $(x,x^*)\in\Rep_\bQ(V)$. 
So, the map $\psi$ is a composition of the inclusion $\psi'$
of vector bundles over $G/P$ and
the conjugation action $\psi''$ (which acts by the same formula as $\psi$)
and the diagram is
$$
G\times_{P} T^*Y
\ \aa{\io}\inj\
T^*[G\tim_P\Rep(\Gr_F(V))]
\
\laa{\phi}
\ Z
\ \aa{\psi'}\sub \
 G\times_P \Rep_\bQ(V)
\raa{\psi''}
 \Rep_\bQ(V)
.$$

\slem
\label{lem:Thom1}
The Thom line bundles
$\Th_\tG(\psi')$, $\Th_\tG(\psi'')$ and $\Th_\tG(\io)$ 
are respectively the line bundles
$$
\Th_\tL[(F_\yy/F_{0})\Rep_\bQ(V)],\ \ 
\Th_\tL(\fg/\fp)\inv
\aand
\Th_\tL(\fp^\perp)=\Th_\tL(\fg/\fp\otimes\omega)
,$$ 
In particular, $\Th(d^*p)\otimes \Th(\tii q)\cong \Theta_L(\iota)\otimes \Theta_L(\psi)$.

\pf
If $\SS$ is one of the first four spaces in the diagram, then
$\fA_\tG(\SS)=\fA_\tL$ since
$\SS=G\tim_P\SS_\bz$ 
for the fiber  $\SS_\bz$ which is an affine space.
In particular, 
for a map $\eta\in\{\io,\psi',\psi''\}$,
the line bundle $\Th_\tG(\eta)$ on $\fA_\tL$ is $\Th_\tL(T(\eta)_\bz)$.

(1)
Vector bundle $T(\io)$
is the normal bundle $N(\iota)$.
According to the   
lemma \ref{Map iota} it is isomorphic to $G\tim_P-$ of the $\tL$-module
$(\fg/\fp)^*\otimes\omega=\fp^\perp\otimes\omega$.

(2)
Similarly, $T(\psi')$ is the normal bundle $N(\psi')$
and the fiber
$T(\psi')_\bz$ is 
$(F/F_0)\Rep_Q(V)$.

(3)
The equality
$\Th_\tG(\psi'')
=\
\Th_\tP[\fg/\fp]\inv
$ is clear.
\qed

\scor

(a)
$
\Th_L[F_\yy/F_0(\Rep_\bQ(V)]
\ =\ 
\Th_L[\Rep_Q(V)- Gr_0^F(\Rep_{Q}(V)]
$.

(b) Consider the case when $Q$ has no loop edges
and the
filtration type $\bv$ is a flag\ie $\bv_k\in I$ for all $k$.
Then $Gr_0\Rep(V)=0$ and
$\Th_L(\psi')=\Th_G(Rep_QV)$.

\pf
(a) 
A filtration on $V$ induces
a family of filtrations, 
compatible  with  
the decomposition
$
\Rep_\bQ(V)
=\
\Rep_Q(V)
\pl
\Rep_{Q^*}(V)
$
and with the $L$-equivariant identification $\Rep_{Q^*}(V)\cong[\Rep_{Q}(V)]^*$.
Therefore, the claim follows from
$$
\fra{F_\yy}{F_0}[(\Rep_{Q}(V)^*]
=\
[\fra{F_{-1}}{F_{-\yy}}(\Rep_{Q}(V)]^*
$$
and the invariance of  Thom line bundles  under duality of
vector bundles.

(b) follows since 
$Gr_0^F\Rep_Q(V)=0$
under the assumption on $Q$. The reason is that $Gr_0^F\Rep_Q(V)=\pl\ \Rep_Q(Gr_pV)$
and all $Gr_p(V)\in\VVI$ are
one-dimensional.
\qed

\end{document}